%%%%%%%%%%%%%%%%%%
%
% T1 paper         
% dec0507.tex RB 
%
%%%%%%%%%%%%%%%%%%%

\magnification=\magstep1
\def\refer#1{ \medskip \parindent=34pt \hang \noindent \rlap{[#1]\hfil} \hskip 30pt }
\baselineskip=14pt

\font\tenmsb=msbm10
\font\sevenmsb=msbm7
\font\fivemsb=msbm5
\newfam\msbfam
\textfont\msbfam=\tenmsb
\scriptfont\msbfam=\sevenmsb
\scriptscriptfont\msbfam=\fivemsb
\def\Bbb#1{{\fam\msbfam\relax#1}}
\def\lam{{\lambda}}

\def\al{{\alpha}}

\def\proof{{\medskip\noindent {\bf Proof. }}}
\def\longproof#1{{\medskip\noindent {\bf Proof #1.}}}
\def\qed{{\hfill $\square$ \bigskip}}
\def\subsec#1{{\bigskip\noindent \bf{#1}.}}

\def\cite#1{{[#1]}}

\def\vareps {\varepsilon}
\def\eps{\varepsilon}

\def\q {\quad} \def\qq {\qquad}
\def\frac#1#2{{#1\over #2}}
\def\del{{\partial}}

\def\ol{\overline}

\def\wh{\widehat}

\def\ni{\noindent }
\def\ms{\medskip}
\def\bs{\bigskip}
\def\cl#1{\centerline{#1}}

\parindent=30pt

\def\Osc{\mathop {\rm Osc\, }}

\def\sgn{{\mathop {{\rm sgn\, }}}}

\def\square{{\vcenter{\vbox{\hrule height.3pt
\hbox{\vrule width.3pt height5pt \kern5pt
 \vrule width.3pt}
\hrule height.3pt}}}}

\def\tfrac#1#2{{\textstyle {#1\over #2}}}

\def\choose#1#2{{\left(\matrix{#1\cr#2\cr}\right)}}
\def\tlint{{- \kern-0.85em \int \kern-0.2em}} 
\def\dlint{{- \kern-1.05em \int \kern-0.4em}}

\def\y0{\langle y_0,\be_i\rangle}
\def\R{{\Bbb R}}
\def\E{{{\Bbb E}\,}}
\def\P{{\Bbb P}}
\def\Q{{\Bbb Q}}
\def\Z{{\Bbb Z}}
\def\F{{\cal F}}

\def\N{{\Bbb N}}

\def\sA {{\cal A}}  
\def\sD {{\cal D}} \def\sE {{\cal E}} \def\sF {{\cal F}}
  
  \def\sL {{\cal L}}

\def\tD {\tilde D}

\def\Xm{X^{(m+1)}}
\def\be{{\epsilon}}

\font\ften=cmr8

\cl{\bf Degenerate Stochastic Differential Equations arising}
\ms
\cl{\bf 
 from Catalytic Branching Networks}

\bs

\centerline{Richard F. Bass${}^{1}$,\q Edwin A. Perkins${}^{2}$}, 
\footnote{$\empty$}{\ften 1. Research supported in part by an NSF grant DMS0601783.}
\footnote{$\empty$}{\ften 2. Research supported in part by an NSERC Research Grant.}
\footnote{$\empty$}{\ften {Keywords:} stochastic differential equations, perturbations, resolvents, Cotlar's lemma}
\footnote{$\empty$}{\ften \q \q \q \q \q catalytic branching, martingale problem, degenerate diffusions}
\footnote{$\empty$}{\ften Classification:  Primary 60H10. Secondary 35R15, 60H30}
\footnote{$\empty$}{\ften Running Head: SDEs for Branching Networks} 

\vskip0.5truein

{ \centerline{\bf Abstract}} 

We establish existence and uniqueness for the martingale problem associated with a system of degenerate SDE's representing a catalytic branching network.  For example, in the hypercyclic case:
$$dX_{t}^{(i)}  
=b_i(X_t)dt+\sqrt{2\gamma_{i}(X_{t})
X_{t}^{(i+1)}X_{t}^{(i)}}dB_{t}^{i},\  X_t^{(i)}\ge 0,\ i=1,\dots, d,$$
where $X^{(d+1)}\equiv X^{(1)}$,
existence and uniqueness is proved when $\gamma$ and $b$ are continuous 
on the positive orthant, $\gamma$ is strictly positive, and  
$b_i>0$ on $\{x_i=0\}$.  The special case $d=2$, $b_i=\theta_i-x_i$ is required in work of [DGHSS] on mean fields limits of block averages for $2$-type branching
 models on a hierarchical group.  
The proofs make use of some new methods, including Cotlar's lemma to establish asymptotic orthogonality of the derivatives of an associated semigroup at different times, 
and a refined  integration by parts technique from [DP1].   As a by-product of the proof we obtain the strong Feller property of the associated resolvent. 

\vfil\eject

\subsec{1. Introduction}
In this paper we establish well-posedness of the martingale problem for certain 
degenerate second order elliptic operators.
The class of operators we consider arises from models of catalytic branching networks including 
catalytic branching, mutually catalytic branching and hypercyclic 
catalytic branching systems (see [DF] for a survey of these systems). 
For example, the
hypercyclic catalytic branching model  is a diffusion on
$\R^{d}_{+},\;d\ge 2$, solving the following
system of stochastic differential equations:
$$dX_{t}^{(i)}  
=(\theta_i-X_{t}^{(i)})dt+\sqrt{2\gamma_{i}(X_{t})
X_{t}^{(i+1)}X_{t}^{(i)}}dB_{t}^{i},\  i=1,\dots, d. \eqno(1.1)
$$
Here $X(t)=(X^{(1)}_t,\dots, X^{(d)}_t)$, addition of the
superscripts is done cyclically so that $X^{(d+1)}_t=X^{(1)}_t$,
$\theta_i>0$, and $\gamma_i>0$. 

Uniqueness results of this type are proved in [DP1] under
H\"older continuity hypotheses on the coefficients.  Our main result here is to
show the uniqueness continues to hold if this is weakened to continuity.  One motivation for 
this problem is that for $d=2$, (1.1) arises in [DGHSS] as 
the mean field limit of the block averages of a system of SDE's on 
a hierarchical group.  The system of SDEs models two types of individuals 
interacting through migration between sites and at each site through 
interactive branching, depending on the masses of the types at 
that particular site.  The branching coefficients $\gamma_i$ of 
the resulting equation for the block averages involves averaging 
the original branching coefficients at a large time (reflecting 
the slower time scale of the block averages) and so is given in terms of 
the equilibrium distribution of 
the original equation.  The authors of [DGHSS]  introduce a renormalization map which gives 
the branching coefficients $\gamma_i$ of the block averages in terms of 
the previous SDE.  They wish to iterate 
this map to study higher order block averages.  Continuity is preserved by 
this map on the interior of $\R_+^d$,
and is conjectured to be preserved at the boundary (see
Conjecture 2.7 of [DGHSS]). It is not known whether H\"older
continuity is preserved (in the interior and on the boundary),
which is why the results of [DP1] are not strong enough to
carry out this program.  The weakened hypotheses also leads to some new methods.

The proofs in this paper are  substantially simpler in the two-dimensional setting required 
for [DGHSS] (see Section 8 below) but as higher dimensional analogues 
of their results are among the ``future challenges" stated there, we thought 
the higher-dimensional results worth pursuing.  

Further motivation for the study of such branching catalytic networks comes from [ES] where a corresponding system of 
ODEs was proposed as a macromolecular precursor to early forms of life. 
There also have been a number of mathematical works on mutually catalytic branching ((1.1) with 
$d=2$ and $\gamma_i$ constant) in spatial settings where a special duality argument ([M, [DP2])
allows a more detailed analysis, and even in spatial analogues of (1.1) for general $d$, but now with much more restricted results due in part to the lack of any uniqueness result ([DFX], [FX]).  
See the introduction of [DP1] for more background material on the model.

Earlier work in [ABBP] and [BP] show uniqueness in the martingale problem for the operator  
${\sA}^{(b,\gamma)}$ on $C^{2}(\R_{+}^{d})$ defined by
$$
{\sA}^{(b,\gamma)}f(x)=\sum_{i=1}^{d}\left(
b_i(x)\frac{\partial f}{\partial
x_{i}}+\gamma_{i}(x)x_{i}\frac{\partial^{2}f}{\partial
x_{i}^{2}}\right)  ,\ \ \ x\in\R_{+}^{d}.
$$
Here $b_i, \gamma_i\;i=1,\dots,d$ are continuous functions on
$\R_+^d$, with $b_i(x)\ge0$ if $x_i=0$, and satisfying some additional regularity
or non-degeneracy condition.  If $b_i(x)=\sum_j x_jq_{ji}$ for some $d\times d$ $Q$-matrix $(q_{ji})$, then such diffusions arise as limit points of rescaled systems of critical branching Markov chains in which $(q_{ji})$ governs the spatial motions of particles and $\gamma_i(x)$ is the branching rate at site $i$ in population $x=(x_1,\dots,x_d)$.  The methods of these papers do not apply to systems such as (1.1) because now the branching rates $\gamma_i$ may be zero.  Although we will still proceed using a Stroock-Varadhan perturbation approach, the process from which we are perturbing will be more involved than the independent squared Bessel process considered in the above references. 

We  will formulate our results in terms of catalytic 
branching networks in which the catalytic reactions are given by 
a finite directed graph $(V, \sE)$ with vertex set $V=\{1,\dots,d\}$ 
and edge set $\sE=\{e_1,\dots,e_k\}$.  
This will include (1.1) and all of the two-dimensional systems arising in [DGHSS].  
As in [DP1] we assume throughout:
\medskip
\noindent {\bf Hypothesis 1.1.} $(i,i)\notin \sE$ for all $i\in V$ and
each vertex is the second element
of at most one edge.\qed
\medskip
The restrictive second part of this hypothesis has been removed by Kliem [K] in the H\"older continuous setting of [DP1].  It is of course no restriction if $|V|=2$ (as in [DGHSS]), and holds in the cyclic setting of (1.1).  

Vertices denote types and an edge $(i,j)\in\sE$ indicates that type $i$ catalyzes the
type $j$ branching. Let $C$ denote the set of
vertices (catalysts) which appear as the first  element of an edge and $R$ denote the set of vertices
that appear as the second element (reactants).  Let $c:R\to C$ be such that for $j\in
R$, $c_j$ denotes the unique $i\in C$ such that $(i,j)\in
\sE$, and for $i\in C$, let $R_i=\{j:(i,j)\in\sE\}$.

Here are the hypotheses on the coefficients:
\ms\ni
{\bf Hypothesis 1.2.} For $i\in V$,
$$
\gamma_{i}:R_{+}^{d}\rightarrow (0,\infty),\quad
b_i:\R_+^d \rightarrow \R,
$$
are continuous
such that $|b_i(x)|\le c(1+|x|)$ on
$\R_+^d$, and $b_i(x)>0$ if $x_i=0$. \qed
\ni The positivity condition on $b_i|_{x_i=0}$ is needed
to ensure the solutions remain in the first orthant. 

\ms

If $D\subset \R^d$, $C_b^2(D)$ denotes the space of twice continuously differentiable bounded functions on $D$ whose first and second order partial derivatives are also bounded.
For $f\in C^2_b(\R_+^d)$, and with the above interpretations, the generators
we study are 
$$\sA f(x)=\sA^{(b,\gamma)} f(x)=\sum_{j\in
R}\gamma_j(x)x_{c_j}x_{j}f_{jj}(x)+\sum_{j\not\in
R}\gamma_j(x)x_{j}f_{jj}(x) +\sum_{j\in V}b_j(x)f_j(x).$$
(Here and elsewhere we use $f_i$ and $f_{ij}$ for the
first and second  partial derivatives of $f$.)

\ms\ni{\bf Definition 1.2.5.}
Let $\Omega=C(\R_+, \R_+^d)$, the continuous functions from $\R_+$
to $\R_+^d$. Let $X_t(\omega)=\omega(t)$ for
$\omega\in \Omega$, and let $(\F_t)$ be the canonical right continuous
filtration generated by $X$.
 If $\nu$ is a probability on $\R^{d}_{+}$, a
probability $\P$ on $\Omega$ solves the
martingale problem $MP({\sA},\nu)$  if under $\P$, the law
of $X_0$ is $\nu$ and for all $f\in
C_b^2(\R_+^d)$,
$$M_f(t)=f( X_t)-f(X_0)-\int_0^t
{\sA}f(X_s)\,ds$$ is a local martingale under $\P$.

A natural  state space for our martingale problem is
$$S=\Big\{x\in \R^{d}_{+}: \prod_{(i,j)\in\sE}(x_{i}+x_j)>0\Big\}.$$
The following result is Lemma 5 of [DP1] -- the H\"older continuity assumed there plays no role in the proof. 

\proclaim Lemma 1.3. If $\P$ is a solution of $MP({\sA},\nu)$,
where $\nu$ is a probability on $\R^{d}_{+}$, then $X_t\in
S$ for all $t>0$ $\P$-a.s.

\ni Here is our main result.

\proclaim Theorem 1.4.
Assume Hypotheses 1.1 and 1.2 hold. Then
for any probability $\nu$ on $S$, there is exactly one solution
to MP$({\sA},\nu).$

The cases required in Theorem~2.2 of [DGHSS] are the three possible directed graphs 
for $V=\{1,2\}$:

\item{(i)} $\sE=\emptyset$;
\item{(ii)} $\sE=\{ (2,1)\}$ or $\sE=\{(1,2)\}$;
\item{(iii)} $\sE=\{(1,2),(2,1)\}$.

The state space here is  $S=\R^2-\{(0,0)\}$.  In addition, [DGHSS] takes $b_i(x)=\theta_i-x_i$ for $\theta_i\ge 0$.  
As discussed in Remark~1 of [DGHSS],  weak uniqueness is trivial 
if either $\theta_i$ is $0$, as that coordinate becomes absorbed at $0$, 
so we may assume $\theta_i>0$.  
In this case  Hypotheses~1.1 and 1.2 hold,  and Theorem 2.2, stated in [DGHSS] 
(the present paper is cited for a proof),  
is immediate from Theorem~1.4 above. 
See Section~8 below for further discussion about our proof and how it simplifies in this two-dimensional setting.
In fact, in Case (i) the result holds for any $\nu$ on all of $\R_+^2$ (as 
again noted in Theorem~2.2 of [DGHSS]) by Theorem A of [BP].

Although it is not required in [DGHSS], 
it is of course natural to ask about uniqueness in cases (ii) and (iii) if $\nu=\delta_{(0,0)}$.  
We do have some partial results when the process starts at the corner $(0,0)$,
but the regularity hypotheses are stronger than in Hypothesis 1.2 and the
techniques are quite different than those used in this paper, so we do not
pursue this here.

Our proof of Theorem 1.4 actually proves a stronger result. We do not
require that the $\gamma_i$ be continuous, but only that their oscillation
not be too large. More precisely, we prove that there exists $\eps_0>0$
such that if (1.2) below holds, then there is exactly one solution of
$MP(\sA,\nu)$. The condition needed is
\medskip 
For each $i=1, \ldots, d$ and each $x\in \R_+^d$ there exists a neighborhood
$N_x$ such that 
$$\hbox{$\Osc_{ N_x}$} \gamma_i<\eps_0,\eqno (1.2)$$ 
where $ \Osc_A f=\sup_A f-\inf_A f$.
\medskip

Our proof of Theorem 1.4 is an $L^2$ perturbation argument,
and some of our argument follows along the lines of [ABBP]. The operators
from which we are perturbing are now different and 
the method of [ABBP] for obtaining
$L^2$ estimates no longer applies.  This leads to some new methodologies.  

The analytic tool we use is Cotlar's lemma, Lemma 2.13, which is
also at the heart of the famous T1 theorem of harmonic analysis. For
a simple application of how Cotlar's lemma can be used, see [Fe], pp.~103--104.

We consider certain operators $T_t$ (defined in (2.18) below) and show
that $$\|T_t\|_2\leq c/t. \eqno (1.3)$$
We require $L^2$ bounds on $\int_0^\infty e^{-\lam t} T_t\, dt$, and (1.3)
is not sufficient to give these. This is where Cotlar's lemma comes in:
we prove $L^2$ bounds on $T_tT_s^*$ and $T_t^*T_s$, and these together
with Cotlar's lemma yield the desired bounds on $\int_0^t e^{-\lam t} T_t\, dt$. The use of Cotlar's lemma to operators arising from a decomposition of the time axis is perhaps noteworthy. In all
other applications of Cotlar's lemma that we are aware of, the corresponding
operators arise from a decomposition of the space variable.
The $L^2$ bounds on $T_tT^*_s$ and $T_t^*T_s$ are the hardest and 
lengthiest parts of the paper. At the heart of these bounds is an integration by parts formula which refines a result used in [DP1] (see the proof of Proposition 17 there) and is discussed in the next section.

In Section 3 we give a proof of Theorem 1.4. The proofs of all
the hard steps, are, however, deferred to later sections.
A brief outline of the rest of paper is given at the end of Section~2.

\ms 

\ni{\bf Acknowledgment.} We would like to thank Frank den Hollander and Rongfeng Sun for helpful
conversations on their related work.

\subsec{2. Structure of the proof}

We first reduce Theorem~1.4 to a local uniqueness result (Theorem~2.1 below).  
Many details are suppressed as this argument is a minor modification of the 
proof of Theorem 4 in [DP1].
  By the localization argument in Section 6.6 of [SV] it suffices 
to fix $x^0\in S$ and show that for some $r_0=r_0(x^0)>0$, there 
are coefficients which agree with $\gamma_i$, $b_i$ on $B(x^0,r_0)$, 
the open ball of radius $r_0$ centered at $x^0$, and for which 
the associated martingale problem has a unique solution for all 
initial distributions.
  Following [DP1], let $Z=\{i\in V:x^0_i=0\}$, $N_1=\cup_{i\in Z\cap C}R_i$, $\bar N_1=N_1\cup(Z\cap C)$, and $N_2=V-\bar N_1$.
  Note that $N_1\cap Z=\emptyset$ because $x^0\in S$.
  Define
$$\tilde \gamma_j(x)=\cases{x_j\gamma_j(x)& if $j\in N_1$;\cr
x_{c_j}\gamma_j(x)& if $j\in (Z\cap C)\cup(N_2\cap R)$;\cr
\gamma_j(x)& if $j\in N_2\cap R^c$,}$$
and note that $\gamma_j^0\equiv \tilde \gamma_j(x^0)>0$ for all $j$ because $x^0\in S$.
  We may now write
$$\eqalign{\sA^{b,\gamma}f(x)=&\sum_{i\in Z\cap C}\Bigl[\sum_{j\in R_i}\tilde\gamma_j(x)x_if_{jj}(x)\Bigr]+\tilde \gamma_i(x)x_i f_{ii}(x)\cr
&+\sum_{j\in N_2}\tilde\gamma_j(x)x_jf_{jj}(x)+\sum_{j\in V}b_j(x)f_j(x).}$$
Let $\delta=\delta(x^0)=\min_{i\in Z}b_i(x^0)>0$, and define
$$\tilde b_j(x)=\cases{b_j(x) & if $j\in N_1$;\cr
b_j(x)\vee {\delta\over 2} & if $j\notin N_1$,}$$
and let $b_j^0=\tilde b_j(x^0)$, so that $b_j^0>0$ for $j\notin N_1$. Although $b_j(x^0)\le 0$ is possible for $j\in N_2\cap Z^c$ (and so $\tilde b_j$ may differ from $b_j$ here), a simple Girsanov argument will allow us to assume that $b_j(x^0)\ge \delta$ for $j\in N_2\cap Z^c$ (see the proof below) and so $\tilde b_j=b_j$ near $x^0$.  With this reduction we see
that by Hypothesis 1.2 and the choice of $\delta$, $\tilde b_j(x)=b_j(x)$ for $x$ near $x^0$.  By changing $\tilde b$ and $\tilde \gamma$ outside a small ball centered at $x^0$ we may assume $\tilde \gamma_j>0$ for all $j$, $\tilde b_j>0$ for $j\notin N_1$, $\tilde \gamma_j,\tilde b_j$ are all bounded continuous and constant outside a compact set,  and 
$$\eps_0\equiv\sum_{j=1}^d \big(\Vert\tilde\gamma-\gamma_j^0\Vert_\infty+\Vert \tilde b_j-b_j^0\Vert_\infty\big)\eqno(2.1)$$
is small.  For these modified coefficients introduce
$$\eqalignno{\tilde\sA f(x)=&\sum_{i\in Z\cap C}\Bigl[\sum_{j\in R_i}\tilde\gamma_j(x)x_if_{jj}(x)\Bigr]+\tilde \gamma_i(x)x_i f_{ii}(x)\cr
&+\sum_{j\in N_2}\tilde\gamma_j(x)x_jf_{jj}(x)+\sum_{j\in V}\tilde b_j(x)f_j(x),&(2.2)}$$
and also define a constant coefficient operator
$$\eqalignno {\sA^0f(x)=&\sum_{i\in Z\cap C}\Bigl[\sum_{j\in R_i}\gamma_j^0x_if_{jj}(x)+b_j^0f_j(x)\Bigr]+\gamma_i^0x_i f_{ii}(x)+b_i^0f_i(x)\cr
&+\sum_{j\in N_2}\gamma^0_jx_jf_{jj}(x)+b_j^0f_j(x)&(2.3)\cr
\equiv&\sum_{i\in Z\cap C}\sA^1_i+\sum_{j\in N_2}\sA^2_j.}$$
As $b_j^0\le 0$ and $\tilde {b_j}|_{x_j=0}\le 0$ is possible for $j\in N_1$ (recall we have modified $\tilde b_j$), the natural state space for the above generators is the larger
$$S^0\equiv S(x^0)=\{x\in \R^d:x_j\ge 0\hbox{ for all }j\notin N_1\}.$$
When modifying $\tilde\gamma_j$ and $\tilde b_j$ it is easy to extend them to this larger space, still ensuring all of the above properties of $\tilde b_j$ and $\tilde \gamma_j$.  If $\nu_0$ is a probability on $S^0$, a solution to the martingale problem $MP(\tilde\sA,\nu_0)$ is a probability $\P$ on $C(\R_+,S^0)$ satisfying the obvious analogue of the definition given for $MP(\sA, \nu)$.
As we have $\sA f(x)=\tilde \sA f(x)$ for $x$ near $x^0$, the localization in [SV] shows that Theorem~1.4 follows from:

\proclaim Theorem 2.1. Assume $\tilde\gamma_j:S(x^0)\to (0,\infty)$,  $\tilde b_j:S(x^0)\to\R$ are 
bounded continuous and constant outside a compact set with $\tilde b_j>0$ for $j\notin N_1$.
For $j\le d$, let $\gamma^0_j>0$, $b^0_j\in\R$, $b_j^0>0$ if $j\notin N_1$, and 
$$M_0=\max_{j\le d}(\gamma_j^0,(\gamma_j^0)^{-1},|b_j^0|)\vee \max_{j\notin  N_1} (b_j^0)^{-1}.\eqno(2.4)$$
There is an $\eps_1(M_0)>0$ so that if $\eps_0\le\eps_1(M_0)$, then for any probability $\nu$ on $S(x^0)$, there is a unique solution to MP$(\tilde\sA,\nu)$.

\ni{\bf Proof of reduction of Theorem 1.4 to Theorem 2.1.}  This proceeds as in the proof of Theorem 4 in [DP1].  The only change is that in Theorem~2.1 we are now assuming $\tilde b_j>0$ and $b_j^0>0$ for all $j\notin N_1$, not just $\tilde b_j\geq 0$ on $\{x_j=0\}$ for $j\notin N_1$ and $b^0_j>0$ for $j\in Z\cap(R\cup C)$ with $b^0_j\ge 0$ for other values of $j\notin N_1$.  If $b_j(x^0) >0$ for all $j\in N_2$, then the proof of Theorem 4 in [DP1] in Case 1 applies without change, and so we need only modify the argument in  Case 2 of the proof of Theorem 4 in [DP1] so that it applies if $b_j(x^0)\le 0$ for some $j\in N_2$.  This means $x_j^0>0$ by our 
(stronger) Hypothesis 1.2 and the Girsanov argument given there now allows us to locally modify $b_j$ so that $b_j(x^0)>0$.  The rest of the argument now goes through as before.  \qed

Turning to the proof of Theorem~2.1, existence is proved as  in Theorem 1.1 of [ABBP]--instead
of the comparison argument given there,
one can use Tanaka's formula and (2.4) to see that solutions must remain in $S(x^0)$.  

We focus on uniqueness from here on. 

The operator $\sA^2_j$ is the generator of a Feller branching diffusion with immigration. We denote its semigroup by $Q_t^j$.  It will be easy to give an explicit representation for the semigroup $P^i_t$ associated with $\sA^1_i$ (see (3.2) below).  An elementary argument shows that the martingale problem associated with $\sA^0$ is well-posed and the associated diffusion has semigroup 
$$P_t=\prod_{i\in Z\cap C}P_t^i\prod_{j\in N_2} Q_t^j,\eqno(2.5)$$ 
and resolvent $R_\lam=\int e^{-\lam t}P_t\,dt$.  Define a reference measure $\mu$ on $S^0$ by 
$$\mu(dx)=\prod_{i\in Z\cap C}\Bigl[\prod_{j\in R_i} dx_j\Bigr]\,x_i^{b^0_i/\gamma_i^0-1}dx_i\times \prod_{j\in N_2}x_j^{b_j^0/\gamma^0_j-1}dx_j=\prod_{i\in Z\cap C}\mu_i\prod_{j\in N_2}\mu_j.$$ 
The norm on $L^2\equiv L^2(S^0,\mu)$ is denoted by $\Vert\cdot\Vert_2$.

The key analytic bound we will need to carry out the Stroock-Varadhan perturbation analysis is the following:

\proclaim Proposition 2.2. There is a dense subspace $\sD_0\subset L^2$ and a $K(M_0)>0$ such that $R_\lam:\sD_0\to C_b^2(S^0)$ for all $\lam>0$ and
$$\eqalignno{\Big[\sum_{i\in Z\cap C}&\Big[\sum_{j\in R_i}\Vert x_i(R_\lam f)_{jj}\Vert_2\Big]+ \Vert x_i(R_\lam f)_{ii}\Vert_2\Big]+\Big[\sum_{j\in N_2}\Vert x_j(R_\lam f)_{jj}\Vert_2\Big]+\Big[\sum_{j\in V}\Vert (R_\lam f)_j\Vert_2\Big]\cr
&\le K\Vert f\Vert_2\quad\hbox{for all }f\in \sD_0\hbox{ and }\lam\ge 1.&(2.6)}$$

Here are the other two ingredients needed to complete the proof of Theorem ~2.1.

\proclaim Proposition 2.3.  Let $\P$ be a solution of $MP(\tilde\sA,\nu)$ where $d\nu=\rho\, d\mu$ for some $\rho\in L^2$ with compact support and set $S_\lam f=\E_{\P}\Bigl(\int_0^\infty e^{-\lam t}f(X_t)\,dt\Bigr).$
If 
$$\vareps_0\le (2K(M_0))^{-1}\wedge (48 d M_0^5)^{-1},\eqno(2.7)$$
then for all $\lam\ge 1$,
$$\Vert S_\lam\Vert:=\sup\{|S_\lam f|:\Vert f\Vert_2\le 1\}\le {2\Vert\rho\Vert_2\over \lam}<\infty.$$

\proclaim Proposition 2.4.  Assume $\{\P^x:x\in S^0\}$ is a collection of probabilities on $C(\R_+,S^0)$ such that:\hfil\break
\ni(i) For each $x\in S^0$, $\P^x$ is a solution of $MP(\tilde\sA,\delta_x)$.\hfil\break
\ni(ii) $(\P^x,X_t)$ is a Borel strong Markov process.\hfil\break
Then for any bounded measurable function $f$ on $S^0$ and any $\lam>0$, 
$$S_\lam f(x)=\E^x\Bigl(\int_0^\infty e^{-\lam t}f(X_t)\,dt\Bigr)$$
is a continuous function in $x\in S^0$.
\ms

\ni{\bf Remark 2.5.} Our proof of Proposition 2.4 will also show the strong Feller property of the resolvent for solutions to
the original MP$(\sA,\nu)$ in Theorem~1.4--see Remark 6.2. 
\ms

Assuming Propositions 2.2--2.4 the proof of Theorem~2.1 is then standard and quite similar to the proof of Proposition 2.1 in Section~7 of [ABBP].  Unlike [ABBP] the state space here is not compact, so we present the proof for completeness. 

\ni{\bf Proof of Theorem 2.1.} Let $\Q_k,k=1,2$, be solutions to $MP(\sA,\nu)$ where $\nu$ is as in Proposition~2.3 and define $S^k_\lam f=\E_k\Bigl(\int e^{-\lam t}f(X_t)\,dt\Bigr)$, where $\E_k$ denotes expectation with respect to $\Q_k$.  Let $f\in C_b^2(S^0)$.  
The martingale problem shows there is a local martingale $M^f$ satisfying
$$f(X_t)=f(X_0)+M^f(t)+\int_0^t \tilde\sA f(X_s)\,ds.\eqno(2.8)$$
Note that for $t>0$,
$$\eqalign{\E_k(\sup_{s\le t}|M^f_s|)&\le 2\Vert f\Vert_\infty+\int_0^t\E_k(|\tilde\sA f(S_s)|)\,ds\cr
&\le 2\Vert f\Vert_\infty+c\int_0^t\E_k\Big(\sum_{j\notin N_1}X^j_s+1\Big)\,ds<\infty,}$$
where the finiteness follows by considering the associated SDE for $X^j$ and using the boundedness of $\tilde b_j$.  This shows that $M^f$ is a martingale under $\Q_k$.  Let $g\in \sD_0$.  Multiply (2.8) by $\lam e^{-\lam t}$ integrate over $t$, take expectations (just as in (7.3) of [ABBP]), and 
set $f=R_\lam g\in C^2_b$ to derive
$$S_\lam^k g=\int R_\lam g\, d\nu+S^k_\lam((\tilde\sA-\sA^0)R_\lam g).$$
Taking the difference of this equation when $k=1,2$, we obtain
$$|(S^1_\lam-S^2_\lam)g|\le \Vert S_1^\lam-S^2_\lam\Vert\ \Vert(\tilde\sA-\sA^0)R_\lam g\Vert_2\le \Vert S_\lam^1-S^2_\lam\Vert\vareps_0K(M_0)\Vert g\Vert_2,$$
where we have used the definition of $\vareps_0$ (in (2.1)) and Proposition~2.2.  Set $\vareps_1(M_0)=(2K(M_0))^{-1}$ to conclude $\Vert S^1_\lam-S^2_\lam\Vert\le {1\over 2}\Vert S^1_\lam-S^2_\lam\Vert$.  Proposition~2.3 implies the above terms are finite for $\lam\ge 1$ and so we have
$$ \Vert S_\lam^1-S_\lam^2\Vert=0 \hbox{ for all }\lam\ge 1.\eqno(2.9)$$

To prove uniqueness we first use Krylov 
selection (Theorem 12.2.4 of [SV]) to see that it suffices to consider Borel 
strong Markov processes $((\Q^x_k)_{x\in S^0},X_t)$, $k=1,2$, where $\Q^x_k$ solves $MP(\tilde\sA,\delta_x)$, and 
to show that $\Q_1^x=\Q_2^x$ for all $x\in S^0$ (see the argument in the proof of Proposition~2.1 of [ABBP], but the situation here is a bit simpler as there is no killing).  If $S^k_\lam$ are the resolvent operators associated with $Q_k$, then (2.9) implies that
$$\eqalign{\int S_\lam^1f(x)\rho(x)d\mu(x)&=\int S_\lam^2f(x)\rho(x)d\mu(x)\cr&\hbox{  for all } f\in L^2,\hbox{compactly supported }\rho\in L^2,\hbox{ and }\lam\ge 1.}$$
For $f$ and $\lam$ as above this implies $S^1_\lam f(x)=S^2_\lam f(x)$ for Lebesgue a.e. $x$ and so 
for all $x$ by Proposition~2.4.  From this one deduces $\Q_1^x=\Q_2^x$ for all $x$ (e.g., see Theorem VI.3.2 of [B97]).   \qed

\ms
It remains to prove Propositions 2.2--2.4.  Propositions 2.3 and 2.4 follow along the lines of Propositions 2.3 and 2.4, respectively, of [ABBP], and are proved in Sections 5 and 6, respectively.  There are some additional complications in the present setting.  Most of the work, however, will go into the proof of Proposition~2.2 where a different approach than those in [ABBP] or [DP1] is followed. In [DP1] a canonical measure formula (Proposition~14 of that work) is used to 
represent and bound derivatives of the semigroups $P^i_tf(x)$ in (2.5) (see Lemma~3.8 below).    This approach will be refined (see, e.g., Lemmas~3.11 and 7.1 below) to give good estimates on the derivatives of the the actual transition densities using an integration by parts formula.   
The formula will convert spatial derivatives on the semigroup or density into differences involving Poisson random variables which can be used to represent the process with semigroup $P_t$ from which we are perturbing.  The construction is described in Lemma~3.4 below.  The integration by parts formula underlies  the proof of Lemma~7.1 and is explicitly stated in the simpler setting of first order derivatives in Proposition~8.1. 

 In [ABBP] we differentiate an explicit eigenfunction expansion for the resolvent of a killed squared Bessel process to get an asymptotically orthogonal expansion.  We have less explicit information about the semigroup $P_t$ of $\sA^0$ and so instead use Cotlar's Lemma (Lemma~2.13 below), to get a different asymptotically orthogonal expansion for the derivatives of the resolvent $R_\lam$--see the proof of Proposition~2.2 later in this section.

\ms

\ni{\bf Notation 2.6.} Set $\underline d=|Z\cap C|+|N_2|=|N_1^c|\le d$.  Here $|\cdot|$ denotes cardinality.

\ms

\noindent{\bf Convention 2.7.} {All
constants appearing in statements of results concerning the semigroup $P_t$ and its associated process may depend on $d$ and the constants
$\{b^0_j,\gamma^0_j:j\le d\}$, but, if $M_0$ is as in (2.4), 
these constants will be uniformly bounded for $M_0\le M$ for any
$M>0$. } 

\ms
We state an easy result on transition densities which will be proved in Section~3.

\proclaim Proposition 2.8.  The semigroup $(P_t,t\ge 0)$, has a jointly continuous transition density $p_t:S^0\times S^0\to[0,\infty)$, $t>0$. This density, $p_t(x,y)$ is $C^3$ on $S^0$ in each variable ($x$ or $y$) separately, and satisfies the following:\hfil\break
\ni(a) $p_t(y,x)=\hat p_t(x,y)$, where $\hat p_t$ is the transition density associated with $\hat \sA^0$ with parameters $\hat\gamma^0=\gamma^0$ and 
$$\hat b_j^0=\cases{-b^0_j& if $j\in N_1$\cr
b_j^0& otherwise.}$$
In particular
$$\int p_t(x,y)\mu(dy)=\int p_t(x,y)\mu(dx)=1.\eqno(2.10)$$
\ni(b) If $D^n_x$ is any $n$th order partial differential operator in $x$ and $0\le n\le 3$, then 
$$\sup_x|D_x^np_t(x,y)|\le c_{2.8}t^{-n-(d-\underline d)-\sum_{i\notin N_1}b_i^0/\gamma_i^0}\prod_{j\in N_2}[1+(y_j/t)^{1/2}]\quad \hbox{ for all }y\in S^0,\eqno(2.11)$$
and
$$\sup_y|D_y^np_t(x,y)|\le c_{2.8}t^{-n-(d-\underline d)-\sum_{i\notin N_1}b_i^0/\gamma_i^0}\prod_{j\in N_2}[1+
(x_j/t)^{1/2}]\quad \hbox{ for all }x\in S^0.\eqno(2.12)$$
\ni(c) If $0\le n\le 3$,
 $$\sup_x\int |D^n_y p_t(x,y)|d\mu(y)\le c_{2.8} t^{-n}.\eqno(2.13)$$
\ni (d) For all bounded Borel $f:S^0\to\R$, $P_tf\in C_b^2(S^0)$, and for $n\le 2$ and $D^n_x$ as in (b),
$$D_x^nP_tf(x)=\int D^n_xp_t(x,y)f(y)d\mu(y)\eqno(2.14)$$
and 
$$\Vert D^nP_t f\Vert_\infty\le c_{2.8}t^{-n}\Vert f\Vert_\infty.\eqno(2.15)$$

\ms

\ni{\bf Notation 2.9.} Throughout $\tD_x$ will denote one of the following first  or second order differential operators:
$$D_{x_j},\ j\le d,\ x_iD^2_{x_jx_j}, \ i\in Z\cap C, j\in R_i,\hbox{ or }x_jD^2_{x_jx_j}, j\notin N_1.$$

\ms
A deeper result is the following bound which sharpens Proposition 2.8.
\ms

\proclaim Proposition 2.10.  For $\tD_x$ as above and all $t>0$, 
$$\sup_x \int|\tD_xp_t(x,y)|\mu(dy)\le c_{2.10}t^{-1}.\eqno(2.16)$$
$$\sup_y \int|\tD_xp_t(x,y)|\mu(dx)\le c_{2.10}t^{-1}.\eqno(2.17)$$

\ni This is proved in Section 4 below.  The case $\tD_x=x_jD^2_{x_jx_j}$ for $j\in Z\cap C$ will be the most delicate.

\ms For $\tD$ as in Notation 2.9 and $t>0$, define an integral operator $T_t=T_t(\tD)$ by
$$T_tf(x)=\int \tD_x p_t(x,y)f(y)\mu(dy), \hbox{ for }f:S^0\to\R \hbox{ for which the integral exists}. \eqno (2.18)$$
By (d) above $T_t$ is  a bounded operator on $L^\infty$, but we will study these operators on $L^2(S^0,\mu)$. We will use the following well known elementary lemma; see [Ba], Theorem IV.5.1, for example, for a proof.

\proclaim Lemma 2.11.  Assume $K:S^0\times S^0\to\R$ is a measurable kernel on $S^0$ such that 
$$\Big\Vert\int |K(\cdot,y)|\mu(dy)\Big\Vert_\infty \le c_1\hbox{ and }\Big\Vert\int |K(x,\cdot)|\mu(dx)\Big\Vert_\infty\le c_2.$$
Then $Kf(x)=\int K(x,y)f(y)\mu(dy)$ is a bounded operator on $L^2$ with norm $\Vert K\Vert\le \sqrt{c_1c_2}$.

\proclaim Corollary 2.12. (a) For any $f\in L^2(\mu)$ and $t,\lam>0$,  $\Vert P_t f\Vert_2\le \Vert f\Vert_2$ and $\Vert R_\lam f\Vert_2\le \lam^{-1}$.\hfil\break
\ni(b) If $g\in C_b^2(S^0)\cap L^2(\mu)$ and $\sA^0g\in L^2(\mu)$, then $t\to P_tg$ is continuous in $L^2(\mu)$.

\proof (a) This is immediate from Lemma 2.11 and (2.10).\hfil\break
\ni(b) By $(MP(\sA^0,\nu ))$, if $0\le s< t$, then
$$\eqalign{ \Vert P_t g-P_s g\Vert_2&=\Big\Vert \int_s^t \sA^0P_r g\,dr\Big\Vert_2\cr
&\le \int_s^t \Vert P_r\sA^0g\Vert_2dr\cr
&\le (t-s)\Vert\sA^0g\Vert_2.}$$
We have used (a) in the last line.\qed

Proposition 2.10 allows us to apply Lemma 2.11 to $T_t$ and conclude
$$T_t \hbox{ is a bounded operator on }L^2\hbox{ with norm }\Vert T_t\Vert\le c_{2.10}t^{-1}\eqno (2.19)$$
Unfortunately this is not integrable near $t=0$ and so we can not integrate this bound to prove Proposition 2.2.  We must take advantage of some cancellation in the integral over $t$ and this is where we use Cotlar's Lemma:

\proclaim Lemma 2.13 (Cotlar's Lemma). Assume $\{U_j:j\in\Z_+\}$ are bounded operators on $L^2(\mu)$  and $\{a(j):j\in\Z\}$ are non-negative real numbers such that 
$$\Vert U_jU^*_k\Vert\vee\Vert U^*_jU_k\Vert\le a(j-k)^2\quad\hbox{all }j,k.\eqno(2.20)$$
Then 
$$\Big\Vert \sum_{j=0}^NU_j\Big\Vert\le A:=\sum_{j=-\infty}^\infty a(j)\quad\hbox{ for all }N.$$

\proof See, e.g., Lemma XI.4.1 in [T].\qed

The subspace $\sD_0$ in Proposition 2.2 will be
$$\sD_0=\{P_{2^{-j}}g:j\in\N, g\in C_b^2(S^0)\cap L^2(\mu), \sA^0g\in L^2(\mu)\}. \eqno(2.21)$$
As we can take $g\in C^2$ with compact support, denseness of $\sD_0$ in $L^2$ follows from Corollary~2.12(b).  To see that $\sD_0$ is a subspace, let $P_{2^{-j_i}}g_i\in \sD_0$ for $i=1,2$ with $j_2\ge j_1$.  If $\tilde g_1=P_{2^{-j_1}-2^{-j_2}}g_1$, then $\tilde g_1$ is in $L^2$ by Corollary 2.12 (a) and also in $C_b^2(S^0)$ by Proposition 2.8(d).  In addition, 
$$\Vert \sA^0 \tilde g_1\Vert_2=\Vert \P_{2^{-j_1}-2^{-j_2}}\sA^0g_1\Vert_2\le \Vert \sA^0g_1\Vert_2<\infty,$$
where we have used Corollary 2.12(a) again.  Hence $P_{2^{-j_1}}g_1=P_{2^{-j_2}}\tilde g_1$ where $\tilde g_1$ satisfies the same conditions as $g_1$.  Therefore 
$$P_{2^{-j_1}}g_1+P_{2^{-j_2}}g_2=P_{2^{-j_2}}(\tilde g_1+g_2)\in \sD_0.$$

\ni We show below how Cotlar's Lemma easily reduces Proposition 2.2 to the following result.

\proclaim Proposition 2.14.  There is an $\eta>0$ and $c_{2.14}$ so that if $\tD_x$ is any of the operators in Notation~2.9, then
$$ \eqalignno{\Vert T_s^*T_tf\Vert_2&\le c_{2.14}s^{-1-\eta/2}t^{-1+\eta/2}\Vert f\Vert_2\hbox{ and }\cr
\Vert T_sT^*_tf\Vert_2&\le c_{2.14}s^{-1-\eta/2}t^{-1+\eta/2}\Vert f\Vert_2\cr
&\hbox{for any }0<t\le s\le 2,\hbox{ and  any bounded Borel $f\in L^2(\mu)$.}&(2.22)}$$

\ni Assuming this result, we can now give the \hfil\break
\ni{\bf Proof of Proposition 2.2.} Fix a choice of $\tD_x$ (recall Notation 2.9), let $\lam\ge 1$, and for $k\in\Z_+$, define
$$U_k=U_k(\tD_x)=\int_{2^{-k}}^{2^{-k+1}}e^{-\lam s}T_s\,ds.$$
By (2.19), $U_k$ is bounded operator on $L^2$.  
Moreover if $k>j$ then 
$$\eqalign{ \Vert U_j^*U_k f\Vert_2&=\Big\Vert \int_{2^{-j}}^{2^{-j+1}}\Bigl[\int_{2^{-k}}^{2^{-k+1}}T_s^*T_tfdt\Bigr]\,ds\Big\Vert\cr
&\le \int_{2^{-j}}^{2^{-j+1}}\Bigl[\int_{2^{-k}}^{2^{-k+1}}c_{2.14}s^{-1-\eta/2}t^{-1+\eta/2}\,dt\Bigr]\,ds\Vert f\Vert_2\cr
&\le c_{2.14} 2^{-(\eta/2)(k-j)}\Vert f\Vert_2.}$$
If $k=j$ a similar calculation where the contributions to the integral from $\{s\ge t\}$ and $\{t\ge s\}$ are evaluated separately shows
$$\Vert U_j^*U_j f\Vert_2\le c_{2.14} \Vert f\Vert_2.$$
Cotlar's Lemma therefore shows that 
$$\Big\Vert \sum_{j=0}^N U_j\Big\Vert \le \sqrt{c_{2.14}}2(1-2^{-\eta/4})^{-1}
:= C(\eta)\hbox{ for all }N.\eqno(2.23)$$

Now let $f=P_{2^{-N}}g\in \sD_0$ where $g$ is as in the definition of $\sD_0$, and for $M\in\N$ set $h=h_M=P_{2^{-N}(1-2^{-M})}g$.  Then
$$\eqalign{\tD_xR_\lam f&=\tD_x\int_0^\infty e^{-\lam t}P_{t+2^{-M-N}}h\,dt\cr  
&=\exp(\lam2^{-M-N})\Bigl[\tD_x\Bigl[\int_{2^{-N-M}}^2 e^{-\lam u}P_uh\,du\Bigr]+\tD_x\Bigl[\int_2^\infty e^{-\lam u}P_uh\,du\Bigr]\Bigr]\cr
&=\exp(\lam2^{-M-N})\Bigl[\sum_{j=0}^{M+N}U_j h+\sum_{k=1}^\infty e^{-\lam k}\int_{k+1}^{k+2}e^{-\lam(u-k)}U_{u-k}(P_k h)\,du\Bigr].}$$
In the last line the bound (2.15) allows us to differentiate through the $t$ integral and (2.14) allows us to differentiate through the $\mu(dy)$ integral and conclude $\tD_xP_uh=T_uh$. A change of variables in the above now gives
$$\tD_xR_\lam f=\exp(\lam2^{-M-N})\Bigl[ \sum_{j=0}^{M+N}U_jh+\sum_{k=1}^\infty e^{-\lam k}U_0(P_kh)\Bigr].$$
So (2.23) shows that
$$\eqalignno{\Vert \tD_xR_\lam f\Vert_2&\le \exp(\lam2^{-M-N})C(\eta)\Bigl[ \Vert h_M\Vert_2+\sum_{k=1}^\infty e^{-\lam k} \Vert P_k h_M\Vert_2\Bigr]\cr
&\le \exp(\lam2^{-M-N})C(\eta)(1-e^{-\lam})^{-1}\Vert h_M\Vert_2.&(2.24)}$$
Corollary 2.12(b) shows that $\Vert h_M\Vert_2=\Vert P_{2^{-N}-2^{-N-M}}g\Vert_2\to \Vert f\Vert_2$ as $M\to\infty$.  Now let $M\to\infty$ in (2.24) to conclude
$$\Vert \tD_xR_\lam f\Vert_2\le C(\eta)(1-e^{-1})^{-1}\Vert f\Vert_2,$$
and the result follows. \qed

For Proposition 2.14, an easy calculation shows that for $0<s\le t$,
$$T^*_sT_tf(x)=\int K^{(1)}_{s,t}(x,y)f(y)d\mu(y)\hbox{ and }T_sT^*_tf(x)=\int K^{(2)}_{s,t}(x,y)f(y)d\mu(y),
\eqno(2.25)$$
where
$$K^{(1)}_{s,t}(y,z)=\int \tD_xp_s(x,y)\tD_xp_t(x,z)d\mu(x),\eqno(2.26)$$
and
$$K^{(2)}_{s,t}(x,y)=\int \tD_xp_s(x,z)\tD_yp_t(y,z)d\mu(z).\eqno(2.27)$$
A simple refinement of Lemma~2.11 (Lemma~3.16  proved at the end of Section~3) will show that (2.22) follows from
$$\eqalignno{\sup_y\int\int|K^{(i)}_{s,t}(x,y')|\, &|K^{(i)}_{s,t}(x,y)|\, d\mu(y')\, d\mu(x)\le c_{2.14}s^{-2-\eta}t^{-2+\eta}\cr
&\hbox{ for all }0<t\le s\le 2\hbox{ and }i=1,2.&(2.28)}$$
This calculation will reduce fairly easily to the case $N_2$ empty and $Z\cap C$ a singleton (see the proof of Proposition~2.14 at the end of Section~4 below).  Here there are essentially $4$ distinct choices of $\tD_x$, making our task one of bounding $8$ different $4$-fold integrals involving first and second derivatives of the transition density $p_t(x,y)$.  Fairly explicit formulae (see (4.7)--(4.9)) are available for all the derivatives except those involving the unique index $j$  in $Z\cap C$, and as a result Proposition~2.14 is easy to prove for all derivatives but those with respect to $j$ (Proposition~4.3).  Even here the first order derivatives are easily handled, leaving $\tD_x=x_jD_{x_jx_j}$.  This is the reason for most of the rather long calculations in Section~7.  In the special case $d=2$, of paramount importance to [DGHSS], one can avoid this case using the identity $\sA^0R_\lam f=\lam R_\lam f-f$, 
as is discussed in Section~8.

We give a brief outline of the rest of the paper.  
Section~3 studies the transition density associated with the resolvent 
in Proposition~2.2 for 
the key special case when $Z\cap C$ is a singleton and $N_2=\emptyset$.
This includes 
the canonical measure formulae for these densities (Lemmas~3.4 and 3.11) and 
the proof of Proposition~2.8.  In addition some important formulae for 
Feller branching processes with immigration, conditional on their value 
at time $t$, are proved (see Lemmas~3.2, 3.14 and Corollary~3.15).  
The section ends with an elementary result (Lemma~3.16) on 
integral operators on $L^2$.  In Section~4, 
the proofs of Propositions~2.14 and 2.10 are reduced to a series of 
technical bounds on the derivatives of the transition densities (Lemmas 4.5, 4.6 and 4.7). 
Most of the work here is in the setting of the key special case considered in Section~3,  and then at 
the end of Section~4 we show how the general case of Proposition~2.14 follows 
fairly easily thanks to the product structure in (2.5).  
Propositions~2.3 and 2.4 are proved in Sections~5 and 6, respectively.  
Lemmas~4.5--4.7 are finally proved in Section 7, thus completing 
the proof of Theorem~2.1 (and 1.4).  The key inequality in Section~7 
is Lemma~7.1 which comes from the integration by parts identity 
for the dominant term (see Proposition~8.1 for a simple special case). In Section~8 we describe how all
of this becomes considerably simpler in 
the $2$-dimensional setting required in [DGHSS].  
\ms

\subsec{3. The basic semigroups}
Unless otherwise indicated, in this section we work with the generator in (2.3) where $Z\cap C=\{d\}$ and $N_2=\emptyset$.  Taking $d=m+1$ to help distinguish this special setting,  this means we work with the generator
$$\sA^1=\Bigl[\sum_{j=1}^mb_j^0{\partial\over \partial x_j}+\gamma_j^0x_{m+1}{\partial^2\over \partial x_j^2}\Bigr]+b_{m+1}^0{\partial\over \partial x_{m+1}}+\gamma_{m+1}^0x_{m+1}{\partial^2\over \partial x_{m+1}^2},$$
with semigroup $P_t$ on the state space $S_m=\R^m\times \R_+$ ($m\in\N$).  Here we write $b=b^0_{m+1}$, $\gamma=\gamma^0_{m+1}$ and assume
$$\gamma_j^0>0, \ b_j^0\in\R\quad\hbox{for }j\le m,\qquad \hbox{and} \qquad \gamma>0, b>0.$$ Our Convention~2.7 on constants therefore means:

\noindent{\bf Convention 3.1.} { Constants appearing in statements of results may depend on $m$ and
$\{b^0_j,\gamma^0_j:\; j\le m+1\}$. 
If $$M_0=M_0(\gamma^0,b^0):= \max_{i\le m+1}(\gamma_i^0\vee
(\gamma_i^0)^{-1}\vee |b_i^0|)\vee 
(b^0_{m+1})^{-1} <\infty,$$
then these constants will be uniformly bounded for $M_0\le M$ for any fixed
$M>0$. } 

Note that $M_0\ge 1$.

It is easy to see that the martingale problem $MP(\sA^1,\nu)$ is well-posed for any initial law $\nu$  on $S_m$.  In fact, we now give an explicit formula for $P_t$.   Let $X_t=(X^{(1)}_t,\dots,X^{(m+1)})$ be a solution to this martingale problem.  By considering the associated SDE, we see that $X^{(m+1)}$ is a Feller branching diffusion (with immigration) with generator
$$\sA_0'=b{d\over dx}+\gamma x{d^2\over dx^2},\eqno(3.1)$$
and is independent of the driving Brownian motions of the first $m$ coordinates.  Let $\P_{x_{m+1}}$ be the law of $X^{(m+1)}$ starting at $x_{m+1}$ on $C(\R_+,\R_+)$.  
By conditioning on $X^{m+1}$ we see that the first $m$ coordinates are then a time-inhomogeneous Brownian motion.  Therefore if $I_t=\int_0^t X_s^{(m+1)}\,ds$ and $p_t(z)=(2\pi t)^{-1/2}e^{-z^2/2t}$, then (see (20) in [DP1])
$$P_tf(x_1,\dots,x_{m+1})\E_{x_{m+1}}\Bigl[\int f(y_1,\dots,y_m,X_t^{(m+1)})\prod_{j=1}^mp_{2\gamma_j^0I_t}(y_j-x_j-b_j^0t)\,dy_j\Bigr].\eqno (3.2)$$

If $x=(x_1,\dots,x_{m+1})=(x^{(m)},x_{m+1})\in S_m$, let $$\mu(dx)=x_{m+1}^{{b\over \gamma}-1}dx=dx^{(m)}\mu_{m+1}(dx_{m+1}).$$  Recall (see, e.g., (2.2) of [BP]) that $X^{(m+1)}$ has a symmetric density $q_t=q_t^{b,\gamma}(x,y)\ (x,y\ge 0)$ with respect to $\mu_{m+1}(dy)$, given by 
$$q_t^{b,\gamma}(x,y)=(\gamma t)^{-b/\gamma}\exp\Bigl\{{-x-y\over \gamma t}\Bigr\}\Bigl[\sum_{m=0}^\infty {1\over m!\Gamma(m+b/\gamma)}\Bigl({x\over \gamma t}\Bigr)^m\Bigl({y\over \gamma t}\Bigr)^m\Bigr],\eqno (3.3)$$
and associated semigroup $Q_t=Q^{b,\gamma}_t$.
Let 
$$\bar r_t(x_{m+1},y_{m+1},dw)=\P_{x_{m+1}}\Bigl(\int_0^t X_s^{(m+1)}\,ds\in dw|X_t^{(m+1)}=y_{m+1}\Bigr),$$ 
or more precisely a version of this collection of probability laws which is 
symmetric in $(x_{m+1},y_{m+1})$ and such that 
$(x_{m+1},y_{m+1})\to
\bar r_t(x_{m+1},y_{m+1},dw)$ is a jointly continuous map with respect
to the weak topology on the space of probability measures.  
The existence of such a version follows from Section IX.3 of [RY].  Indeed, Corollary 4.3 of the above states that if $\gamma=2$, then
$$\eqalignno{L(\lambda,x,y)&:=\int\exp\Bigl\{{-\lam^2\over 2}w\Bigr\}\bar r_1(x,y,dw)\cr
&=\cases{ {\lam\over \sinh \lam}\exp\Bigl\{\Bigl({x+y\over 2}\Bigr)(1-\lam \coth\lam)\Bigr\}I_\nu\Bigl({\lam\sqrt{xy}\over \sinh\lam}\Bigr)/I_\nu(\sqrt{xy})& if $xy>0$;\cr
\Bigl({\lam\over \sinh \lam}\Bigr)^{b/2}\exp\Bigl\{\Bigl({x+y\over 2}\Bigr)(1-\lam\coth\lam)\Bigr\}& if $xy=0$.}&(3.4)}
$$
Here $\nu={b\over 2}-1$ and $I_\nu(z)=\sum_{m=0}^\infty {1\over m!}{1\over \Gamma(m+\nu+1)}\Bigl({z\over 2}\Bigr)^{2m+\nu}$ is the modified 
Bessel function of the first kind of index $\nu>-1$.  
The continuity and symmetry of $L$ in $(x,y)$ gives the required continuous and symmetric version of $\bar r_1(x_{m+1},y_{m+1})$.  A scaling argument (see the proof of Lemma~3.2 below) gives the required version of $\bar r_t$ for general $\gamma>0$.  

Now define $r_t(x_{m+1},y_{m+1},dw)=q_t^{b,\gamma}(x_{m+1},y_{m+1})\bar r_t(x_{m+1},y_{m+1},dw)$, so that
$$(x_{m+1},y_{m+1})\to r_t(x_{m+1},y_{m+1},dw)\hbox{ is symmetric and weakly continuous}\eqno(3.5)$$
and
$$\eqalignno {\int\int&\psi(y_{m+1},w)r_t(x_{m+1},y_{m+1},dw)\mu_{m+1}(dy_{m+1})\cr
&=\E_{x_{m+1}}(\psi(X_t^{(m+1)},I_t))\hbox{ for all }x_{m+1}\ge 0\hbox{ and Borel }\psi:\R_+^2\to\R_+.&(3.6)}$$
(Weakly continuous means continuity with respect to the weak topology on the
space of probability measures.)
Combine (3.6) and (3.2) to conclude that $X$ has a transition density with respect to $\mu(dy)$ given by
$$\eqalignno{p_t(x,y)&=\int_0^\infty \prod_{j=1}^mp_{2\gamma_j^0w}(y_j-x_j-b_j^0t)r_t(x_{m+1},y_{m+1},dw)\cr
&=p_t(x^{(m)}-y^{(m)},x_{m+1},0,y_{m+1})\equiv p_t^0(x^{(m)}-y^{(m)},x_{m+1},y_{m+1}).&(3.7)}$$
Moreover if we set $b^0=(b_1^0,\dots,b_m^0)\in\R^m$ and write $p_t^{b^0}(x,y)$ for $p_t(x,y)$, then (3.5) implies
$$p^{b^0}_t(x,y)=p_t^{-b^0}(y,x)\hbox { for all }x,y\in S_m.\eqno(3.8)$$

The next result is a refinement of Lemma 7(b) of [DP1].

\proclaim Lemma 3.2. For any $p>0$ there is a $c_{3.2}(p)$ such that 
$$\int w^{-p}\bar r_t(x,y,dw)\le c_{3.2}(x+y+t)^{-p}t^{-p}\hbox{ for all }x,y\ge 0\hbox{ and }t>0.$$

\proof Assume first $\gamma=2$, $t=1$ so that we may use (3.4) to conclude (recall $\nu={b\over 2}-1$)
$$\eqalignno{{L(\lambda,x,y)\over L(\lambda,x,0)}&=\cases{\Bigl({\lam\over \sinh\lam}\Bigr)^{-\nu}{I_\nu\Bigl({\lam\over\sinh\lam}\sqrt{xy}\Bigl)\over I_\nu(\sqrt{xy})}\exp\{{y\over 2}(1-\lam\coth\lam)\}& if $x>0$;\cr
\exp\{{y\over 2}(1-\lam\coth\lam)\}& if $x=0$.}}$$
A bit of calculus shows 
$$\lam\coth\lam\ge 1,\ \alpha(\lam):={\lam\over \sinh\lam}\in[0,1]\hbox{ for all }\lam\ge 0,\eqno(3.9)$$
with the first inequality being strict  if $\lam>0$.  
The above series expansion shows that $I_\nu(\alpha z)\le \alpha^\nu I_\nu(z)$ for all $z\ge0,\ \alpha\in[0,1]$, and so using (3.9) in the above ratio bound, we get 
$$L(\lam,x,y)\le L(\lam,x,0)\hbox{ for all }\lam,x,y\ge 0.\eqno(3.10)$$ 
We have
$$\eqalign{\int w^{-p}\bar r_1(x,y,dw)&=2p\int_0^\infty \bar r_1(x,y,[0,u^2])u^{-2p-1}du\cr
&\le 2p\sqrt e \int_0^\infty L(u^{-1},x,y)u^{-2p-1}du\cr
&\le 2p\sqrt e \int_0^\infty L(u^{-1},x,0)u^{-2p-1}du\quad(\hbox{by }(3.10))\cr
&=c_p\int_0^\infty \Bigl({u^{-1}\over \sinh u^{-1}}\Bigr)^{b/2}\exp\{{x\over 2}(1-u^{-1}\coth u^{-1})\}u^{-2p-1}du\cr
&=c_p\int_0^\infty \Bigl({w\over \sinh w }\Bigr)^{b/2}\exp\{{-x\over 2}(w\coth w-1)\}w^{2p-1}dw\cr
&\le c_p\Bigl[\int_0^1\exp\{-cxw^2\}w^{2p-1}dw\cr
&\qquad\qquad+\int_1^\infty\exp\{-cxw-wb/2\}2^{b/2}w^{b/2+2p-1}dw\Bigr],\cr}$$
where in the last line $c>0$ and we have used (3.9), $\inf_{0\le w\le 1}{w\coth w-1\over w^2}=c_1>0$, and $\inf_{w\ge 1}{w\coth w-1\over w}=c_2>0$. For $x\le 1$ we may bound the above by (recall Convention 3.1)
$$c_p\Big[\int_0^1w^{2p-1}\,dw+\int_1^\infty (e^{-w}2w)^{b/2}w^{2p-1}\,dw\Big]\le c_1(p),$$
and for $x\ge 1$ we may, using  $(2we^{-w})^{b/2}\le 1$ for $w\ge 1$, bound it by
$$c_p\Bigl[\int_0^1\exp\{-cxw^2\}w^{2p-1}dw+\int_1^\infty e^{-cxw}w^{2p-1}dw\Bigr]\le c_2(p)x^{-p}.$$
These bounds show that $\int w^{-p}\bar r_1(x,y,dw)\le c(p)(1+x)^{-p}$ 
and so by symmetry in $x$ and $y$ we get
$$\int w^{-p}\bar r_1(x,y,dw)\le c(p)(1+x+y)^{-p}\hbox{ for all }x,y\ge 0.$$

For general $\gamma$ and $t$, $\hat X_s={2\over t\gamma}X_{ts}$ is as above with $\hat \gamma=2$ and $\hat b={2b\over \gamma}$.  We have $\int_0^tX_sds=\Bigl({t^2\gamma\over 2}\Bigr)\int_0^1\hat X_u\,du$, and so, using the above case, 
$$\eqalign{\int w^{-p}\bar r_t^{b,\gamma}(x,y,dw)&=\Bigl({t^2\gamma\over 2}\Bigr)^{-p}\int w^{-p}\bar r_1^{\hat b,2}
\Big({2x\over t\gamma},{2y\over t\gamma},dw\Big)\cr
&\le c(p) t^{-p}\Bigl({t\gamma\over 2}\Bigr)^{-p}\Bigl(1+{2x\over t\gamma}+{2y\over t\gamma}\Bigr)^{-p}\cr
&\le c(p)t^{-p}(t+x+y)^{-p}.}$$\qed

\ms

We observe that there exist $c_1,c_2$ (recall Convention~3.1 is in force) such that
$$c_1m^{b/\gamma-1}m!\leq \Gamma(m+ b/\gamma)\leq c_2m^{b/\gamma-1}m!
\hbox{ for all }m\in\N. \eqno (3.11)$$
To see this, suppose first $b/\gamma\equiv r\geq 1$ and use Jensen's inequality to obtain
$${\Gamma(m+r)\over m!}=\int x^{r}x^{m-1}e^{-x}{dx\over \Gamma(m)}m^{-1}\ge 
\Bigl(\int x \cdot x^{m-1}e^{-x}{dx\over \Gamma (m)}\Bigr)^{r}m^{-1}=m^{r-1}.$$
Next suppose $r\in [1,2]$ and again use Jensen's inequality to see that
$$\eqalign{{\Gamma(m+r)\over m!}&=\int x^{r -1}x^me^{-x}{dx\over \Gamma (m+1)}\cr&\le \Bigl(\int x\cdot x^me^{-x}{dx\over \Gamma (m+1)}\Bigr)^{r -1}=(m+1)^{r -1}\le 2m^{r -1}.}$$
These two inequalities imply (3.11) by using the identity $\Gamma(m+r+1)=(m+r)\Gamma(m+r)$ a finite number of times.

\ms

\proclaim Lemma 3.3. There is a $c_{3.3}$ so that for all $t>0$, $x_{m+1},y_{m+1}\ge 0$:\hfil\break
\ni(a) $q_t^{b,\gamma}(x_{m+1},y_{m+1})\le c_{3.3}[t^{-b/\gamma}+1_{(b/\gamma<1/2)}(x_{m+1}\wedge y_{m+1})^{1/2-b/\gamma} t^{-1/2}]$.

\ni (b) For all $t>0$, $(x,y)\rightarrow p_t(x,y)$ is continuous on $S_m^2$ and \hfil\break $\sup_{x,y\in S_m}p_t(x,y)\le c_{3.3}t^{-m-b/\gamma}$.

\ni (c) $\E_{x_{m+1}}(\exp(-\lam X_t^{(m+1)}))=(1+\lam\gamma t)^{-b/\gamma}\exp(-x_{m+1}\lam/(1+\lam\gamma t))$ for all $\lam>-(\gamma t)^{-1}$. 

\ni (d) If $0<p<b/\gamma$ then
$$\E_{x_{m+1}}((X_t^{(m+1)})^{-p})\le c_{3.3}\Bigl({p\over (b/\gamma)-p}+1\Bigr)(x_{m+1}+t)^{-p}.$$

\ni (e) $\E_{x_{m+1}}((X_t^{(m+1)})^2)\le c_{3.3}(x_{m+1}+t)^2$.

\ni (f) For any $p>0$, $\E_{x_{m+1}}\Bigr(\Bigr(\int_0^t X_s^{(m+1)}\,ds\Bigr)^{-p}\Bigr)\le c_{3.2}(p)t^{-p}(t+x_{m+1})^{-p}.$

\ni (g) $\sup_{y\ge 0}\Bigl|\int_0^\infty(x-y)x^{b/\gamma}D^2_{x}q_t(x,y)dx\Bigr|\le c_{3.3}$. 

\proof (a) 
If $q(x,y)=e^{-x-y}\sum_{m=0}^\infty {(xy)^m\over m!\Gamma(m+b/\gamma)}$, then $q_t(x,y)=(\gamma t)^{-b/\gamma}q(x/\gamma t,y/\gamma t)$ and it suffices to show
$$q(x,y)\le c_2(1+1_{(b/\gamma<1/2)}(x\wedge y)^{1/2-b/\gamma})\hbox{ for all }x,y\ge 0.\eqno(3.12)$$
By (3.11) and Stirling's formula we have
$$\eqalign{q(x,y)&\le c_1^{-1}e^{-x-y}\sum_{m=1}^\infty {(2m)!\over m!m!}2^{-2m}m^{1-b/\gamma}{(2\sqrt{xy})^{2m}\over (2m)!}+e^{-x-y}\Gamma (b/\gamma)^{-1}\cr
&\le c\Bigl[1+e^{-x-y+2\sqrt{xy}}\sum_{m=1}^\infty m^{1/2-b/\gamma}{(2\sqrt{xy})^m\over m!}e^{-2\sqrt{xy}}\Bigr]\cr
&\le c\Bigl[ 1+e^{-(\sqrt x-\sqrt y)^2}(2\sqrt{xy}+1)^{1/2-b/\gamma}\Bigr],}$$
where in the last line we used an elementary Poisson expectation calculation (for $b/\gamma\ge 1/2$ see Lemma~3.3 of [BP] ).  If $b/\gamma\ge 1/2$, the above is bounded and (3.12) is immediate.  Assume now that $p=1/2-b/\gamma>0$ and $x\ge y$.  Then the above is at most
$$\eqalign{c(1+e^{-(\sqrt x-\sqrt y)^2}(\sqrt x\sqrt y)^p)&\le c(1+{\sqrt y}^p({\sqrt y}^p+({\sqrt x}-\sqrt y)^p)e^{-(\sqrt x-\sqrt y)^2})\cr
&\le c(1+y^p+y^{p/2})\le c(1+y^p).}$$
This proves (3.12) and hence (a).

\ni(b) The continuity  follows easily from (3.7), the continuity of $q_t(\cdot,\cdot)$, the weak continuity of $\bar r_t(\cdot,\cdot,dw)$, Lemma 3.2 and 
dominated convergence.  Using Lemma 3.2 and (a) in (3.7), we obtain (recall Convention~3.1)
$$\eqalign{p_t(x,y)&\le c(t+x_{m+1}+y_{m+1})^{-m/2}q_t^{b,\gamma}(x_{m+1},y_{m+1})t^{-m/2}\cr
&\le c[t^{-m/2}t^{-b/\gamma}t^{-m/2}+1(_{b/\gamma<1/2)}(x_{m+1}+y_{m+1})^{b/\gamma-1/2}t^{-m/2+1/2-b/\gamma}\cr
&\qquad\qquad\qquad\times(x_{m+1}\wedge y_{m+1})^{1/2-b/\gamma}t^{-1/2}t^{-m/2}]\cr
&\le ct^{-m-b/\gamma}.}$$

\ni(c) This is well-known, and is easily derived from (3.3).

\ni(d) The expectation we need to bound equals
$$\eqalign{p\int_0^\infty &v^{-1-p}\P_{x_{m+1}}(X^{(m+1)}_t<v)\,dv\cr
&\le p\int_0^{x_{m+1}\vee t}v^{-1-p}e\E_{x_{m+1}}(e^{-X^{(m+1)}_t/v})\,dv+p\int_{x_{m+1}\vee t}^\infty v^{-1-p}\,dv\cr
&=e p\int_0^{x_{m+1}\vee t}v^{-1-p+b/\gamma}(v+\gamma t)^{-b/\gamma}e^{-x_{m+1}/(v+\gamma t)}\,dv+(x_{m+1}\vee t)^{-p}}$$
We have used (c) in the last line.
Set $c_r=\sup _{x\ge 0}x^re^{-x}$ and $M_1=M_0^2$.  If $x_{m+1}\ge t$, use 
$$(v+\gamma t)^{-b/\gamma}e^{-x_{m+1}/(v+\gamma t)}\le c_{b/\gamma}x_{m+1}^{-b/\gamma}\le (c_{M_1}+1)x_{m+1}^{-b/\gamma}$$
to bound the above expression for $\E_{x_{m+1}}((X^{(m+1)})^{-p})$ by
$$\eqalign{(c_{M_1}+1)&ep x_{m+1}^{-b/\gamma}\int_0^{x_{m+1}}v^{-1-p+b/\gamma}\,dv+x_{m+1}^{-p}\cr
&\le \Bigl({(c_{M_1}+1)ep\over b/\gamma -p}+1\Bigr)x_{m+1}^{-p}.}$$
On the other hand if $x_{m+1}<t$, the above expression for $\E_{x_{m+1}}((X^{(m+1)})^{-p})$is trivially at most
$$ep\int_0^tv^{-1-p+b/\gamma}\gamma^{-b/\gamma}t^{-b/\gamma}\,dv+t^{-p}=\Bigl({ep\gamma^{-b/\gamma}\over b/\gamma-p}+1\Bigr)t^{-p}.$$
The result follows from these two bounds.  

\ni (e) This is standard (e.g., see Lemma 7(a) of [DP1]).

\ni (f) Multiply the bound in Lemma~3.2 by $q_t(x,y)y^{b/\gamma-1}$ and integrate over $y$.  

\ni(g) Let $b'=\gamma/b$.
Since $q_t(x,y)= t^{-b/\gamma}q_1(x/t,y/t)$, a simple change of variables shows the quantity we need to bound is
$$\eqalign{\sup_{y\ge 0}\Bigl|&\int_0^\infty (x-\gamma y)x^{b'}D^2_xq_1(x,\gamma y)dx\Bigr|\cr
&=\sup_{y\ge 0}\Bigl|\sum_{m=0}^\infty e^{-y}{y^m\over m!}\int_0^\infty (x-y)x^{b'}D^2_x(e^{-x}x^m)\,dx/\Gamma(m+b')\Bigr|.}$$
Carrying out the differentiation and resulting Gamma integrals, we see the absolute value of the above summation equals
$$\eqalign{\Bigl|&\sum_{m=0}^\infty e^{-y}{y^m\over m!}\Bigl[(m+1+b')(m+b')-2m(m+b')+m(m-1)\cr
&\qq -y(m+b')+2my-ym(1-{b'\over m+b'-1})1_{(m\ge 1)}\Bigr]\cr
&\le\sum_{m=0}^\infty e^{-y}{y^m\over m!}\Bigl[b'(1+b')+yb'|-1+1_{(m\ge 1)}m(m+b'-1)^{-1}|\Bigr]\cr
&\le b'(1+b')+b'\sum_{m=1}^\infty e^{-y}{y^{m+1}\over (m+1)!}|1-b'|{m+1\over m+b'-1}+e^{-y}yb'\le c_{3.3}.}$$\qed

Now let $\{\P^0_x:x\ge 0\}$ denote the laws of the Feller branching
 process with generator $\sL^0 f(x)=\gamma xf''(x)$. If
$\omega\in C(\R_+,\R_+)$ let
$\zeta(\omega)=\inf\{t>0:\omega(t)=0\}$.  There is a unique
$\sigma$-finite measure $\N_0$ on
$$C_{ex}=\{\omega\in C(\R_+,R_+):\omega(0)=0,\
\zeta(\omega)>0,\ \omega(t)=0\ \forall t\ge \zeta(\omega)\}$$
such
that for each $h>0$, if $\Xi^h$ is a Poisson point process on
$C_{ex}$ with intensity $h\N_0$, then
$$X=\int_{C_{ex}}\nu\, \Xi^h(d\nu) \hbox{ has law } \P^0_h;\eqno(3.13)
$$
 see,
e.g., Theorem II.7.3 of [P] which can be projected down to the 
above situation
by considering the total mass function. Moreover for each $t>0$ we have (Theorems II.7.2(iii) and II.7.3(b) of [P])
$$\N_0(\{\nu:\nu_t>0\})=(\gamma t)^{-1}
\eqno(3.14)$$
and so we may define a probability on $C_{ex}$
by
$$ \P^*_t(A)=\frac{\N_0(\{\nu\in A:\nu_t>0\})}{\N_0(\{\nu:\nu_t>0)\}}.\eqno(3.15)$$
The above references in [P] also give the well-known
$$\P^*_t(\nu_t>x)=e^{-x/\gamma t},\eqno(3.16)$$
and so this together with (3.14) implies
$$\int_{C_{ex}} \nu_t\,  d\N_0(\nu)=1. \eqno(3.17)$$

The representation (3.13) leads to the following decomposition of $X^{(m+1)}$ from Lemma 10 of [DP1].  As it is consistent with the above notation, we will use $X^{(m+1)}$ to denote a Feller branching diffusion (with immigration) starting at $x_{m+1}$ and with generator given by (3.1), under the law $\P_{x_{m+1}}$.  

\proclaim Lemma 3.4. Let  $0\leq \rho \leq 1$.\hfil\break
\ni(a) We may assume
$$
X^{(m+1)}=X'_0+X_1,\eqno(3.18)$$
where $X'_0$ is  a diffusion with generator
$\sA'_0$ starting at $\rho
x_{m+1}$, $X_1$ is a diffusion with generator $\gamma x f''(x)$
starting at $(1-\rho)x_{m+1}\ge 0$, and $X_0', X_1$ are independent.
In addition, we may assume
$$X_1(t)=\int_{C_{ex}}\nu_t\,\Xi(d\nu)=\sum_{j=1}^{N_\rho(t)}e_j(t),\eqno(3.19)$$
where $\Xi$ is independent of $X_0'$ and is a Poisson point process on $C_{ex}$ with intensity
${(1-\rho)x_{m+1}}\N_0$, $\{e_j,j\in\N\}$ is an i.i.d.\  sequence with
common law $\P^*_t$, and $N_\rho(t)=\Xi(\{\nu:\nu_t>0\})$ is a Poisson random variable
(independent of the $\{e_j\}$) with mean $\frac{
(1-\rho)x_{m+1}}{t\gamma_{m+1}^0}$.\hfil\break
\noindent(b) We also have
$$\eqalignno{
\int_0^t X_1(s)ds &=\int_{C_{ex}} \int_0^t \nu_s\, ds\, 1_{(\nu_t\not
=0)}\Xi(d\nu)+\int_{C_{ex}} \int_0^t \nu_s\, ds\, 1_{(\nu_t
=0)}\Xi(d\nu)\cr
&  \equiv\sum_{j=1}^{N_\rho(t)}
r_j(t) +I_1(t),&(3.20)\cr
\int_0^tX_s^{(m+1)}ds&=\sum_{j=1}^{N_\rho(t)}r_j(t)+I_2(t),&(3.21)}
$$
where $r_j(t)=\int_0^t e_j(s)ds$ and 
$I_2(t)=I_1(t)+\int_0^tX'_0(s)ds$.

\ni{\bf Remark 3.5.} A double application of the decomposition in Lemma 3.4(a), first with general $\rho$ and then $\rho=0$ shows we may write
$$X_t^{(m+1)}=\tilde X'_0(t)+\sum_{j=1}^{N'_0(t)}e^2_j(t)+\sum_{j=1}^{N_{\rho}(t)}e^1_j(t),\eqno(3.22)$$
where $\tilde X'_0$ is as in Lemma 3.4(a) with $\rho=0$, $\{e^1_j(t), e^2_k(t),j,k\}$ are independent exponential variables with mean $(1/\gamma t)$, $N'_0(t)$, $N_{\rho}(t)$ are independent Poisson random variables
with 
means $\rho x_{m+1}/\gamma t$ and $(1-\rho)x_{m+1}/\gamma t$, respectively, and $(\tilde X_0(t), \{e^1_j(t), e^2_k(t)\}, N'_0(t), N_{\rho}(t))$ are jointly independent.  
The group of two sums of exponentials in (3.21) may correspond to $X_1(t)$  in (3.18) and (3.19), and  so we may use this as the decomposition in Lemma~3.4(a) with $\rho=0$.  
Therefore we may take $N_0$ to be $N'_0+N_\rho$, and hence may couple these decompositions so that
$$N_\rho\le N_0.\eqno(3.23)$$

The decomposition in Lemma 3.4 also gives a finer interpretation of the series expansion (3.2) for $q_t^{b,\gamma}(x,y)$, as we now show. Note that the decomposition (from (3.18),(3.19)),
$$X^{(m+1)}(t)=X'_0(t)+\sum_{j=1}^{N_\rho(t)}e_i(t),$$
where $X'_0$, $N_\rho(t)$ and $\{e_j(t)\}$ satisfy the distributional assumptions in (a), uniquely determines the joint law of $(X^{(m+1)}_t,N_\rho(t))$.  This can be seen by conditioning on $N_\rho(t)=n$.  Both this result and method are used in the following.

\proclaim Lemma 3.6. Assume $\phi:\R_+\times\Z_+\to\R_+$. If $N_\rho (t)$ is as in Lemma~3.4, then 
\item{(a)} $$\eqalign{\E_x(\phi(X_t^{(m+1)},N_0(t)))=\int_0^\infty&\sum_{n=0}^\infty\phi(y,n)( n!\Gamma(n+b/\gamma))^{-1}\cr           
&\times\Bigl({x\over \gamma t}\Bigr)^n\Bigl({y\over \gamma t}\Bigr)^n\exp\Bigl\{-{x\over \gamma t}-{y\over \gamma t}\Bigr\}(\gamma t)^{-b/\gamma}d\mu_{m+1}(y).\cr}$$
\item{(b)} $$\eqalign{\E_x(\phi(X_t^{(m+1)},N_{1/2}(t)))
=&\int_0^\infty\sum_{n=0}^\infty\sum_{k=0}^n {\choose n k} 2^{-n}\phi(y,k)( n!\Gamma(n+b/\gamma))^{-1}\cr          
&\times\Bigl({x\over \gamma t}\Bigr)^n\Bigl({y\over \gamma t}\Bigr)^n\exp\Bigl\{-{x\over \gamma t}-{y\over \gamma t}\Bigr\}(\gamma t)^{-b/\gamma}d\mu_{m+1}(y).\cr}$$

\proof (a) Set $x=0$ in (3.3) to see that $X'_0(t)$ has Lebesgue density
 $$\exp\{{-y\over\gamma t}\}(y/\gamma t)^{b/\gamma -1}(\gamma t)^{-1}\Gamma(b/\gamma)^{-1},$$ that is, has a gamma distribution with parameters $(b/\gamma,\gamma t)$.  It follows from Lemma 3.4(a), (3.16) and the joint independence of 
 $(\{e_j(t)\},X'_0(t),N_0(t))$ that, conditional on $N_0(t)=n$, $X_t^{(m+1)}$ has a gamma distribution with parameter $(n+b/\gamma,\gamma t)$.  This gives (a).

\ni(b) Apply (3.22) with $\rho =1/2$ to see that the decomposition in (3.18) for $\rho=1/2$ is given by (3.22) with $X'_0(t)=\tilde X'_0(t)+\sum_{j=1}^{N'_0(t)}e^2_j(t)$ and $e_j(t)=e^1_j(t)$. As in (a), conditional on $(N'_0(t), N_{1/2}(t))=(j,k)$, $X_t^{(m+1)}$ has a gamma distribution with parameters $(j+k+b/\gamma,\gamma t)$.  A short calculation (with $n=j+k$) now gives (b).\qed

\medskip

\ni{\bf Notation 3.7} Let $D^n$ denote any $n$th order partial differential operator on $S_m$ and 
let $D^n_{x_i}$ denote the $n$th partial derivative with respect to $x_i$. 

If $X\in C(\R_+,\R_+)$, $\nu^i\in C_{ex}$, $G:\R^2_+\to\R$ and $t\ge 0$, let
$$\Delta_t G(X,\nu^1)=\Delta^1_tG(X,\nu^1)=G\Big(\int_0^t X_s+\nu^1_s\,ds, X_t+\nu^1_t\Big)-G\Big(\int_0^1 X_s\,ds, X_t\Big),$$
$$\eqalign{\Delta^2_tG(X,\nu^1,\nu^2)=&G\Big(\int_0^t X_s+\nu^1_s+\nu^2_s\,ds, X_t+\nu^1_t+\nu^2_t\Big)-G\Big(\int_0^t X_s+\nu^1_s\,ds, X_t+\nu^1_t\Big)\cr
&-G\Big(\int_0^t X_s+\nu^2_s\,ds, X_t+\nu^2_t\Big)+G\Big(\int_0^t X_s\,ds, X_t\Big)\cr
=&\Delta_tG(X+\nu^1,\nu^2)-\Delta_t G(X,\nu^2),\cr}$$
$$\eqalign{\Delta_t^3G(X,\nu^1,\nu^2,\nu^3)=&G\Big(\int_0^t X_s+\nu^1_s+\nu^2_s+\nu^3_s\,ds, X_t+\nu^1_t+\nu^2_t+\nu^3_t\Big)\cr
&-\sum_{i=1}^3G\Big(\int_0^t X_s+\nu^1_s+\nu^2_s+\nu^3_s-\nu^i_s\,ds, X_t+\nu^1_t+\nu^2_t+\nu^3_t-\nu^i_t\Big)\cr
&+\sum_{i=1}^3 G\Big(\int_0^t X_s+\nu^i_s,X_t+\nu^i_t\Big)-G\Big(\int_0^t X_s\,ds,X_t\Big).}$$

\proclaim Lemma 3.8. If $f:S_m\to\R$ is a bounded Borel  function and $t>0$, then $P_tf\in C_b^3(S_m)$ and for $n\le 3$ 
$$\Vert D^nP_tf\Vert_\infty\le c_{3.8} \Vert f\Vert_\infty t^{-n}.\eqno(3.24)$$
Moreover if $f\in C_b(S_m)$, then for $n\le 3$, 
$$D^n_{x_{m+1}}P_tf(x)=\E_{x_{m+1}}\Bigl[ \int \Delta_t^n(G_{t, x^{(m)}}f)(X,\nu^1,\dots, \nu^n)\prod_{j=1}^nd\N_0((\nu^j))\Bigl],\eqno(3.25)$$
where 
$$G_{t,x^{(m)}}f(I,X)=\int_{\R^m}f(z_1,\dots,z_m,X)\prod_{j=1}^mp_{2\gamma_j^0 I}(z_j-x_j-b_j^0t)dz_j.\eqno (3.26)$$

\proof Proposition 14 and Remark 15 of [DP1] show that $P_tf\in C_b^2(S_m)$ and give (3.24) and (3.25) for $n\le 2$.  The proof there shows how to derive the $n=2$ case from the $n=1$ case and similar reasoning, albeit with more terms to check, allows one to derive the $n=3$ case from the $n=2$ case. \qed

Recall that $Q_t$ is the semigroup of $X^{(m+1)}$, the squared Bessel diffusion with transition density given by (3.3).  

\proclaim Corollary 3.9.  If $g:\R_+\to \R$ is a bounded Borel function  and $t>0$, then $Q_tg\in C_b^3(\R_+)$ and for $n\le 3,$ 
$$\Vert D^nQ_tf\Vert_\infty\le c_{3.8} \Vert g\Vert_\infty t^{-n}.\eqno(3.27)$$

\proof  Apply Lemma 3.8 to $f(y)=g(y_{m+1})$.  \qed

\ni{\bf Notation 3.10} For $t,\delta>0$, $x^{(m)}\in\R^m$, $y\in S_m$, $1\le j\le m$, $I,X\ge 0$, define
$$G^\delta _{t,x^{(m)},y}(I,X)=\prod_{i=1}^m p_{\delta+2\gamma_i^0 I}(x_i-y_i+b_i^0 t) q_\delta^{b,\gamma}(y_{m+1},X).$$

\proclaim Lemma 3.11. (a) For each $t>0$ and $y\in S_m$, the functions $x\to p_t(x,y)$ and $x\to p_t(y,x)$ are in $C_b^3(S_m)$, and if $D_x^n$ denotes 
any  $n$th order partial differential operator in the $x$ variable, then
$$|D^n_yp_t(x,y)|+|D_x^np_t(x,y)|\le c_{3.11}t^{-n-m-b/\gamma}\hbox{ for all }x,y\in S_m\hbox{ and }0\le n\le 3.\eqno(3.28)$$
\ni(b) $\sup_x\int|D_x^np_t(x,y)|\mu(dy)\le c_{3.8}t^{-n}$ for all $t>0$ and $0\le n\le 3$.

\ni (c) {\it For $n\le 3$ and $t>0$, $y\to D_x^np_t(x,y)$ and $x\to D_y^n p_t(x,y)$ are in $C_b(S_m)$.} 

\ni(d) {\it For $n=1,2,3$, $$D^n_{x_{m+1}}p_t(x,y)=\lim_{\delta\downarrow 0}\E_{x_{m+1}}\Bigl(\int \Delta_t^nG^\delta_{t,x^{(m)},y}(X,\nu^1,\dots,\nu^n)\prod_{i=1}^nd\N_0(\nu^i)\Bigr).$$}

\proof (a) By the Chapman-Kolmogorov equations, 
$p_t(x,y)=\E_x(p_{t/2}(X_{t/2},y))$.  
The result for $x\to p_t(x,y)$ now follows from Lemma 3.3(b) and Lemma 3.8 with $f(x)=p_{t/2}(x,y)$.  
By (3.8) it follows for $y\to p_t(x,y)$.

\ni(b) For $n=1,2,3$, $N\in\N$ and $x\in S_m$, let $f(y)=\sgn(D_x^np_t(x,y))1_{(|y|\le N)}$.  Then
$$\eqalign{\int|D_x^np_t(x,y)|1_{(|y|\le N)}d\mu(y)&=\Bigl|\int D_x^np_t(x,y) f(y)d\mu(y)\Bigr|\cr
&=|D_x^nP_t f(x)|,}$$
where the last line follows by dominated 
convergence, the uniform bound in (3.28) and the fact that $f$ has compact support.  An application of (3.24) implies 
$$\int|D_x^np_t(x,y)|1_{(|y|\le N)}d\mu(y)\le c_{3.8}t^{-n},$$ and the result follows upon letting $N\to\infty$.

\ni (c) $$\eqalignno{ D_x^np_t(x,y)&=D^n_x\int p_{t/2}(x,z)p_{t/2}(z,y)d\mu(z)\cr
&=\int D_x^np_{t/2}(x,z)p_{t/2}(z,y)d\mu(z)&(3.29)}$$
by dominated convergence and the uniform bounds in (a).  The integrability of \hfil\break
 $D_x^np_{t/2}(x,z)$ with respect to $\mu(dz)$ (from (b)) and the fact that $p_{t/2}(z,\cdot)\in C_b(S_m)$ allow us to deduce the continuity of (3.29) in $y$ from dominated convergence.  
Now use symmetry,  i.e.,  (3.8), to complete the proof. 

\ni (d) If $y,z\in S_m$, $\delta>0$, let 
$$f^{y,\delta}(z)=\prod_1^mp_\delta(z_i-y_i)q_\delta(y_{m+1},z_{m+1}).$$
$f$ is bounded and continuous in $z$ by Lemma 3.3(a) (with a bound depending on $y,\delta$).   Let $n\le 3$.  The uniform bounds in (a) and integrability of $f^{y,\delta}$ allow us to apply dominated 
convergence to differentiate through the integral  and conclude  
$$\eqalignno{D^n_{x_{m+1}}P_tf^{y,\delta}(x)&=\int\int D^n_{x_{m+1}}p_t(x,z)f^{y,\delta}(z)\,d\mu(z)\cr
&\rightarrow D^n_{x_{m+1}}p_t(x,y)\hbox{ as }\delta\downarrow 0.&(3.30)}$$
In the last line we have used (c).  Now note that if $I,X\ge 0$, then
$$\eqalignno{G_{t,x^{(m)}}f^{y,\delta}(X,I)&\equiv\int_{\R^m}f^{y,\delta}(z,X) \prod_{j=1}^mp_{2\gamma_j^0 I}(z_j-x_j-b_j^0t)dz_j\cr
=&q_\delta(y_{m+1},X)\prod_{j=1}^mp_{\delta+2\gamma_j^0 I}(y_j-x_j-b_j^0t)=G^\delta_{t,x^{(m)},y}(I,X).&(3.31)}$$
Use this in (3.25) to conclude that 
$$D^n_{x_{m+1}}P_t f^{y,\delta}(x)=\E_{x_{m+1}}\Bigl(\int\int \Delta_t^n G^\delta_{t,x^{(m)},y}(X,\nu^1,\dots,\nu^n)\prod_{i=1}^nd\N_0(\nu^i)\Bigr).\eqno(3.32)$$
Combine (3.32) and (3.30) to derive (d).\qed

\proclaim Lemma 3.12. (a) For each $t>0$ and $y\in \R_+$, the functions $x\to q_t(x,y)$ and $x\to q_t(y,x)$ are in $C_b^3(\R_+)$, and
$$|D^n_yq_t(x,y)|\le c_{3.12}t^{-n-b/\gamma}[1+\sqrt{y/t}]\hbox{ for all }x,y\ge 0\hbox{ and }0\le n\le 3.\eqno(3.33)$$

\ni(b) $\sup_{x\ge 0}\int|D_x^nq_t(x,y)|\mu_{m+1}(dy)\le c_{3.12}t^{-n}$ for all $t>0$ and $0\le n\le 3$.

\proof  This is a minor modification of the proofs of Lemma 3.11 (a),(b).  Use the bound (from Lemma~3.3(a))
$$\eqalign{q_{t/2}(\cdot,y)&\le c_{3.3}(t^{-b/\gamma}+1_{(b/\gamma<1/2)}y^{1/2-b/\gamma}t^{-1/2})\cr
&\le c_{3.3}t^{-b/\gamma}(2+\sqrt{y/t})}$$ in place of Lemma~3.3(b), and (3.27) in place of (3.24), in the above argument. Of course this is much easier and can also be derived by direct calculation from our
series expansion for $q_t$.\qed

\ms

Expectation under $\P^*_t$ is denoted by $\E^*_t$ and we let $e(s),s\ge 0$, denote the canonical excursion process under this probability.

\proclaim Lemma 3.13.  If $0<s\le t$, then $$\E^*_t(e(s)|e(t))=(s/t)((s/t)e(t)+2\gamma(t-s))\le e(t)+2\gamma s,\quad \P^*_t-\hbox{ a.s.}$$

\proof Let $\P^0_h$ denote the law of the diffusion with generator $\gamma x {d^2\over dx^2}$ starting at $h$. If $f\in C_b(\R_+)$ has compact support, then Proposition~9 of [DP1] and (3.14) show that 
$$\eqalignno{\E^*_t(e(s)f(e(t)))&=\lim_{h\downarrow 0}h^{-1}\E^0_h(X_s f(X_t)1_{(X_t>0)})\gamma t\cr 
&=\lim_{h\downarrow 0} {t\over s} \int\int {\gamma sq^{0,\gamma}(h,y)\over h}yq^{0,\gamma}_{t-s}(y,x)f(x)\, dy\,  dx.&(3.34)}$$
Here we extend the notation in (3.3) by letting $q^{0,\gamma}$ denote the absolutely continuous part of the transition kernel for $\E^0$ (it also has an atom at $0$).  We have also extended the convergence in [DP1] slightly as the functional $e(s)f(e(t))$ is not bounded but this extension is justified by a uniform integrability argument--the approximating functionals are $L^2$ bounded.  
By (2.4) of [BP]
$$\eqalign{{\gamma s\over h}q^{0,\gamma}(h,y)&=\exp\Bigl\{{-h-y\over \gamma s}\Bigr\}\sum_{m=0}^\infty {1\over (m+1)!}\Bigl({h\over\gamma s}\Bigl)^m{1\over m!}\Bigl({y\over \gamma s}\Bigr)^m(\gamma s)^{-1}\cr
&\to \exp(-y/\gamma s)(\gamma s)^{-1}\quad\hbox{ as }h\downarrow 0,}$$
and also ${\gamma s\over h}q^{0,\gamma}_s(h,y)\le (\gamma s)^{-1}$.  Dominated convergence allows us to take the limit in (3.34) through the integral and deduce that
$$\E_t^*(e(s)f(e(t)))={t\over s}\int\int \exp(-y/\gamma s)(\gamma s)^{-1}yq^{0,\gamma}_{t-s}(y,x)\,dy f(x)\,dx.$$
By (3.16) we conclude that
$$\eqalign{\E^*_t(e(s)|e(t)=x)&=(t/s)^2e^{x/\gamma t}\int e^{-y/\gamma s}y q_{t-s}^{0,\gamma }(y,x)\,dy.}$$
By inserting the above series expansion for $q^{0,\gamma}_{t-s}(y,x)$ and calculating the resulting gamma integrals, the result follows. 

\proclaim Lemma 3.14.  Let $\nu=b/\gamma-1$, and $\kappa_\nu(z)={I'_\nu\over I_\nu}(z)z+1$ for $z\ge 0$, where  $\kappa_\nu(0)\equiv \lim_{z\downarrow 0}\kappa_\nu(z)=\nu+1$. Then
$$\eqalign{\E_{z_{m+1}}\Bigl(\int_0^tX^{(m+1)}_s\,ds\Bigl|X_t^{(m+1)}=y\Bigr)&=\kappa_\nu\Bigl({2\sqrt{z_{m+1}y}\over t\gamma}\Bigr){t^2\gamma\over 6}+{(z_{m+1}+y)t\over 3}\cr
&\le c_{3.14}[t^2+t(z_{m+1}+y)].}$$

\proof Write $X$ for $X^{(m+1)}$ and $z$ for $z_{m+1}$. The scaling argument used in the proof of Lemma~3.2 allows us to assume $t=1$, $\hat\gamma=2$ and $\hat b=2b/\gamma$.  Dominated convergence implies that
$$\E_z\Bigl(\int_0^1 X_s\,ds\Bigl|X_1=y)=\lim_{\lam\to 0+}-\lam^{-1}{d\over d\lam}\E_z\Bigl(\exp\Bigl((-\lam^2/2)\int_0^1 X_s\,ds\Bigr)\Bigl| X_1=y\Bigl).$$

The right side can be calculated explicitly from the two formulae in (3.4) and after some calculus we arrive at $$\E_z\Bigl(\int_0^1X_s\,ds\Bigl|X_1=y\Bigr)={\kappa_\nu(\sqrt{zy})\over 3}+{(z+y)\over 3},$$
where $\nu=\hat b/2-1=b/\gamma-1$.  This gives the required equality.

To obtain the bound (recall Convention 3.1) it suffices to show
$$\kappa_\nu(z)\le c(\eps)(1+z)\hbox{ for all }z\ge 0\hbox{ and }\nu+1\in[\eps,\eps^{-1}].$$
Set $\alpha=\nu+1$ and recall that
$$I_\nu(z)=(z/2)^\nu\sum_{n=0}^\infty {1\over n!\Gamma(n+\alpha)}(z/2)^{2n},$$
and
$$zI'_\nu(z)=2(z/2)^\nu\sum_{n=1}^\infty {1\over (n-1)!\Gamma(n+\alpha)}(z/2)^{2n}+\nu I_\nu(z).$$
Taking ratios of the above and setting $w=(z/2)^2$, we see that it suffices to show
$$\sum_{n=1}^\infty{w^n\over (n-1)!\Gamma(n+\alpha)}\le c_0(\eps)\Bigl[\sum_{n=0}^\infty {w^{n+1/2}\over n!\Gamma(n+\alpha)}\Bigr]\hbox{ for all }w\ge 0,\eps\le\alpha\le \eps^{-1}.\eqno(3.35)$$
We claim in fact that 
$$\eqalignno{{w^n\over (n-1)!\Gamma(n+\alpha)}\le {1\over 2}(\alpha^{-1/2}\vee 1)\Bigl[&{w^{n-1/2}\over (n-1)!\Gamma(n-1+\alpha)}+{w^{n+1/2}\over n!\Gamma(n+\alpha)}\Bigr]\cr
&\qquad\quad\hbox{for all }n\ge 1, w\ge 0, \alpha\in[\eps,\eps^{-1}].&(3.36)}$$
Assuming this, (3.35) will follow with $c_0(\eps)=\eps^{-1/2}$ by summing (3.36) over $n\ge 1$. The proof of (3.36) is an elementary application of the quadratic formula once factors of $w^{n-1/2}$ are canceled.  
\qed

\proclaim Corollary 3.15.  If $N_0(t)$ is as in Lemma 3.4 then
$$\E_{z_{m+1}}\Bigl(\int_0^t X_s^{(m+1)}\,ds\Bigl| X^{(m+1)}_t, N_0(t)\Bigl)\le c_{3.15}[t^2(1+N_0(t))+t(X^{(m+1)}_t+z_{m+1})].$$

\proof We write $X$ for $X^{(m+1)}$ and $z$ for $z_{m+1}$.  Recall the decomposition (3.21):
$$\int _0^t X_s\,ds =\int_0^tX'_0(s)\,ds+\int\int_0^t\nu_sds1_{(\nu_t=0)}\Xi(d\nu)+\sum_{i=1}^{N_0(t)}\int _0^t e_i(s)ds,$$
where the second integral is independent of $(\{e_j(t),j\in\N\},N_0(t)=\int_0^t1_{(\nu_t>0)}\Xi(d\nu), X'_0)$ by elementary properties of Poisson point processes and the independence of $X'_0$ and $\Xi$.  Therefore
$$\eqalign{\E_z&\Bigl(\int_0^t X_sds\Bigl|X'_0(t),N_0(t),\{e_j(t),j\in\N\}\Bigr)\cr
&=\E_z\Bigl(\int_0^t X'_0(s)ds\Bigl|X'_0(t)\Bigr)+\E_z\Bigl(\int\int_0^t\nu_sds1_{(\nu_t=0)}\Xi(d\nu)\Bigl)\cr
&\phantom{=\E_z\Bigl(\int_0^t X'_0(s)}+\sum_{i=1}^{N_0(t)}\E_z\Bigl(\int_0^te_i(s)ds\Bigl|X'_0(t),N_0(t),\{e_j(t),j\in\N\}\Bigr)\cr
&=\E_z\Bigl(\int_0^t X'_0(s)ds\Bigl|X'_0(t)\Bigr)+\E_z\Bigl(\int\int_0^t\nu_sds\Xi(d\nu)\Bigl)+\sum_{i=1}^{N_0(t)}\E_z\Bigl(\int_0^t e_i(s)ds|e_i(t)\Bigr).}$$
In the last line we have used the independence of $X'_0$ and $\Xi$, of $N_0(t)$ and $\{e_j,j\in \N\}$, and the joint independence of the $\{e_j\}$ (see Lemma 3.4).  Now use Lemmas 3.13 and 3.14 to bound the last and first terms, respectively, and note the second term is bounded by the mean of $\int_0^tX_1(s)\,ds$, where  $X_1$ is as in (3.19).  This bounds the above by
$$\eqalign{c_{3.14}&[t^2+tz+tX'_0(t)]+\E_z\Bigl(\int_0^t X_1(s)ds\Bigr)+\sum_{i=1}^{N_0(t)}(e_i(t)t+\gamma t^2)\cr
&\le c[t^2+tX_t+zt]+\gamma N_0(t)t^2.\cr}$$
Condition the above on $\sigma(N_0(t),X_t)$ to complete the proof.\qed

We now return to the general setting of Propositions 2.2 and 2.8.

\ni{\bf Proof of Proposition 2.8.}  From (2.5) we may write
$$p_t(x,y)=\prod_{i\in Z\cap C}p^i_t(x_{(i)},y_{(i)})\prod_{j\in N_2}q^j_t(x_j,y_j),$$
where $p^i_t$ are the transition densities from Lemma~3.11 and $q^j_t$ are the transition
densities from Lemma~3.12. The joint continuity and smoothness in each variable is immediate from
these properties for each factor (from Lemma~3.11 and (3.3)).  (a) is also immediate from (3.8). The first part of  (b) is also clear from the above factorization
and the upper bounds in Lemmas~3.11(a) and 3.12(a).  The second part of (b) is then immediate from (a). 
(c) also follows from Lemmas~3.11(b) and 3.12(b) and a short calculation.  

(d) is an exercise in differentiating through the integral.  As we will be doing a lot of this in the future we risk boring the reader by outlining the proof here and refer to this argument for such manipulations hereafter.  Let $f$ be 
a bounded Borel function on $S^0$ and $0\le n$.  Note first that the right-hand side of (2.14), $I(x)$, (finite by (c)) is continuous in $x$.
To see this, choose a unit vector $e_j$, set $x'_i=x_i$ if $i\neq j$ and $x'_j$ variable,  and note that for $h>0$,
$$\eqalign{|I(x+he_j)-I(x)|&\le \int \int_{x_j}^{x_j+h}|D_{x_j}D^n_xp_t(x',y)|dx'_jf(y)\mu(dy)\cr
&\le \Vert f\Vert_\infty c_{2.8}t^{-n-1}h.}$$
We have used (c) in the above.

Let $f_N(y)=f(y)1_{(|y|\le N)}$.  By the integrability in (c), the left-hand side of (2.14) equals
$$\lim_{N\to\infty}\int D_x^np_t(x,y)f_N(y)d\mu(y)=\lim_{N\to\infty}D_x^n\int p_t(x,y)f_N(y)d\mu(y),$$
where the differentiation through the integral over a compact set is justified by the bounds in (b) and dominated convergence.  The bound in (c) shows this convergence is uniformly bounded (in $x$).  
For definiteness assume $n=2$ and $D^2_x=D^2_{x_i x_j}$ for $i\neq j$.  By the above convergence and dominated convergence we get
$$\eqalign{\int_0^{x_j'}\int _0^{x_i'}\int D^2_xp_t(x,y)f(y)d\mu(y)dx_idx_j&=\lim_{N\to\infty}\int_0^{x_j'}\int _0^{x_i'}D^2_x\int p_t(x,y)f_N(y)d\mu(y)dx_idx_j\cr
&=\int [p_t(x,y)|_{x_i=0}^{x_i=x'_i}|_{x_j=0}^{x^j=x'_j}]f(y)d\mu(y)}$$
Now differentiate both sides with respect $x'_j$ and then $x'_i$ and use the Fundamental Theorem of Calculus and the continuity of $\int D_x^2 p_t(x,y) f(y) d\mu(y)$, noted above, to obtain (2.14).  This shows $P_tf\in C_b^2$ as continuity in $x$ was established above.  Finally (2.15) is now immediate from (2.13) and (2.14).\qed

We end this section with a corollary to Lemma 2.11.

\proclaim Lemma 3.16. Assume $K:S_0\times S_0\to \R$ satisfies
$$|K|^2:= \Big\Vert\int\int |K(x,y')|\, |K(x,\cdot)|\mu(dy')\mu(dx)\Big\Vert_\infty<\infty.\eqno(3.37)$$
Then $Kf(\cdot)=\int K(\cdot,y)f(y)\mu(dy)$ is bounded on $L^2(\mu)$ and its norm satisfies $\Vert K\Vert\le |K|$.  

\proof Let $K^*K(\cdot,\cdot)$ denote the integral kernel associated with the operator $K^*K$.
The hypothesis implies that 
$$\Big\Vert \int |K^*K(y, \cdot)|\, \mu(dy)\Big\Vert_\infty\leq |K|^2.$$
By Lemma 2.11 and the fact that $K^*K$ is symmetric, we have
$$\Vert K\Vert^2=\Vert K^*K\Vert \leq |K|^2.$$
\qed

\subsec{4. Proofs of Propositions 2.14 and 2.10}

By Lemma 3.16, Proposition 2.14 will follow from (2.28).  We restate this latter inequality explicitly.  Recall that $\tD_x$ is one of the first or second order partial 
differential operators listed in Notation~2.9.  (2.28) then becomes
$$\eqalignno{\sup_y\int&\Bigl[\int\Bigl|\int\tD_zp_s(z,x)\tD_zp_t(z,y')\mu(dz)\Bigr|\mu(dy')\Bigl|\int\tD_zp_s(z,x)\tD_zp_t(z,y)\mu(dz)\Bigr|\Bigr]\mu(dx)\cr
&\le c_{2.14}s^{-2-\eta}t^{_2+\eta}\hbox{ for all }0<t\le s\le 2,&(4.1)}$$
and
$$\eqalignno{\sup_y\int&\Bigl[\int\Bigl|\int\tD_xp_s(x,z)\tD_{y'}p_t(y',z)\mu(dz)\Bigr|\mu(dy')\Bigl|\int\tD_xp_s(x,z)\tD_yp_t(y,z)\mu(dz)\Bigr|\Bigr]\mu(dx)\cr
&\le c_{2.14}s^{-2-\eta}t^{_2+\eta}\hbox{ for all }0<t\le s\le 2,&(4.2)}$$
We have stated these conditions for bounded times for other potential uses; in our case we will verify (2.28) for all $0<t\le s$.
Recall also our 
Convention 3.1  for constants applies to $c_{2.14}$ and $\eta$.  

Until otherwise indicated, we continue to work in the setting of the last section and use the notation introduced there.  In particular, Convention~3.1 will be in force and the differential operators in Notation~2.9 are 
$$\tD_x=D_{x_i},\  i\le m+1,\ \hbox{or }\tD_x=x_{m+1}D^2_{x_i},\ i\le m+1.\eqno(4.3)$$

In [DP1] a number of bounds were obtained on the derivatives of the semigroup $P_tf$; (3.24) in the last section was one such bound.  Propositions~16 and 17 of [DP1] state 
there exists a $c_{4.1}$ such that for all $\tD_x$ as in Notation~2.9, $t>0$ and bounded Borel function $f$,
$$\sup_{x\in S_m}|\tD_xP_tf(x)|\le c_{4.1}t^{-1}\Vert f\Vert_\infty.\eqno(4.4)$$
Although these results are stated for $m=1$ in [DP1], the same argument works in $m+1$ dimensions (see,  for example,
 Proposition 20 of [DP1]).  As a simple consequence of this result we get:

\proclaim Lemma 4.1.  For all $t>0$, and all $\tD_x$ as in (4.3) 
$$\sup_{x\in S_m}\int|\tD_{x}p_t(x,y)|\mu(dy)\le c_{4.1}t^{-1},\eqno (4.5)$$
and
$$\sup_{x_{m+1}\ge 0}\int [x_{m+1}|D^2_{x_{m+1}}q_t(x_{m+1},y)|+|D_{x_{m+1}}q_t(x_{m+1},y)|\, ]\, \mu_{m+1}(dy)\le c_{4.1}t^{-1}.\eqno (4.6)$$

\proof Apply (4.4) and (2.14) to $f(y)=\sgn(\tD_xp_t(x,y))$ to obtain (4.5), and to $f(y)=\sgn(\tD_{x_{m+1}}q_t(x_{m+1},y_{m+1}))$ to obtain (4.6).\qed

One of the ingredients we will need is a bound like (4.5) but with the integral on the left with respect 
to $x$ instead of $y$.  For derivatives with respect to $x_i, i\le m$, this is straightforward as we now show.

By differentiating through the integral in (3.7) we find for $i\le m$,
$$D_{x_i}p_t(x,y)=\int_0^\infty {(y_i-x_i-b_i^0t)\over 2\gamma^0_iw}\prod_{j=1}^mp_{2\gamma_j^ow}(y_j-x_j-b_j^0t)r_t(x_{m+1},y_{m+1},dw),\eqno(4.7)$$
 $$\eqalignno{D^2_{x_i}p_t(x,y)=\int_0^\infty\Bigl[ &{(y_i-x_i-b_i^0t)^2\over 2\gamma^0_iw}-1\Bigr](2\gamma^0_iw)^{-1}\cr
 &\quad\times\prod_{j=1}^mp_{2\gamma_j^0w}(y_j-x_j-b_j^0t)r_t(x_{m+1},y_{m+1},dw),&(4.8)}$$
and
$$\eqalignno{D^3_{x_i}p_t(x,y)=\int_0^\infty \Bigl[&{(y_i-x_i-b_i^0t)^3\over (2\gamma^0_iw)}-{ 3(y_i-x_i-b_i^0t)}\Bigr](2\gamma_0^iw)^{-2}\cr
&\quad\times\prod_{j=1}^mp_{2\gamma_j^0w}(y_j-x_j-b_j^0t)r_t(x_{m+1},y_{m+1},dw).&(4.9)}$$
Integration through the integral is justified by the bounds in Lemma~3.2 and dominated convergence.

\proclaim Lemma 4.2. For all $t>0$, and $i\le m$,

\ni(a) $\int|D_{z_i} p_t(z,y)|dz^{(m)}\le c_{4.2}t^{-1/2}(t+z_{m+1}+y_{m+1})^{-1/2}q_t(y_{m+1},z_{m+1})$.

\ni(b) For all $0\le p\le 2$, $\int z_{m+1}^p|D_{z_i}p_t(z,y)|\mu(dz)\le c_{4.2}t^{-1/2}(t+y_{m+1})^{p-1/2}$.

\ni(c) For  $\tD_z=z_{m+1}D^2_{z_i}$ or $D_{z_i}$,

$$\sup_y\int |\tD_z p_t(z,y)|\mu(dz)\le c_{4.2}t^{-1},\eqno(4.10)$$
and
$$\sup_z\int |\tD_z p_t(z,y)|\mu(dy)\le c_{4.2}t^{-1}.\eqno(4.11)$$

\proof (a)  By (4.7) the integral in (a) is
$$\eqalign{\int& \Bigl|\int_0^\infty\Bigl[ {(y_i-z_i-b_i^0t)\over 2\gamma^0_iw}\Bigr]
 \prod_{j=1}^mp_{2\gamma_j^0w}(y_j-z_j-b_j^0t)r_t(z_{m+1},y_{m+1},dw)\Bigr|dz^{(m)}\cr
 &\le c\int_0^\infty w^{-1/2}r_t(y_{m+1},z_{m+1},dw)\cr
 &\le ct^{-1/2}(t+x_{m+1}+y_{m+1})^{-1/2}q_t(y_{m+1},z_{m+1}).}$$
  where in the first inequality we used the symmetry of $r_t$ (recall (3.5)) and in the last inequality we have used Lemma 3.2.  
 
 \ni(b) Integrate the inequality in (a) to bound the integral in (b) by
 $$\eqalign{c &t^{-1/2}\int  z_{m+1}^p(t+z_{m+1}+y_{m+1})^{-1/2}q_t(y_{m+1},z_{m+1})\mu_{m+1}(dz_{m+1})\cr
 &\le ct^{-1/2}(t+y_{m+1})^{-1/2}\E_{y_{m+1}}((X_t^{(m+1)})^p)\cr
 &\le ct^{-1/2}(t+y_m)^{p-1/2}.}$$
 In the last line we used Lemma~3.3(e).
 
 \ni(c) For $\tD_z=D_{z_i}$ (4.10) follows from (b) upon taking $p=0$.  The other cases are similarly proved, now using (4.8) for the second order derivatives.  ((4.11) is also immediate from Lemma~4.1.)  \qed

\ni Consider first (2.28) for $\tD_x=D_{x_i}$ or $x_{m+1}D^2_{x_i}$ for some $i\le m$.
 
\proclaim Proposition 4.3. If $\tD_x=D_{x_i}$ or $x_{m+1}D^2_{x_i}$ for some $i\le m$, then (2.28) holds with $\eta=1/2$.  

\proof Consider (4.1) for $\tD_x=x_{m+1}D^2_{x_i}$.  We may as well take $i=1$.  Assume $0<t\le s$  
and let
$$J=\int\Big|\int \tD_zp_s(z,x)\tD_zp_t(z,y')\mu(dz)\Bigr|\mu(dy').$$
Then
$$\eqalign{J\le& \int\Bigl|\int z_{m+1}D^2_{y'_1}p_s( y'_1,z_2,\dots,z_{m+1},x)z_{m+1}D^2_{z_1}p_t(z,y')\mu(dz)\Bigr|\mu(dy')\cr
&+\int\Bigl|\int z_{m+1}[D^2_{y'_1}p_s(y'_1,z_2,\dots,z_{m+1},x)-D^2_{z_1}p_s(z,x)]z_{m+1}D^2_{z_1}p_t(z,y')\mu(dz)\Bigr|\mu(dy')\cr
\equiv&J_1+J_2.}$$
To evaluate $J_1$ do the $dz_1$ integral first and use (4.8) to see
$$\int z_{m+1}D^2_{z_1}p_t(z,y')\,dz_1=0,$$
and so $J_1=0$. (Lemma~3.2 handles integrability issues.)  

Let $J_2'$ and  $J_2''$
 denote the contribution to the integral defining $J_2$ from $\{z_1\le y _1'\}$  and $\{z_1\ge y'_1\}$, respectively.  
Then by (4.8),
$$\eqalign{J'_2\le& \int \int z_{m+1}^2\Bigl|\int\int 1_{(z_1\le z'\le y'_1)}D^3_{z'}p_s(z',z_2,\dots,z_{m+1},x)\cr
&\times \Bigl(\int_0^\infty\Bigl[{(y'_1-z_1-b_1^0t)^2\over 2\gamma_1^0w}-1\Bigr](2\gamma_1^0 w)^{-1}\prod_{j=1}^mp_{2\gamma_j^0w}(y'_j-z_j-b_j^0t)\Bigr)\cr
&\phantom{\times \Bigl(\int_0^\infty\Bigl[{(y'_1-z_1-b_1^0t)^2\over 2\gamma_1^0w}-1\Bigr]}r_t(z_{m+1},y'_{m+1},dw)dz^{(m)}dz'\Bigr|\mu_{m+1}(dz_{m+1})\mu(dy').}$$
Do the $z_1$ integral first in the above and if $X={-y'_1+z'+b_1^0t\over \sqrt{2\gamma_1^0w}}$, note the absolute value of this integral is
$$\eqalignno{\Bigl| \int_{-\infty}^{z'}&\Bigl[{(y'_1-z_1-b_1^0t)^2\over 2\gamma_1^0w}-1\Bigr]p_{2\gamma_1^0w}(y'_1-z_1-b_1^0t)\,dz_1\Bigr|\cr
&=\Bigl|\int_{-\infty}^X[v^2-1]p_1(v)\,dv\Bigr|\cr
&\le c|X|e^{-X^2/2}\cr
&\le c|y'_1-z'-b_1^0t|p_{2\gamma_1^0w}(y'_1-z'-b_1^0t),&(4.12)}$$
where the first inequality follows by an elementary calculation--consider $|X|\ge 1$ and $|X|<1$ separately, note that $\int_0^\infty (v^2-1)p_1(v)\,dv=0$, and in the last case use $$\Bigl|\int_{-\infty}^Xp_1(v)[v^2-1]dv\Bigr|=\int_0^{|X|}p_1(v)[1-v^2]dv.$$  
Take the absolute value inside the remaining integrals, then integrate 
over  $dy'_2...dy'_m$, and use (4.9) to express the third order derivative.  This and (4.12) lead to
$$\eqalign{J'_2\le & \int\int z^2_{m+1}\Bigl\{\int\int 1_{(z'\le y'_1)}\Bigl[\int_0^\infty\Bigl[{|x_1-z'-b_1^0s|^3\over (2\gamma_1^0w')}+{|x_1-z'-b_1^0s|}\Bigr](2\gamma_1^0w')^{-2}\cr
&\phantom{\int\int z^2_{m+1}}\times p_{2\gamma_1^0w'}(x_1-z'-b_1^0s)\prod_{j=2}^mp_{2\gamma_j^0w'}(x_j-z_j-b_1^0s)dz_jr_s(z_{m+1},x_{m+1},dw')\Bigr]\cr
&\phantom{\int\int z^2_{m+1}}\times \int_0^\infty c{|-y'_1+z'+b_1^0t|\over{2\gamma_1^0w}}
p_{2\gamma_1^0 w}(y'_1-z'-b_1^0t)r_t(z_{m+1},y'_{m+1},dw)dy'_1dz'\Bigr\}\cr
&\phantom{\int\int z^2_{m+1}\Bigl\{\int\int 1_{(z'\le y'_1)}\Bigl[\int_0^\infty\Bigl[{|x_1-z'-b_1^0s|^3\over (2\gamma_1^0w')^2}}\times\mu_{m+1}(dz_{m+1})\mu_{m+1}(dy'_{m+1}).}$$
Do the trivial integral
over  $dz_2\dots dz_m$ and then consider the $dy'_1\,dz'$ integral of the resulting integrand.  If $z''=z'-x_1+b_1^0s$ and $y''_1=y'_1-x_1-b_1^0(t-s)$, this integral equals
$$\eqalign{\int&\int 1_{(z''\le y''_1+b_1^0t)}\Bigl[{|z''|^3\over (2\gamma_1^0w')}+{|z''|}\Bigr](2\gamma_1^0w')^{-2}p_{2\gamma_1^0w'}(z''){|y''_1-z''|\over 2\gamma_1^0 w}\cr
&\qq\qq\times p_{2\gamma_1^0w}(y''_1-z'')dy''_1\,dz''\cr
&\le c(2\gamma_1^0 w')^{-3/2}(2\gamma_1^0w)^{-1/2}.}$$
Use this in the above bound on $J'_2$ and the symmetry of $r_s(z_{m+1},x_{m+1},\cdot)$ in $(z_{m+1},x_{m+1})$ (recall (3.5)), and conclude
$$\eqalign{J'_2&\le c \int\int\int\int z_{m+1}^2( w')^{-3/2}w^{-1/2}r_s(x_{m+1},z_{m+1},dw')r_t(z_{m+1},y'_{m+1},dw)\cr
&\qq\qq\times\mu_{m+1}(dy'_{m+1})\mu_{m+1}(dz_{m+1})\cr
&=c\int\int z_{m+1}^2(w')^{-3/2}\E_{z_{m+1}}\Bigl(\Bigl(\int_0^tX_r^{(m+1)}dr\Bigr)^{-1/2}\Bigr)r_s(x_{m+1},z_{m+1},dw')\mu_{m+1}(dz_{m+1})\cr
&\le c \int\int (w')^{-3/2}z_{m+1}^2t^{-1/2}(t+z_{m+1})^{-1/2}r_s(x_{m+1},z_{m+1},dw')\mu_{m+1}(dz_{m+1})\cr
&\phantom{\le c \int\int (w')^{-3/2}z_{m+1}^2t^{-1/2}(t+z_{m+1})^{-1/2}r_s(x_{m+1},z_{m+1},dw')}\hbox{(by Lemma 3.3(f))}\cr
&\le ct^{-1/2} \E_{x_{m+1}}\Bigl(\Bigl(X^{(m+1)}_s\Bigr)^{3/2}\Bigl(\int_0^s X_r^{(m+1)}dr\Bigr)^{-3/2}\Bigr)\ \qq\hbox{(recall (3.6))}\cr
&\le t^{-1/2}\E_{x_{m+1}}\Bigl(\Bigl(X^{(m+1)}_s\Bigr)^2\Bigr)^{3/4}\E_{x_{m+1}}\Bigl(\Bigl(\int_0^s X_r^{(m+1)}dr\Bigr)^{-6}\Bigr)^{1/4}.
}$$
Another application of Lemma 3.3(e,f) now shows 
$$J_2'\le ct^{-1/2}(x_{m+1}+s)^{3/2}(x_{m+1}+s)^{-6/4}s^{-6/4}=ct^{-1/2}s^{-3/2}.$$
Symmetry (switching $z_1$ and $y'_1$ in the integral defining $J'_2$ amounts to switching the sign of $b_1^0$) gives the same bound on $J''_2$ and hence for $J$.  This implies that the left-hand side of (4.1) is at most
$$\eqalign{c&t^{-1/2}s^{-3/2}\sup_y\int\int z_{m+1}^2|D^2_{z_1}p_s(z,x)||D^2_{z_1}p_t(z,y)|\mu(dx)\mu(dz)\cr
&\le ct^{-1/2}s^{-5/2}\sup_y\int z_{m+1}|D^2_{z_1}p_t(z,y)|\mu(dz)\hbox{ (by Lemma 4.1)}\cr
&\le ct^{-3/2}s^{-5/2},}$$
where (4.10) is used in the last line.

This completes the proof of (4.1) with $\eta=1/2$ for $\tD_x=x_{m+1}D^2_{x_1}$.  The proof for $D_{x_1}$ is similar and a bit easier.  Finally, very similar arguments (the powers change a bit in the last part of the bound on $J'_2$) will verify (4.2) for these operators.\qed

Recall the notation $\hat p_t$ from Proposition~2.8.

\proclaim Lemma 4.4. If $D_x^n$ and $D^k_y$ are $n$th and $k$th order partial differential operators in $x$ and $y$, respectively, then for all $t>0$, $k,n\le 2$, $D^n_{y}D^k_{x}p_t(x,y)$ exists, is bounded, is continuous in each variable separately, and equals
$$\int D^k_{y}\hat p_{t/2}(y,z)D^n_{x}p_{t/2}(x,z)\,\mu(dz).$$
For $n\le 1$, $D_y^nD_x^kp_t(x,y)$ is jointly continuous.  

\proof From (3.29) we have 
$$D^n_{x}p_t(x,y)=\int D^n_{x}p_{t/2}(x,z)\hat p_{t/2}(y,z)\mu(dz).$$
Apply (2.14), with $f(z)=D^n_{x}p_{t/2}(x,z)$ and $\hat p_{t/2}$ in place of $p_t$, to differentiate with respect to $y$ through the integral and derive the above identity.  Uniform boundedness  in $(x,y)$, and continuity in each variable separately follows from the boundedness in Lemma 3.11(a), the $L^1$-boundedness in Lemma~3.11(b) and 
dominated convergence. If $n=1$ the uniform boundedness of the above derivative implies continuity in $y$ uniformly in $x$ and hence joint continuity. \qed 

We now turn to the verification of (2.28) and Proposition 2.10 for $\tD_x=D_{x_{m+1}}$ or $x_{m+1}D^2_{x_{m+1}}$.  The argument here seems to be much harder (at least in the second order case) and so we will reduce its proof to three technical bounds whose proofs are deferred to Section 7.  Not surprisingly these proofs will rely on the representations in Lemmas~3.4 and 3.11 as well as the other explicit expressions obtained in Section~3 such as Lemmas 3.6 and 3.14.

\proclaim Lemma 4.5. There is a $c_{4.5}$ such that for all $t>0$:\hfil\break
\ni(a) If $-(2M_0^2)^{-1}\le p\le 1/2$ then for all $j\le m$ and all $y\in S_m$,
$$\int z_{m+1}^p|D_{z_{m+1}}p_t(z,y)|\mu(dz)\le c_{4.5}t^{-1/2}(t+y_{m+1})^{p-1/2},\eqno(4.13)$$
and for all $z\in S_m$,
$$\int y_{m+1}^p|D_{z_{m+1}}p_t(z,y)|\mu(dy)\le c_{4.5}t^{-1/2}(t+z_{m+1})^{p-1/2}.\eqno(4.14)$$
\ni(b) If $0\le q\le 2$ and $-(2M_0^2)^{-1}\le p\le 2$, then for all $j\le m$ and all $y\in S_m$
$$\int|y_j-z_j|^qz_{m+1}^p|D^2_{z_{m+1}}p_t(z,y)|\mu(dz)\le c_{4.5}t^{q/2-1}(t+y_{m+1})^{p+q/2-1}.\eqno(4.15)$$\hfil\break
\ni(c) If $0\le q\le 2$, then for all $j\le m$ and $z\in S_m$,
$$\int|y_j-z_j|^qz_{m+1}|D^2_{z_{m+1}}p_t(z,y)|\mu(dy)\le c_{4.5}[t^{q-1}+t^{q/2-1}z_{m+1}^{q/2}].\eqno(4.16)$$
\ni(d) $$\sup_y\int z_{m+1}^{3/2}|D^3_{z_{m+1}}p_t(z,y)|\mu(dz)\le c_{4.5}t^{-3/2}.\eqno(4.17)$$
\ni(e) If $0\le p\le 1/2$, then for all $(y^{(m)},z_{m+1})\in S_m$, 
$$\int  z_{m+1}y_{m+1}^p|D_{y_{m+1}}D^2_{z_{m+1}}p_t(z,y)|dz^{(m)}\mu_{m+1}(dy_{m+1})\le c_{4.5}t^{p-2},\eqno(4.18)$$
and for all $j\le m$,
$$\int z_{m+1}y_{m+1}^p|D_{z_j}D^2_{z_{m+1}}p_t(z,y)|dz^{(m)}\mu_{m+1}(dy_{m+1})\le c_{4.5}t^{p-2}.\eqno(4.19)$$
\ni(f) If $0\le p\le 3/2$, then for all $j\le m$,
$$\sup_y\int z_{m+1}^p|D_{z_j}D^2_{z_{m+1}}p_t(z,y)|\mu(dz)\le c_{4.5}t^{p-3}.\eqno(4.20)$$

\proclaim Lemma 4.6. There is a $c_{4.6}$ such that for all $0<t\le s$,
$$t^{b/\gamma}\int_0^{\gamma t }\Bigl[\int |D_{z_{m+1}}p_s(z,y)|dz^{(m)}\Bigr]dz_{m+1}\le c_{4.6}t/s,\eqno(4.21)$$
and
$$t^{b/\gamma}\int_0^{\gamma t }\Bigl[\int |D^3_{z_{m+1}}p_s(z,y)|dz^{(m)}\Bigr]dz_{m+1}\le c_{4.6}ts^{-3}\le c_{4.6}s^{-2}\eqno(4.22)$$

\proclaim Lemma 4.7. There is a $c_{4.7}$ such that if $1\le p\le 2$, then for all $t>0$, $w>0$ and $y^{(m)}\in\R^m$,
$$\eqalignno{\int\int (&1_{(y_{m+1}\le w\le z_{m+1})}+1_{(z_{m+1}\le w\le y_{m+1})})z_{m+1}^p|D^2_{z_{m+1}}p_t(z,y)|\mu(dz)\mu_{m+1}(dy_{m+1})\cr
&\le c_{4.7}[1_{(w\le \gamma t)}t^{p-2+b/\gamma}+1_{(w>\gamma t)}t^{-1/2}w^{p-3/2+b/\gamma}].&(4.23)}$$

\ms

Assuming these results we now verify (2.28) and Proposition~2.10 for $\tD_x=D_{x_{m+1}}$ or $x_{m+1}D^2_{x_{m+1}}$.  The analogue of Proposition~2.10 is immediate.

\proclaim Proposition 4.8. For $\tD_x$ as in (4.3), and all $t>0$, 
$$\sup_x \int|\tD_xp_t(x,y)|\mu(dy)\le c_{4.8}t^{-1}.\eqno(4.24)$$
$$\sup_y \int|\tD_xp_t(x,y)|\mu(dx)\le c_{4.8}t^{-1}.\eqno(4.25)$$

\proof This is immediate from Lemma~4.1, Lemma~4.2(c), Lemma~4.5(a) (with $p=0$) and Lemma~4.5(b) (with $q=0$ and $p=1$).\qed

\proclaim Proposition 4.9.  (2.28) holds for $\tD_x=D_{x_{m+1}}$ and $\eta=(2M_0^2)^{-1}$.  

\proof Define $\eta$ as above and fix $0<t\le s\le 2$.  We will verify (4.1) even with the absolute values taken inside all the integrals.  By (4.14) with $p=0$,
$$\int|D_{z_{m+1}}p_t(z,y')|\mu(dy')\le c_{4.5}t^{-1/2}(t+z_{m+1})^{-1/2}\le c_{4.5}t^{-1+\eta}z_{m+1}^{-\eta}.$$
In the last inequality we used $\eta\le 1/2$ (recall $M_0\ge 1$).  Therefore the left side of (4.1) (even with absolute values inside all the integrals) is at most
$$\eqalignno{\sup_y&\ c_{4.5}\int\Bigl[\int |D_{z_{m+1}}p_s(z,x)|z_{m+1}^{-\eta}\mu(dz)t^{-1+\eta}\cr
&\phantom{\sup_yc_{4.5}\int\Bigl[\int |D_{z_{m+1}}}\times\Bigl[\int  |D_{z'_{m+1}}p_s(z',x)|\,|D_{z'_{m+1}}p_t(z',y)|\mu(dz')\Bigr]\Bigr]\mu(dx)\cr
&\le c_{4.5}^2 s^{-1-\eta}t^{-1+\eta}\sup_y\int\Bigl[\int |D_{z_{m+1}}p_s(z',x)|\mu(dx)\Bigr]|D_{z_{m+1}}p_t(z',y)|\mu(dz'),&(4.26)}$$
where we have used (4.13) with $p=-\eta$ in the last line.  Now apply (4.14) to the above integral in $x$ and then (4.13) to the integral in $z'$, both with $p=0$, and conclude that (4.26) is at most
$$cs^{-2-\eta}t^{-2+\eta},$$
as required.  The derivation of (4.2) (with absolute values inside in the integral) is almost the same. One starts with (4.13) with $p=0$ to bound the integral in $y'$ as above, and then uses (4.14) with $p=-\eta$ to bound the resulting integral in $z$. \qed

It remains to verify (2.28) for $\tD_x=x_{m+1}D_{x_{m+1}}^2$.  This is the hard part of the proof and we will not be able to take the absolute values inside 
the integrals in (2.28).  

\proclaim Lemma 4.10.  For $j=1,2$, $\int |D^j_{z_{m+1}}p_t(z,y)|dz^{(m)}<\infty$ for all $z_{m+1}>0$ and $y\in S_m$, and
$$\eqalignno{&\int D^j_{z_{m+1}}p_t(z,y)dz^{(m)}=D^j_{z_{m+1}}\int p_t(z,y)dz^{(m)}=D^j_{z_{m+1}}q_t(z_{m+1},y_{m+1})\cr
&\phantom{\int D^j_{z_{m+1}}p_t(z,y)dz^{(m)}=D^j_{z_{m+1}}\int p_t(z,y)dz^{(m)}}\hbox{ for all }z_{m+1}>0\hbox{ and }y\in S_m.&(4.27)}$$

\proof We give a proof as this differentiation is a bit delicate, and the result is used on a number of occasions.  Set $j=2$ as $j=1$ is slightly easier. Fix $y\in S_m$ and $t>0$.  By Lemma~4.5(b,d),
$$\int |D^2_{z_{m+1}}p_t(z,y)|+1_{(z_{m+1}>\eps)}|D^3_{z_{m+1}}p_t(z,y)|\mu(dz)<\infty\hbox{ for all }\eps>0.\eqno(4.28)$$
We claim
$$z_{m+1}\to\int D^2_{z_{m+1}}p_t(z,y)dz^{(m)}\equiv F(z_{m+1})\hbox{ is continuous on }\{z_{m+1}>0\}.\eqno(4.29)$$
Note that $F(z_{m+1})<\infty$ for almost all $z_{m+1}>0$ by (4.28).  The Fundamental Theorem of Calculus and Lemma~3.11(a) imply that if $z'_{m+1}>z_{m+1}>0$,
 then
$$\eqalign{\int& |D^2_{z_{m+1}}p_t(z^{(m)},z'_{m+1},y)-D^2_{z_{m+1}}p_t(z^{(m)},z_{m+1},y)|\,dz^{(m)}\cr
&\le \int \int_{z_{m+1}}^{z'_{m+1}}|D^3_wp_t(z^{(m)},w,y)|dw\,dz^{(m)}\to 0\hbox{ as }z'_{m+1}\to z_{m+1}\hbox{ or }z_{m+1}\to z'_{m+1},}$$
where dominated convergence and (4.28) are used to show the convergence to 0.
This allows us to first conclude that 
$$\int |D^2_{z_{m+1}}p_t(z,y)|dz^{(m)}<\infty\hbox{ for all }z_{m+1}>0,\eqno(4.30)$$
and in particular,
$F(z_{m+1})$ is finite for all $z_{m+1}>0$, and also that $F$ is continuous.

The differentiation through the integral now proceeds as in the proof of Proposition~2.8(d) given in Section~3 (using the Fundamental Theorem of Calculus). The last equality follows from (3.7) and the definition of $r_t$.\qed

\proclaim Lemma 4.11. There is a $c_{4.11}$ so that for $\tD_z\equiv\tD_{z_{m+1}}=D_{z_{m+1}}$ or $z_{m+1}D^2_{z_{m+1}}$, all $t>0$ and all $y'\in S_m$,
$$\Bigl|\int \tD_zp_t(z,y')\mu(dz)\Bigr|\le c_{4.11}(t+y'_{m+1})^{-1}\eqno(4.31)$$

\proof Use (4.27) to see that 
$$\eqalign{\int\tD_zp_t(z,y')\mu(dz)
&=\int \tD_{z_{m+1}}q_t(z_{m+1},y'_{m+1})\mu_{m+1}(dz_{m+1}).}$$
Changing variables, we must show
$$\Bigl|\int \tD_zq_t(z,y)z^{b/\gamma-1}dz\Bigr|\le c_{4.10} (t+y)^{-1}.\eqno(4.32)$$
The arguments are the same for either choice of $\tD_z$ so let 
us take $\tD_z=D_z$ for which the algebra is slightly easier. Let $w=y/\gamma t$ and $x=z/\gamma t$.  
By differentiating the power series (3.3) the left side of (4.32) is then
$$\eqalign{\Bigl|\int_0^\infty &(\gamma t)^{-1}e^{-w}\sum_{m=0}^\infty {w^m\over m!\Gamma(m+b/\gamma)} e^{-x}[mx^{m-1}-x^m] x^{b/\gamma-1}dx\Bigr|\cr
&=\Bigl|(\gamma t)^{-1}\Bigr(\sum_{m=1}^\infty e^{-w}{w^m\over m!\Gamma(m+b/\gamma)}[m\Gamma(m-1+b/\gamma)-\Gamma(m+b/\gamma)]\Bigr)-(\gamma t)^{-1}e^{-w}\Bigr|\cr
&=\Bigl|(\gamma t)^{-1}\Bigl(\sum_{m=1}^\infty e^{-w}{w^m\over m!}{1-b/\gamma\over m+b/\gamma-1}\Bigr)-(\gamma t)^{-1}e^{-w}\Bigr|\cr
&\le c_1(\gamma t)^{-1}[\E(N(w)^{-1}1_{(N(w)\ge 1)})+e^{-w}],}$$
where $N(w)$ is a Poisson random variable with mean $w$, $c_1$ satisfies Convention~3.1, and we have used $(n+b/\gamma-1)^{-1}\le cn^{-1}$ for all $n\in\N$.  An elementary calculation (e.g. Lemma~3.3 of [BP]) bounds the above by 
$$c_2t^{-1}(1\wedge w^{-1}+e^{-w})\le 2c_2t^{-1}[1\wedge w^{-1}]\le c_3(t+y)^{-1}\Bigr].$$ (4.32) follows. \qed

\proclaim Proposition 4.12.  (2.28) holds for $\tD_x=x_{m+1}D^2_{x_{m+1}}$ and $\eta=1/2$.

\proof Consider general $0<\eta<1$ for now.  (2.28) will follow from
$$\sup_x\int\Bigl|\int \tD_z p_s(z,x)\tD_zp_t(z,y')\mu(dz)\Bigr|\mu(dy')\le c_{4.11}s^{-1-\eta}t^{-1+\eta}\hbox{ for }0<t\le s\le 2,\eqno(4.33)$$
and
$$\sup_x\int\Bigl|\int \tD_x p_s(x,z)\tD_{y'}p_t(y',z)\mu(dz)\Bigr|\mu(dy')\le c_{4.11}s^{-1-\eta}t^{-1+\eta}\hbox{ for }0<t\le s\le 2.\eqno(4.34)$$
To see (4.1), multiply both sides of  (4.33) by 
$$\int \int |\tD_zp_s(z,x)|\,|\tD_zp_t(z,y)|\mu(dz)\mu(dx).$$
After taking a supremum over $y$, the resulting left-hand side is an upper bound for the left-hand side of (4.1).  For the resulting right-hand side,  use (4.24) to first bound the integral in $x$, uniformly in $z$,  by $c_{4.8}s^{-1}$ and (4.25) to then bound the integral in $z$, uniformly in  $y$ by $c_{4.8}t^{-1}$.  This gives (4.1) and similar reasoning derives (4.2) from (4.34). 

Next use Lemma~4.11 to see that
$$\eqalign{\int&\Bigl|\int \tD_{y'}p_s(y',x)\tD_zp_t(z,y')\mu(dz)\Bigr|\mu(dy')\cr
&\le \int|\tD_{y'}p_s(y',x)|c_{4.11}(t+y'_{m+1})^{-1}\mu(dy')\cr
&\le \int |D^2_{y'_{m+1}}p_s(y',x)|\mu(dy')\le cs^{-2}\le cs^{-1-\eta}t^{-1+\eta}.}\eqno(4.35)$$
In the next to last inequality we have used Lemma~4.5(b) with $q=p=0$. Therefore
the triangle inequality shows that (4.33) (with perhaps a different constant) will follow from
$$\eqalignno{&\int\Bigl|\int (\tD_z p_s(z,x)-\tD_{y'}p_s(y',x))\tD_zp_t(z,y')\mu(dz)\Bigr|\mu(dy')\le c_{4.12}s^{-1-\eta}t^{-1+\eta}\cr
&\phantom{\int\Bigl|\int (\tD_z p_s(z,x)-\tD_{y'}p_s(y',x))\tD_zp_t(z,y')\mu(dz)}\hbox{ for }0<t\le s\le 2.&(4.36)}$$
The analogous reduction for (4.34) is easier as (2.14) with $f\equiv 1$ implies
$$ \int \tD_xp_s(x,y')\tD_{y'}p_t(y',z)\mu(dz)=0.$$
Use this in place of (4.35) and again apply the triangle inequality to see that (4.34) will follow from
$$\eqalignno{&\int\Bigl|\int (\tD_xp_s(x,z)-\tD_xp_s(x,y'))\tD_{y'}p_t(y',z)\mu(dz)\Bigr|\mu(dy')\le c_{4.12}s^{-1-\eta}t^{-1+\eta}\cr
&\phantom{\int\Bigl|\int (\tD_zp_s(x,z)-\tD_xp_s(x,y'))\tD_{y'}p_t(y',z)\mu(dz)}\hbox{ for }0<t\le s\le 2.&(4.37)}$$

Having reduced our problem to (4.36) and (4.37), we consider (4.36) first and take $\eta=1/2$ for the rest of the proof.  The left-hand side of (4.36) is bounded by
$$\eqalignno{\int&\Bigl|\int (z_{m+1}(D_{z_{m+1}}^2p_s(z,x)-D^2_{y'_{m+1}}p_s(z^{(m)},y'_{m+1},x))z_{m+1}D^2_{z_{m+1}}p_t(z,y')\mu(dz)\Bigr|\mu(dy')\cr
&+\int\Bigl|\int (z_{m+1}D^2_{y'_{m+1}}p_s(z^{(m)},y'_{m+1},x)-y'_{m+1}D^2_{y'_{m+1}}p_s(y',x))\cr
&\qq\qq\qq \times z_{m+1}D^2_{z_{m+1}}p_t(z,y')\mu(dz)\Bigr|\mu(dy')\cr
&:= T_{a,1}+T_{a,2}.&(4.38)}$$
Use the Fundamental Theorem of Calculus (recall Proposition~2.8 for the required regularity) to see that
$$\eqalignno{T_{a,1}\le \int\int\int(1_{(z_{m+1}<w<y'_{m+1})}&+1_{(y'_{m+1}<w<z_{m+1})})z_{m+1}|D^3_wp_s(z^{(m)},w,x)|\cr
&\times z_{m+1}|D^2_{z_{m+1}}p_t(z,y')|\,dw\,\mu(dz)\mu(dy').)}$$
Now recall from (3.7) that $p_t(z,y)=p^0_t(z^{(m)}-y^{(m)},z_{m+1},y_{m+1})$.  
First do the \break
$\mu(dy')\,\mu_{m+1}(dz_{m+1})$ integrals and change variables to $y''=z^{(m)}-y'^{(m)}$ in this integral to see
that 
$$\eqalignno{T_{a,1}\le \int\int\int\int\int&(1_{(z_{m+1}<w<y'_{m+1})}+1_{(y'_{m+1}<w<z_{m+1})})\cr
&\times z^2_{m+1}|D^2_{z_{m+1}}p^0_t(y'',z_{m+1},y'_{m+1})|dy''\mu_{m+1}(dz_{m+1})\mu(dy_{m+1})\cr
&\times |D^3_wp_s(z^{(m)},w,x)|dz^{(m)}dw\cr
\le c_{4.7}\int_{\gamma t}^\infty &\int t^{-1/2}w^{3/2}|D^3_wp_s(z^{(m)},w,x)|dz^{(m)}\,w^{b/\gamma-1}dw\cr
&+c_{4.7}\int_0^{\gamma t}\int t^{b/\gamma}|D^3_wp_s(z^{(m)},w,x)|dz^{(m)}dw.&(4.39)}$$
In the last line we have used Lemma~4.7 with $p=2$. Now use Lemma~4.5(d) to bound the first term by $ct^{-1/2}s^{-3/2}$, and use (4.22) to bound the second term by $cs^{-2}\le ct^{-1/2}s^{-3/2}$. We have proved
$$T_{a,1}\le c_1 t^{-1/2}s^{-3/2}.\eqno(4.40)$$

Note that 
$$\eqalignno{T_{a,2}&\le \int\Bigl|\int (z_{m+1}(D^2_{y'_{m+1}}p_s(z^{(m)},y'_{m+1},x)-D^2_{y'_{m+1}}p_s(y',x))\cr
&\qq\qq\qq \times z_{m+1}D^2_{z_{m+1}}p_t(z,y')\mu(dz)\Bigr|\mu(dy')\cr
&\quad+\int\Bigl|\int(z_{m+1}-y'_{m+1})D^2_{y'_{m+1}}p_s(y',x)z_{m+1}D^2_{z_{m+1}}p_t(z,y')\mu(dz)\Bigr|\mu(dy')\cr
:=& T_{a,3}+T_{a,4}.&(4.41)}$$
By Lemma~4.10,
$$\int D^2_{z_{m+1}}p_t(z,y')dz^{(m)}=D^2_{z_{m+1}}q_t(z_{m+1},y'_{m+1}),$$ and so using Lemma~3.3(g) we have
$$\eqalignno{T_{a,4}&=\int\bigl|\int(z_{m+1}-y'_{m+1})z_{m+1}^{b/\gamma}D^2_{z_{m+1}}q_t(z_{m+1},y'_{m+1})dz_{m+1}\Bigr||D^2_{y'_{m+1}}p_s(y',x))|\mu(dy')\cr
&\le c_{3.3}\int|D^2_{y'_{m+1}}p_s(y',x))|\mu(dy')\le c_4s^{-2}\le c_4t^{-1/2}s^{-3/2},&(4.42)}$$
where we have used Lemma~4.5(b) in the last line with $p=q=0$.  

For $T_{a,3}$ use the Fundamental Theorem of Calculus to write
$$\eqalign{T_{a,3}&=\int\Bigl|\int \Bigl[\int_0^1\sum_{j=1}^m(z_j-y'_j)D_jD^2_{y'_{m+1}}p^0_s(y'^{(m)}-x^{(m)}+r(z^{(m)}-y'^{(m)}),y'_{m+1},x_{m+1})\,dr\Bigr]\cr
&\times z^2_{m+1}D^2_{z_{m+1}}p_t^0(z^{(m)}-y'^{(m)},z_{m+1},y'_{m+1})dz^{(m)}\mu_{m+1}(dz_{m+1})\Bigr|dy'^{(m)}\mu_{m+1}(dy'_{m+1}).}$$
Now take the absolute values inside the integrals and summation, do the integral in $r$ last, and for each $r$ carry out the linear change of variables for the other ($2m$-dimensional) Lebesgue integrals: $(u,w)=(z^{(m)}-y'^{(m)},y'^{(m)}-x^{(m)}+r(z^{(m)}-y'^{(m)}))$ (noting that $|dz^{(m)}dy'^{(m)}|\le 2^m|du
\, dw|$ for all $0\le r\le 1$.
This shows
$$\eqalignno{T_{a,3}&\le c\sum_{j=1}^m\int\int\Bigl[\int\int|u_j|z^2_{m+1}|D^2_{z_{m+1}}p^0_t(u,z_{m+1},y'_{m+1})|du\,\mu_{m+1}(dz_{m+1})\Bigr]\cr
&\qquad\times |D_{w_j}D^2_{y'_{m+1}}p^0_s(w,y'_{m+1},x)|\,dw\,\mu_{m+1}(dy'_{m+1})&(4.43)\cr
&\le c\sum_{j=1}^m t^{-1/2}\int\int(t^{3/2}+(y'_{m+1})^{3/2})|D_{w_j}D^2_{y'_{m+1}}p^0_s(w,y'_{m+1},x_{m+1})|\,dw\,\mu_{m+1}(dy'_{m+1}).}$$
For the last inequality we have used Lemma~4.5(b) with $q=1$ and $p=2$.  Now use Lemma 4.5(f) with $p=0$ (for the $t^{3/2}$ term) and then with $p=3/2$ (for the $(y'_{m+1})^{3/2}$ term) to conclude that 
$$T_{a,3}\le c [ts^{-3}+t^{-1/2}s^{-3/2}]\le c_3t^{-1/2}s^{-3/2}\eqno(4.44)$$
Combining (4.40), (4.42) and (4.44) now gives (4.36) (with $\eta=1/2$). 

The left side of (4.37) is at most
$$\eqalignno{\int&\Bigl|\int x_{m+1}(D^2_{x_{m+1}}p_s(x,z)-D^2_{x_{m+1}}p_s(x,z^{(m)},y'_{m+1}))y'_{m+1}D^2_{y'_{m+1}}p_t(y',z)\mu(dz)\Bigr|\mu(dy')\cr
&+\int\Bigl|\int x_{m+1}(D^2_{x_{m+1}}p_s(x,z^{(m)},y'_{m+1})-D^2_{x_{m+1}}p_s(x,y'))y'_{m+1}D^2_{y'_{m+1}}p_t(y',z)\mu(dz)\Bigr|\mu(dy')\cr
&:= T_{b,1}+T_{b,2}.&(4.45)}$$
The Fundamental Theorem of Calculus gives (Lemma~4.4 gives the required regularity)
$$\eqalign{T_{b,1}&\le \int\int\int\int\int (1_{(z_{m+1}<w<y'_{m+1})}+1_{(y'_{m+1}<w<z_{m+1})})x_{m+1}\cr
&\phantom{\int\int\int}\times |D_wD^2_{x_{m+1}}p_s(x,z^{(m)},w)|y'_{m+1}|D^2_{y'_{m+1}}p_t(y',z)|\cr
&\qq \qq dw\, dz^{(m)}\, dy'^{(m)}\mu_{m+1}(dz_{m+1})\mu_{m+1}(dy'_{m+1}).}$$
Re-express $p_t$ in terms of $p^0_t$ and set $y''=y'^{(m)}-z^{(m)}$ to conclude
$$\eqalignno{T_{b,1}&\le \int\int\int (1_{(z_{m+1}<w<y'_{m+1})}+1_{(y'_{m+1}<w<z_{m+1})})\cr
&\phantom{\int\int\int (}\times y'_{m+1}|D^2_{y'_{m+1}}p^0_t(y'',y'_{m+1},z_{m+1}))|\,dy''\mu_{m+1}(dy'_{m+1})\mu_{m+1}(dz_{m+1})\cr
&\phantom{\int\int\int (}\times x_{m+1}|D_wD^2_{x_{m+1}}p_s(x,z^{(m)},w)|dw\,dz^{(m)}\cr
&\le c_{4.7}\int_{\gamma t}^\infty\int t^{-1/2}w^{1/2}w^{b/\gamma-1}x_{m+1}|D_wD^2_{x_{m+1}}p_s(x,z^{(m)},w)|dz^{(m)}dw\cr
&\quad+c_{4.7}\int_0^{\gamma t}\int t^{b/\gamma-1}  x_{m+1}|D_wD^2_{x_{m+1}}p_s(x,z^{(m)},w)|dz^{(m)}dw.&(4.46)}$$
In the last line we used Lemma 4.7 with $p=1$.  Now use (4.18) with $p=1/2$ to bound the first  term by $c_{4.5}c_{4.7}t^{-1/2}s^{-3/2}$.  By Lemma 4.4, the second term in (4.46) is at most
$$\eqalign{c&t^{-1}\int\Bigl[\int_0^{\gamma t}\int t^{b/\gamma} |D_w\hat p_{s/2}(z^{(m)},w,z')|dz^{(m)}dw\Bigr]\, x_{m+1}|D^2_{x_{m+1}}p_{s/2}(x,z')|\mu(dz')\cr
&\le ct^{-1}(t/s)s^{-1}\le cs^{-2}\le ct^{-1/2}s^{-3/2},}$$
where we have used (4.21) then (4.5). (We are applying (4.21) to $\hat p_{s/2}$.) We have shown
$$T_{b,1}\le ct^{-1/2}s^{-3/2}.\eqno(4.47)$$

For $T_{b,2}$, an argument similar to that leading to (4.43) bounds $T_{b,2}$ above by
$$\eqalign{c\sum_{j=1}^m\int\int&\Bigl[\int\int|u_j|y'_{m+1}|D^2_{y'_{m+1}}p^0_t(u,y'_{m+1},z_{m+1})|du\,\mu_{m+1}(dz_{m+1})\Bigr]\cr
&\times x_{m+1}|D_{w_j}D^2_{x_{m+1}}p^0_s(w,x_{m+1},y'_{m+1})|dw\,\mu_{m+1}(dy'_{m+1})\cr
\le c\int \int&c_{4.5}[1+t^{-1/2}(y'_{m+1})^{1/2}]x_{m+1}|D_{w_j}D^2_{x_{m+1}}p^0_s(w,x_{m+1},y'_{m+1})|dw\,\mu_{m+1}(dy'_{m+1}).}$$
In the last line we have used the identity $p^0_t(u,y'_{m+1},z_{m+1})=p_t(0,y'_{m+1},-u,z_{m+1})$ and then Lemma~4.5(c) with $q=1$.  Finally
use (4.19) with $p=0$ and $p=1/2$ to bound the above by $c(s^{-2}+t^{-1/2}s^{-3/2})\le ct^{-1/2}s^{-3/2}$.  Use this and (4.47) in (4.45) to complete the proof of (4.37).\qed

Having obtained (2.28) and Proposition 2.10 for the special case $N_2$ null and $Z\cap C=\{d\}$, we now turn to the general case.  In the rest of this section we work in the general setting of Propositions~2.2 and 2.14.
\ms

\ni{\bf Proof of Proposition 2.14.}  We need to establish (2.28) (thanks to Lemma 3.16), and first do this for the special case when our transition density is $q_t=q_t^{b,\gamma}$, that is $Z\cap C$ empty and $N_2$ a singleton.  Let $p_t$ be the transition density considered above with $m=1$.  Recall from (3.7) that
$$\int p_t(x,z)dz_1=\int p_t(x,z)dx_1=q_t(x_2,z_2).\eqno(4.48)$$
Let $\tD_{y_2}=D_{y_2}$ or $y_2D_{y_2y_2}$.  We claim that we can differentiate through the above integrals and so
$$\int \tD_{x_2}p_t(x,z)dz_1=\tD_{x_2}q_t(x_2,z_2)\hbox{ for almost all }z_2>0 \hbox{ and all }x,\eqno(4.49)$$
and 
$$\int \tD_{x_2}p_t(x,z)dx_1=\tD_{x_2}q_t(x_2,z_2)\hbox{ for all }x_2>0 \hbox{ and }z.\eqno(4.50)$$
Lemma~4.10 implies (4.50).  The proof of (4.49) uses $\int|\tD_{x_2}p_t(x,z)|dz_1<\infty$ for a.a. $z_2>0$ by Lemma~3.11(b), and then proceeds using the Fundamental Theorem of Calculus as in the proof of Proposition~2.8(d) in Section~3.  (The stronger version of (4.49) also holds but this result will suffice.)  

Consider first (4.2) for $q_t$.  Let $0<t\le s\le 2$.  By (4.49) and (4.50), we have for all $x_2,y_2>0$,
$$\eqalignno{\Bigl|\int &\tD_{x_2}q_s(x_2,z_2)\tD_{y_2}q_t(y_2,z_2)\mu_2(dz_2)\Bigr|\cr
&=\Bigl|\int\Bigl[\int \tD_{x_2}p_s(x,z)dx_1\Bigr]\Bigl[\int \tD_{y_2}p_t(y,z)dz_1\Bigr]\mu_2(dz_2)\Bigr|\cr
&\le \int\Bigl|\int \tD_{x_2}p_s(x,z)\tD_{y_2}p_t(y,z)\mu(dz)\Bigr|dx_1.&(4.51)}$$
Similarly, for all $x_2>0$,  
$$\eqalignno{\int\Bigl|&\int \tD_{x_2}q_s(x_2,z_2)\tD_{y'_2}q_t(y'_2,z_2)\mu_2(dz_2)\Bigr|\mu_2(dy'_2)\cr
&=\inf_{x_1}\int\Bigl|\int \Bigl[\int \tD_{x_2}p_s(x,z)dz_1\Bigr]\Bigl[\int \tD_{y'_2}p_t(y',z)dy'_1\Bigr]\mu_2(dz_2)\Bigr|\mu_2(dy'_2)\cr
&\le \inf_{x_1}\int\Bigl|\int \tD_{x_2}p_s(x,z)\tD_{y'_2}p_t(y',z)\mu(dz)\Bigr|\mu(dy').&(4.52)}$$
Integrability issues are handled by Lemma~4.10 and Proposition~4.8.  The infimum in the second line can be omitted as the expression following does not depend
on $x_1$.
 Multiply (4.52) and (4.51) and integrate with respect to $\mu_2(dx_2)$ to see that for any $y_2>0$,
$$\eqalign{\int&\Bigl|\int\tD_{x_2}q_s(x_2,z_2)\tD_{y'_2}q_t(y'_2,z_2)\mu_2(dz_2)\Bigr|\cr
&\qquad\times \Bigr|\int \tD_{x_2}q_s(x_2,z_2)\tD_{y_2}q_t(y_2,z_2)\mu_2(dz_2)\Bigl|\mu_2(dy'_2)\mu_2(dx_2)\cr
&\le \int\Bigl|\int \tD_{x_2}p_s(x,z)\tD_{y'_2}p_t(y',z)\mu(dz)\Bigr|\Bigl|\int \tD_{x_2}p_s(x,z)\tD_{y_2}p_t(y,z)\mu(dz)\Bigr|\mu(dy')\mu(dx)\cr
&\le cs^{-2-\eta}t^{-2+\eta},}$$
the last by Proposition~4.12.  This gives (4.2) for $q_t$.  A similar argument works for (4.1).

Next we consider (2.28) in the general case.  Write $x=((x_{(i)})_{i\in Z\cap C},( x_j)_{j\in N_2})$ so that (from (2.5)) 
$$p_t(x,y)=\prod_{j\in Z\cap C}p^j_t(x_{(j)},y_{(j)})\prod_{j\in N_2}q_t^j(x_j,y_j).\eqno(4.53)$$
For $j\in (Z\cap C)\cup N_2$, let $x_{\hat j}$ denote $x$ but with $x_{(j)}$ (if $j\in Z\cap C$) or $x_j$ (if $j\in N_2$) omitted, $\mu_{\hat j}=\prod_{i\neq j}\mu_i$, and let $p^{\hat j}_t(x_{\hat j},y_{\hat j})$ denote the above product of transition densities but with the $j$th factor (which may be a $p^j_t$ or a $q^j_t$) omitted.  Consider (4.2) and let $\tD_x\equiv\tD_{x_{(j)}}$ be one of the differential operators in Notation~2.9 acting on the variable $j'\in \{j\}\cup R_j$ for some $j\in Z\cap C$.  (The case $j'=j\in N_2$ is considered below.)  In this case (4.53) shows that the left-hand side of (4.2) equals
$$\eqalign{\sup_y&\int\Bigl[\int\Bigl|\int \tD_{x_{(j)}}p^j_s(x_{(j)},z_{(j)})\tD_{y'_{(j)}}p_t^j(y'_{(j)},z_{(j)})p_s^{\hat j}(x_{\hat j},z_{\hat j})p_t^{\hat j}(y'_{\hat j},z_{\hat j})\cr
&\phantom{\int\Bigl[\int\Bigl|\int \tD_{x_{(j)}}p^j_s(x_{(j)},z_{(j)})\tD_{y'_{(j)}}}\mu_j(dz_{(j)})\mu_{\hat j}(dz_{\hat j})\Bigr|\mu_j(dy'_{(j)})\mu_{\hat j}(dy'_{\hat j})\cr
&\quad \times\Bigl|\int \tD_{x_{(j)}}p^j_s(x_{(j)},z_{(j)})\tD_{y_{(j)}}p_t^j(y_{(j)},z_{(j)})p_s^{\hat j}(x_{\hat j},z_{\hat j})p_t^{\hat j}(y_{\hat j},z_{\hat j})\cr
&\phantom{\int\Bigl[\int\Bigl|\int \tD_{x_{(j)}}p^j_s(x_{(j)},z_{(j)})\tD_{y'_{(j)}}}\mu_j(dz_{(j)})\mu_{\hat j}(dz_{\hat j})\Bigr|\Bigr]\mu_j(dx_{(j)})\mu_{\hat j}(d\hat x_j).\cr}$$
Take the absolute values inside the two $\mu_{\hat j}(dz_{\hat j})$ integrals (giving an upper bound) and pull the 
$p^{\hat j}$ terms out of the $\mu_j(dz_{(j)})$ integrals.  
Now we can integrate the $p^{\hat j}$ integrals using (3.8) by first integrating over $y'_{\hat j}$, then the first $z_{\hat j}$ integral, then the $x_{\hat j}$ integral, and finally the second $z_{\hat j}$ integral.  
This shows that the left-hand side of (4.2) is at most
$$\eqalign{\sup_{y_{(j)}}\int\int&\Bigl|\int \tD_{x_{(j)}}p^j_s(x_{(j)},z_{(j)})\tD_{y'_{(j)}}p_t^j(y'_{(j)},z_{(j)})\mu_j(dz_{(j)})\Bigr|\mu_j(dy'_{(j)})\cr
&\times \Bigl|\int \tD_{x_{(j)}}p^j_s(x_{(j)},z_{(j)})\tD_{y_{(j)}}p_t^j(y_{(j)},z_{(j)})\mu_j(dz_{(j)})\Bigr|\mu_{(j)}(dx_{(j)})\cr
\le & cs^{-2-\eta}t^{-2+\eta}.}$$
In the last line we used (4.2) for $p^j$ (i.e., Proposition~4.12). For $j'=j\in N_2$ we would use (4.2) for $q_t$, which was established above.  This completes the proof of (4.2) and the proof for (4.1) is similar.  \qed

\ms 

\ni{\bf Proof of Proposition 2.10.}   By Proposition~4.8 the required result holds  for each $p^j_t$ factor and then using (4.49) and (4.50), one easily verifies it for $q_t$ (as was done implicitly in the previous proof).  The general case now follows easily from the product structure (4.53) and a short calculation which is much simpler than that given above.\qed

\subsec {5. Proof of Proposition 2.3}

Assume $\P$ is a solution of $M(\tilde\sA,\nu)$ where $d\nu=\rho \,d\mu$ is as in Proposition~2.3, and assume (2.7) throughout this section. Without loss of generality
 we may
work on a probability space carrying a $d$-dimensional Brownian motion $B$ and realize $\P$ as the law of a solution $X$ of 
$$\eqalign{X^i_t&=X^i_0+\int_0^t\sqrt{2\tilde \gamma_i(X_s)X^i_s}\,dB^i_s+\int_0^t\tilde b_i(X_s)\,ds,\quad i\notin N_1,\cr
X^j_t&=X_0^j+\int_0^t\sqrt{2\tilde \gamma_j(X_s)X^i_s}\,dB^j_s+\int_0^t\tilde b_j(X_s)\,ds,\quad j\in R_i, i\in Z\cap C.}$$
Set $[s]_n=(([ns]-1)/n)\lor 0$, define $\tilde \gamma_j(X_s)\equiv\gamma_j^0$, $\tilde b_j(X_s)\equiv b_j^0$ if $s<0$, and consider the unique solution $X^n$ 
to 
$$\eqalignno{X^{n,i}_t&=X^i_0+\int_0^t\sqrt{2\tilde \gamma_i(X_{[s]_n})X^{n,i}_s}\,dB^i_s+\int_0^t\tilde b_i(X_{[s]_n})\,ds,\quad i\notin N_1,&(5.1)\cr
X^{n,j}_t&=X_0^j+\int_0^t\sqrt{2\tilde \gamma_j(X_{[s]_n})X^{n,i}_s}\,dB^j_s+\int_0^t\tilde b_j(X_{[s]_n})\,ds,\quad j\in R_i, i\in Z\cap C.}$$
Note that
$$\eqalignno{&\hbox{ for $k\ge 0$, on $[{k\over n},{k+1\over n}]$ and conditional on $\sF_{k/n}$, $X^n$ has generator $\sA^0$ but with $\gamma^0$, $b^0$}\cr
&\hbox{replaced with the (random) $\gamma^{k}\equiv \tilde\gamma(X_{(k-1)/n})$, $b^k\equiv \tilde b(X_{(k-1)/n})$.}&(5.2)}$$
We see in particular that pathwise uniqueness of $X^n$ follows from the classical Yamada-Watanabe theorem.   An easy stochastic calculus argument, using Burkholder's inequalities and the boundedness
of $\tilde \gamma$, $\tilde b$ and $X_0$, shows that 
$$\E((X^{n,i}_t)^p)\le c_p(1+t^p)\hbox{  for all }t\ge 0\hbox{ and }p\in \N.\eqno(5.3)$$
Here $c_p$ may depend on the aforementioned bounds and is independent of $n$ (although we will
not need the latter).  

By making only minor modifications in the proof of Lemma 5.1 in [ABBP] we have:

\proclaim Lemma 5.1.  For any $T>0$, $\sup_{t\le T}\Vert X^n_t-X_t\Vert\to 0$ in probability as $n\to\infty$.

For $k\in\Z_+$, let
$$\mu_k(dx)=\prod_{i\in Z\cap C}\Bigl(\prod_{j\in R_i}dx_j\Bigr)x_i^{b_i^k/\gamma_i^k-1}dx_i\times\prod_{j\in N_2}x_j^{b_j^k/\gamma_j^k-1}dx_j,$$
and let $p^k_t(x,y)$ denote the (random) transition density with respect to $\mu_k$ of the diffusion described in (5.2) operating on the interval $[k/n,(k+1)/n]$.  Proposition~2.8(b) with $n=0$ implies that 
$$p^k_t(x,y)\le c_{5.1}t^{-c_{5.1}}\prod_{j\in N_2}(t^{1/2}+y_j^{1/2})\le c_{5.1}t^{-c_{5.1}}\prod_{j\in N_2}(1+y_j^{1/2}), \hbox{ for }x,y\in S^0, 0<t\le 1,\eqno(5.4)$$
where as usual $c_{5.1}$ may depend on $M_0, d$ but not on $k$.  We are also using (2.7) here to bound $b^k_i/\gamma_i^k$.

Let $S^n_\lam f=\E\Bigl(\int_0^\infty e^{-\lam t}f(X^n_t)\,dt\Bigr)$ and define $\Vert S_\lam^n\Vert=\sup\{|S^n_\lam f|:\Vert f\Vert_2\le 1\}$, where as usual the $L^2$ norm refers to the fixed measure $\mu$.  

\proclaim Lemma 5.2. If (2.7) holds, then $\Vert S^n_\lam\Vert<\infty$ for all $\lam>0$, $n\in\N$.  

\proof It suffices to consider $|S^n_\lam f|$ for non-negative $f\in L^2(\mu)$.  Let 
$$\delta=\sup_{i\notin N_1, x}\Bigl|{\tilde b^i(x)\over \tilde\gamma^i(x)}-{b_i^0\over \gamma_i^0}\Bigr|.\eqno(5.5)$$
A bit of algebra using (2.7) shows that
$$\delta\le {\vareps_0M_0^2\over M_0^{-1}-\vareps_0}\le 2\vareps_0M_0^3.\eqno(5.6)$$
Let $\E^k_x$ denote expectation starting at $x$ with respect to the law of the diffusion with (random) transition density $p^k$, and let $R^k_\lam$ and $r^k_\lam(x,y)$ denote the corresponding resolvent and resolvent density with respect to $\mu_k$.  Then (5.2) shows that
$$\eqalignno{S_\lam^nf&=\sum_{k=0}^\infty e^{-\lam k/n}\E\Bigl(\E^k_{X^n_{k/n}}\Bigl(\int_0^{1/n}e^{-\lam t}f(X_t)\,dt\Bigr)\Bigr)\cr
&\le\int R_\lam^0f(x)\rho(x)d\mu(x)+\sum_{k=1}^\infty e^{-\lam k/n}\E(\E(R_\lam^k f (X^n_{k/n})|\sF_{(k-1)/n}))\cr
&\le \Vert R^0_\lam f\Vert_2\Vert\rho\Vert_2+\sum_1^\infty e^{-\lam k/n}\E(\int R_\lam ^kf(x)p_{1/n}^{k-1}(X^n_{(k-1)/n},x)\mu_{k-1}(dx))\cr
&\le \lam^{-1}\Vert \rho\Vert_2\Vert f\Vert_2+\sum_1^\infty e^{-\lam k/n}\E(I^n_k),&(5.7)}$$
where
$$I^n_k=\int\int r_\lam^k(x,y)p^{k-1}_{1/n}(X^n_{(k-1)/n},x)\mu_{k-1}(dx)f(y)\prod_{i\notin N_1}(y_i^\delta+y_i^{-\delta})\mu(dy)$$
and we have used Corollary~2.12(a) to see $\Vert R^0_\lam f\Vert_2\le \lam^{-1}\Vert f\Vert_2$.  

As usual we suppress dependence on $d$ and $M_0$ in our constants (which may change from line to line) but will record dependence on $n$.  Use (5.4) and (5.5) to see that 
$$\eqalignno{\int &r_\lam^k(x,y)p^{k-1}_{1/n}(X^n_{(k-1)/n},x)\mu_{k-1}(dx)&(5.8)\cr
&\le c_n\int r_\lam^k(x,y)\prod_{j\in N_2}(1+\sqrt{x_j})\prod_{i\notin N_1}(x_i^{2\delta}+x_i^{-2\delta})\mu_k(dx)\cr
&\le c_n \int r_\lam^k(x,y)\prod_{j\notin N_1}(x_i^{-2\delta}+x_i^{1/2+2\delta})\mu_k(dx).}$$
Let $\hat\E^k_x$ denote expectation with respect to the  diffusion with transition kernel $\hat {p}^k_t(x,y)=p^k_t(y,x)$ (as in Proposition~2.8(a)) and $\hat r^k_\lam(x,y)=r^k_\lam(y,x)$ be the associated resolvent density with respect to $\mu_k$.  Then
(5.8) is bounded by
$$\eqalign{c_n&\int_0^\infty e^{-\lam s}\hat\E^k_y\Bigl(\prod_{i\notin N_1}((X_s^i)^{-2\delta}+(X_s^i)^{1/2+2\delta})\Bigr)\,ds\cr
&=c_n \int_0^\infty e^{-\lam s}\prod_{i\notin N_1}\hat\E^k_y\Bigl((X_s^i)^{-2\delta}+(X_s^i)^{1/2+2\delta})\Bigr)\,ds\hbox{  (by (2.5))}\cr
&\le c_n\int_0^\infty e^{-\lam s}\prod_{i\notin N_1}\Bigl(s^{-2\delta}+[\hat\E^k_y(X_s^i)]^{1/2+2\delta}\Bigl)\,ds.}$$
In the next to last line we have used the conditional independence of $\{X^i:i\notin N_1\}$ (recall (2.5)). In the last line we have used (5.6) and (2.7) to see that $1/2+2\delta\le 1$, and also used Lemma~3.3(d).  We can apply this last result because for all $i\notin N_1$ and all $k\in \Z_+$, 
$$\eqalignno{b_i^k/\gamma_i^k-2\delta&\ge {M_0^{-1}-\vareps_0\over M_0+\vareps_0}-2\delta&(5.9)\cr
&\ge (4M_0^2)^{-1}-4\vareps_0M_0^3\ge (8M_0^2)^{-1},\cr}$$
thanks to (5.6) and (2.7).  For $i\notin N_1$ we have $\E^k_y(X^i_s)\le y_i+M_0s$, and by (2.7) and (5.6) we have $2d\delta\le 1/12$.  Therefore if $f(s)=s^{-2\delta}+(M_0s)^{1/2+2\delta}$ (clearly $f\ge 1$) we may now bound (5.8) by
$$\eqalign{c_n&\int_0^\infty e^{-\lam s}\prod_{i\notin N_1}(f(s)+y_i^{1/2+2\delta})\,ds\cr
&\le c_n\int_0^\infty e^{-\lam s}f(s)^{|N_1^c|}\prod_{i\notin N_1}(1+(y_i^{1/2+2\delta}/f(s)))\,ds\cr
&\le c_n\int _0^\infty e^{-\lam s}f(s)^{|N_1^c|}\,ds \prod_{i\notin N_1}(1+y_i^{1/2+2\delta})
\le c_{n,\lam}\prod_{i\notin N_1}(1+y_i^{1/2+2\delta}).\cr}$$
In the definition of $I_k^n$ use H\"older's inequality and then this bound  on (5.8) on one of the resulting squared factors  to see that
$$\eqalignno{I_k^n&\le c\Vert f\Vert_2\Bigl[\int\Bigl[\int r^k_\lam(x,y)p_{1/n}^{k-1}(X^n_{(k-1)/n},x)\mu_{k-1}(dx)\Bigr]^2\prod_{i\notin N_1}(y_i^{2\delta}+y_i^{-2\delta})\mu(dy)\Bigr]^{1/2}\cr
&\le c_{n,\lam} \Vert f\Vert_2 \Bigl[\int\int r_\lam^k(x,y)p^{k-1}_{1/n}(X^n_{(k-1)/n},x)\prod_{i\notin N_1}(1+y_i^{1/2+2\delta})\cr
&\qq\qq\qq\times\prod_{i\notin N_1}(y_i^{2\delta}+y_i^{-2\delta})\mu(dy)\mu_{k-1}(dx)\Bigr]^{1/2}\cr
&\le  c_{n,\lam} \Vert f\Vert_2 \Bigl[\int p^{k-1}_{1/n}(X^n_{(k-1)/n},x)\int r_\lam^k(x,y)\prod_{i\notin N_1}(1+y_i^{1/2+2\delta})\cr
&\qq\qq\times\prod_{i\notin N_1}(y_i^{3\delta}+y_i^{-3\delta})\mu_k(dy)\mu_{k-1}(dx)\Bigr]^{1/2}\cr
&\le  c_{n,\lam} \Vert f\Vert_2 \Bigl[\int p^{k-1}_{1/n}(X^n_{(k-1)/n},x) J_{k,\lam}(x)\mu_{k-1}(dx)\Bigr]^{1/2}.&(5.10)}$$
Here
$$\eqalign{J_{k,\lam}(x)&=\int r_\lam^k(x,y)\prod_{i\notin N_1}(y_i^{1/2+5\delta}+y_i^{-3\delta})\mu_k(dy)\cr
&=\int_0^\infty e^{-\lam s}\E^k_x\Bigl(\prod_{i\notin N_1}((X_s^i)^{1/2+5\delta}+(X_s^i)^{-3\delta})\Bigr)\,ds\cr
&\le c\int_0^\infty e^{-\lam s}\prod_{i\notin N_1}\E^k_x\Bigl(X_s^i+(X_s^i)^{-3\delta}\Bigr)\,ds\cr
&\le c\int_0^\infty e^{-\lam s} \prod_{i\notin N_1}[x_i+M_0s+s^{-3\delta}]\,ds.\cr}$$
In the next to last line we have again used the independence of $X^i,i\notin N_1$ under $\E^k_x$, and the bound $1/2+5\delta\le 1$ which follows from (5.6) and (2.7).  In the last line we have again used 
Lemma~3.3(d) whose applicability can again be checked as in (5.9).
An elementary calculation on the above bound, again using (5.6) and (2.7) to see that $3\delta d\le 1/8$, now shows that 
$$J_{k,\lam}(x)\le c_\lam \prod_{i\notin N_1}(1+x_i),$$
and therefore by (5.10),
$$\eqalign{I_k^n&\le c_{n,\lam}\Vert f\Vert_2\Bigl[\int p_{1/n}^{k-1}(X^n_{(k-1)/n},x)\prod_{i\notin N_1}(1+x_i)\mu_{k-1}(dx)\Bigr]^{1/2}\cr
&\le c_{n,\lam}\Vert f\Vert_2\prod_{i\notin N_1}\Bigl(1+X^{n,i}_{(k-1)/n}+{M_0\over n}\Bigr).}$$
Take expectations in the above and use some elementary inequalities to conclude
that
$$\eqalign{\E(I_k^n)&\le c_{n,\lam}\Vert f\Vert_2(1+M_0/n)^d\E\Bigl(\prod_{i\notin N_1}(1+X^{n,i}_{(k-1)/n})\Bigr)\cr
&\le c_{n,\lam} \Vert f\Vert_2\sum_{i\notin N_1}\E((1+X^{n,i}_{(k-1)/n})^d)\cr
&\le c_{n,\lam} \Vert f\Vert_2(1+(k/n)^d),}$$
the last by (5.3).  Put this bound into (5.7) to complete the proof.  \qed

\proclaim Remark 5.3. {\rm  The above argument is considerably longer than its counterpart 
(Lemma 5.3) in [ABBP].  This is in part due to the non-compact state
space in the above leading to the unboundedness of the densities on
this state space (recall the bound (5.4)).  More significantly, it is also because the
argument in [ABBP] is incomplete.  In (5.4) of [ABBP] the norm on $f$
actually depends on $k$ and is not the norm on the canonical $L^2$
space.  The argument above, however, will also give a correct proof of
Lemma 5.3 of [ABBP]--in fact the compact state space there leads to
considerable simplification.
}

\ni{\bf Proof of Proposition 2.3.} This now proceeds by making only minor changes 
in the proof of Proposition 2.3 in Section 5 of [ABBP].  One uses the
above Lemmas 5.1, 5.2 and Proposition 2.2. We only point out the
(trivial) changes required.  For $f\in C^2_b(S^0)$ one uses It\^o's
Lemma and (5.1) to obtain the semimartingale decomposition of
$f(X^n_t)$.  The local martingale part is a martingale as in the proof
of Theorem~2.1 in Section~2 (use (5.3)).  Corollary~2.12(a) is used, 
instead of the eigenfunction
expansion in [ABBP], to conclude that the constant coefficient
resolvent $R_\lam$ has bound $\lam^{-1}$ as an operator on $L^2$.  The
rest of the proof proceeds as in [ABBP] where the bound $\vareps_0\le
(2K(M_0))^{-1}$ is used to get the final bound, first on $S^n_\lam
(|f|)$, and then on $S_\lam (|f|)$ by Fatou's lemma.\qed

\subsec{6. Proof of Proposition 2.4}  

Let $(\P^x,X_t)$ ($x\in S^0$) be as in the statement of Proposition~2.4.
Throughout this section, for any Borel set $A$ we let $T_A=T_A(X)=\inf\{t: X_t\in A\}$
and $\tau_A=\tau_A(X)=\inf\{t: X_t\notin A\}$, be the first entrance and exit times,
respectively, and let $|A|$ denote the Lebesgue measure of $A$. We
say a function $h$ is harmonic in $D=B(x,r)\cap S$ if  $h$ is bounded on 
$\ol D$ and $h(X_{t\land \tau_D})$ is a right continuous martingale with respect to $\P_x$ for each $x$.

The key step in the proof of Proposition 2.4 is the following.

\proclaim Proposition 6.1. Let $z\in S$. There exist positive constants $r, c_{6.1}$ and $\al$,
depending on $z$, such that if $h$ is harmonic
in $B(z,r)\cap S$, then 
$$|h(x)-h(z)|\leq c_{6.1}\Big(\frac{|x-z|}{r}\Big)^\al \Big(\sup_{B(z,r)\cap S} |h|
\Big), \qq x\in B(z,r/2)\cap S. \eqno (6.1)$$

\proof By relabeling the axes we may assume that $S^0=\{x\in\R^d:x_i\ge 0 \hbox{ for }i>J_0\}$.  If $z$ is in the interior of $S^0$, the result is easy, because the generator is locally uniformly elliptic, and follows by
the  first paragraph of the proof of Theorem 6.4 of [ABBP].
So suppose $z\in \del S^0$. Then $J_0<d$ and we may assume, again by reordering the axes, that there is a $K\in\{J_0,\dots,d-1\}$ so that $z_i=0$ for all $i>K$and $z_i>0$ if $J_0<i\le K$.  Assume (set $\min\emptyset=1$)
$$0<r<\min_{J_0<i\le K}{z_i\over 2}.\eqno(6.2)$$
Since our result only depends on the values of 
of $h$ in $B(x,r)\cap S$, we may change the diffusion and drift coefficients
of the generator of $X$ outside $\ol{B(z,r)}\cap S$ as we please.
By changing the coefficients in this way and again relabeling the axes if necessary, we may suppose 
that our
generator is 
$$\eqalignno{\tilde\sA f(x)&= \sum_{i=1}^J \sigma_i(x) f_{ii}(x)+\sum_{i=J+1}^K
\sigma_i(x)x_{a(i)}f_{ii}(x)\cr
&\qq +\sum_{i=K+1}^d \sigma_i(x) x_if_{ii}(x)+\sum_{i=1}^d b_i(x)f_i(x),&(6.3)\cr}$$
where $J\leq K$, $a(i)\in \{K+1, \ldots, d\}$ for $i=J+1, \ldots, K$,
each $\sigma_i$ is continuous and bounded above and below by positive
constants, each $b_i$ is continuous and bounded, and each $b_i$
for $i>K$ is bounded below by a positive constant. We have extended our coefficients to the possibly larger space $S^1=\{x\in\R^d:x_i\ge 0\hbox { for all }i\ge K\}$ as this is the natural state space for $\tilde\sA$.  As $B(z,r)\cap S^0=B(z,r)\cap S^1$ (by (6.2)) this will not affect the harmonic functions we are dealing with.  
For $0\le \delta<1$ let
$$\eqalign{ 
Q_n(\delta)&=\prod_{i=1}^J [z_i-2^{-n/2},z_i+2^{-n/2}] \times
\prod_{i=J+1}^K [z_i-2^{-n},z_i+2^{-n}] \times \prod_{i=K+1}^d [\delta 2^{-n},2^{-n}]\cr
R_n(\delta)&=\prod_{i=1}^J [z_i-\tfrac32 \cdot 2^{-n/2},z_i+\tfrac32 \cdot 2^{-n/2}] \times
\prod_{i=J+1}^K [z_i-\tfrac32 \cdot 2^{-n},z_i+\tfrac32 \cdot 2^{-n}] \cr
&\qq\qq  \times \prod_{i=K+1}^d [\delta 2^{-n},2^{-n}].\cr}$$
Take $n\geq 1$ large enough so that $Q_n(0)\subset B(z,r/2)\cap S^1$. 

We will first show there exist $c_2, \delta>0$ independent of $n$ such that
$$\P_x(T_{R_{n+1}(\delta)}<\tau_{Q_n(0)})\geq c_2, \qq x\in Q_{n+1}(0).
\eqno (6.4)$$
We may assume there exist independent one-dimensional Brownian motions $B^i_t$
such that
$$X_t^i=X_0^i+M_t^i+A_t^i,$$
where $dA_t^i=b_i(X_t)\, dt$ and 
$$\eqalign{dM_t^i&=(2\sigma_i(X_t))^{1/2}\, dB_t^i, \qq  i\leq J,\cr 
dM_t^i&=(2\sigma_i(X_t)X_t^{a(i)})^{1/2}\, dB_t^i, \qq J+1\leq i\leq K, \cr  
dM_t^i&=(2\sigma_i(X_t)X_t^{i})^{1/2}\, dB_t^i, \qq K+1\leq i\leq d. \cr}$$  
Since the $b_i$ and $\sigma_i$ are bounded, there exists $t_0$  small such that for all $x\in Q_{n+1}(0)$ 
$$\eqalign{\P_x(\sup_{s\leq t_02^{-n}}|M_s^i|> \tfrac14 \cdot 2^{-n/2})&\leq \tfrac1{8d}, \qq i\leq J, \cr
\P_x(\sup_{s\leq t_02^{-n}}|M_s^i|> \tfrac14 \cdot 2^{-n})&\leq \tfrac1{8d}, \qq J+1\leq i\leq d, \cr}$$
and
$$\sup_{s\leq t_02^{-n}}|A_s^i|\leq \tfrac14 \cdot 2^{-n}, \qq 1\leq i\leq d.$$
The first and last bounds are trivial, and the second inequality is easily proved by first noting
$$\sup_{x\in Q_{n+1}(0)}\E_x\Bigl(\sum_{i=1}^d\int_0^{t_02^{-n}}X_s^i\,ds\Bigr)\le c 2^{-2n},$$
and then using the Dubins-Schwarz Theorem and Markov's inequality. 
Hence
$$\sup_{x\in Q_{n+1}(0)}\P_x(\tau_{R_{n+1}(0)}<t_02^{-n})\leq \tfrac14.\eqno(6.5)$$

By Lemma 6.2 of [ABBP] there exists $\delta$ such that if $U$ is uniformly distributed
on $[t_0/2,t_0]$, then 
$$\sup_{x\in S^1}\P(X_U^i\leq \delta)\leq 1/4d, \qq K+1\leq i\leq d.$$
Scaling shows that
$$\sup_{x\in S^1}\P_x(X^i_{U2^{-n}}\leq \delta 2^{-n})\leq 1/4d, \qq K+1\leq i\leq d.$$
Therefore by (6.5), for any $x\in Q_{n+1}(0)$ with $\P^x$-probability at least $1/2$, $X_{U2^{-n}}\in R_{n+1}(\delta)$.
Since $R_{n+1}(0)\subset Q_n(0)$ and $U2^{-n}\leq t_02^{-n}$, this and (6.5) proves (6.4).

Take $\delta$ even smaller if necessary so that $|Q_n(0)-Q_n(\delta)|
<\frac14 |Q_n(0)|$ and $\delta\le \sigma_i\le \delta^{-1}$, $|b_i|\le \delta^{-1}$ for all $i$.
Next we show that if $G\subset Q_n(0)$ and $|G|\geq |Q_n(0)|/3$, then there is a $c_3(\delta)>0$, independent of $n$,  so that
$$\P_x(T_G<\tau_{Q_n(\delta/2)})\geq c_3, \qq x\in R_{n+1}(\delta). \eqno (6.6)$$
Let
$$Y_t^i=\cases{2^{n}X^i_{2^{-n}t},&\qq $i>J$,\cr
2^{n/2}X^i_{2^{-n}t},&\qq $i\leq J$.\cr}$$
It is straightforward (cf.~[ABBP], Proof of Theorem 6.4) to see that for $t\le \tau_{Q_0(\delta/2)}$, $Y_0\in Q_0(\delta/2)$,
$Y_t$ solves
$$dY_t^i =\wh \sigma_i(Y_t)\, d\wh B_t^i+\wh b_i(Y_t)\, dt,$$
where the  $\wh B^i$ are independent one-dimensional Brownian motions, the $\wh b_i$
are bounded above, and the $\wh \sigma_i$ are bounded above and below
by positive constants depending on $\delta$ but not $n$. (6.6) now follows by Proposition 6.1 of [ABBP]
(with a minor change to account for the fact that $|G\cap Q_n(\delta/2)|/|Q_n(\delta/2)|$ 
is greater than $\frac1{12}$ rather than  $\frac12$).

Combining (6.4) and (6.6) and using the strong Markov property, we see that if $c_4=c_2c_3$, then
$$\P_x(T_G<\tau_{Q_n(0)})\geq c_4>0, \qq x\in Q_{n+1}(0).$$ 
Suppose $h$ is harmonic on $Q_{n_0}$ for
some $n_0$. Our conclusion will follow by setting  $\alpha\allowbreak=\log(1/\rho)/\log 2$ if we show there exists $\rho<1$
such that
$$\hbox{$\Osc_{Q_{n+1}(0)} h\leq \rho \Osc_{Q_n(0)} h$}, \qq n\geq n_0,\eqno (6.7)$$
where $\Osc_A h=\sup_A h-\inf_A h$. 
Take $n\geq n_0$ and by looking at $c_{5}h+c_{6}$, we may suppose $\sup_{Q_n(0)} h=1$ and $\inf_{Q_n(0)} h=0$. By looking at $1-h$ if necessary, we may
suppose $|G|\geq \frac12 |Q_n(0)|$, where $G=\{x\in Q_n(0): h(x)\geq 1/2\}$.
By Doob's optional stopping theorem
$$h(x)\geq \E_x[h(X_{T_G}); T_G<\tau_{Q_n(0)}]
\geq \tfrac12 \P_x(T_G<\tau_{Q_n(0)})\geq c_4/2, \qq x\in Q_{n+1}(0).$$
Hence $\Osc_{Q_{n+1}(0)} h\leq 1-c_4/2$, and (6.7) follows with $\rho=1-c_4/2$. \qed

\longproof{of Proposition 2.4}
We can now proceed as in the proof of Theorem 6.4 of
[ABBP].  To obtain the analogue of (6.14) in [ABBP], we note from (2.2) that if $x\in\partial S^0$, at least one coordinate can bounded below by a squared Bessel process with positive drift starting at zero.  

\ms
\ni{\bf Remark 6.2.} Essentially the same argument shows that if for each $x\in S$, $\P^x$ is  solution of  MP$(\sA,\delta_x)$ as in Theorem~1.4 (it will be Borel strong Markov by Theorem~1.4), then the resolvent $S_\lam$ maps bounded Borel measurable functions to continuous functions.  After localizing the problem, one is left with a generator in the same form as (6.3) and so the proof proceeds as above. 
\qed

\subsec{7.  Proofs of Lemmas 4.5, 4.6 and 4.7}

We work in the setting and with the notation from Sections~3 and 4.
Recall, in particular, the Poisson random variables $N_\rho(t)$ from Lemma~3.4.

\proclaim Lemma 7.1. There is a $c_{7.1}$ such that for all $0\le q\le 2$, $1\le j\le m$, $y\in S_m$, $0<t$,  and $z'=z_{m+1}/\gamma t>0$: 

\ni(a) If for $x\ge 0$ and $n\in\Z_+$,
$$\psi_1(z',n)=(z'+1)^{q/2-1}\Big[1_{(n\le 1)}+1_{(n=1)}{z'}^{-1}+1_{(n\ge 2)}{z'}^{-2}\Big(n+\Big(n-{z'\over 2}\Big)^2\Big)\Big],$$
and
$$\eqalignno{\psi_2(z',x,n)=&1_{(n\le 2)}(1+{z'}^{-n})(1+{z'}^{q/2}+x^{q/2})\cr
&+1_{(n\ge 3)}{z'}^{-3}(|n-z'|^3+3n|n-z'|+n)({z'}^{q/2}+n^{q/2}+x^{q/2}),}$$
then
$$\eqalignno{\int &|y_j-z_j|^q|D^3_{z_{m+1}}p_t(z,y)|dz^{(m)}&(7.1)\cr
&\le c_{7.1}t^{q-3}\liminf_{\delta\to 0}\E_{z_{m+1}}(q_\delta(y_{m+1},X_t^{(m+1)})[\psi_1(z',N_{1/2}(t))+\psi_2(z',X^{(m+1)}_t/t,N_0(t))]).}$$

\ni(b) If for $x,n$ as in (a),
$$\eqalignno{\phi(z',x,n)&=1_{(n\le 1)}(1+x^{q/2}+(z')^{q/2})+1_{(n=1)}(z')^{-1}(1+(z')^{q/2}+x^{q/2})\cr
\phantom{}&+1_{(n\ge 2)}(z')^{-2}(n+(n-z')^2)((z')^{q/2}+n^{q/2}+x^{q/2}),}$$
then
$$\eqalignno{\int &|y_j-z_j|^q|D^2_{z_{m+1}}p_t(z,y)|dz^{(m)}\cr
&\le c_{7.1}t^{q-2}\liminf_{\delta\to 0}\E_{z_{m+1}}(q_\delta(y_{m+1},X_t^{(m+1)})\phi(z',\Xm_t/t,N_0(t))).&(7.2)}$$

\ni(c) $$\eqalignno{\int &|D_{z_{m+1}}p_t(z,y)|dz^{(m)}\cr
&\phantom{}\le c_{7.1}t^{-1}\liminf_{\delta\to 0}\E_{z_{m+1}}\Bigl(q_\delta(y_{m+1},X_t^{(m+1)})\Bigl((1+z')^{-1}+{|N_0-z'|\over z'}\Bigr)\Bigr).&(7.3)}$$
In addition for all $z\in S_m$,
$$\eqalignno{\int &|D_{z_{m+1}}p_t(z,y)|dy^{(m)}\cr
&\phantom{}\le c_{7.1}t^{-1}\liminf_{\delta\to 0}\E_{z_{m+1}}\Bigl(q_\delta(y_{m+1},X_t^{(m+1)})\Bigl((1+z')^{-1}+{|N_0-z'|\over z'}\Bigr)\Bigr).&(7.4)}$$
\proof The proof of (a) is lengthy and the reader may want to first take a look at the simpler proof of (c) given in Section~8.\hfil\break
\ni(a) By Lemma~3.11(d), Fatou's lemma and symmetry we have
$$\eqalignno{\int &|y_j-z_j|^q|D^3_{z_{m+1}}p_t(z,y)|dz^{(m)}\cr
&\le \liminf_{\delta\to0}\int |y_j-z_j|^q\Bigl|\E_{z_{m+1}}\Bigl(\Delta^3_t G^\delta_{t,z^{(m)},y}(X,\nu^1,\nu^2,\nu^3)[1_{(\nu^i_t=0\hbox{ for } i=1,2,3)}\cr
&\quad+31_{(\nu^1_t>0,\nu^2_t=\nu^3_t=0)}
+31_{(\nu^1_t>0,\nu^2_t>0,\nu^3_t=0)}+1_{(\nu^i_t>0\hbox{ for }i=1,2,3)}]\cr
&\phantom{\le \liminf_{\delta\to0}\int |y_j-z_j|^q|\Bigl|\E_{z_{m+1}}\Bigl(\Delta^3_t G^\delta_{t,z^{(m)},y}(X,\nu^1,\nu^2,\nu^3)}\quad\times\prod_{i+1}^3\N_0(d\nu^i)\Bigr)\Bigr|\,dz^{(m)}\cr
&:=\liminf_{\delta\to 0}E_1^\delta+3E_2^\delta+3E_3^\delta+E_4^\delta.&(7.5)}$$

Consider $E_1^\delta$ first.  An explicit differentiation shows
$$|D^k_up_u(w)|\le cp_{2u}(w)u^{-k}\hbox{ for }k=1,2,3,\eqno(7.6)$$
which implies
$$|D^k_I G^\delta_{t,z^{(m)},y}(I,X)|\le c_2q_\delta(y_{m+1},X)I^{-k}\prod_{i+1}^m p_{2\delta+4\gamma_i^0I}(z_i-y_i+b^0_it),\hbox{ for }k=1,2,3.\eqno(7.7)$$
Use this (with $k=3$) with the Fundamental Theorem of Calculus to see that on $\{\nu^i_t>0\hbox{ for }i=1,2,3\}$, 
$$\eqalignno{\int&|y_j-z_j|^q|\Delta^3_t G^\delta_{t,z^{(m)},y}(X,\nu^1,\nu^2,\nu^3)|dz^{(m)}\cr
&\le c_2 \int\int|y_j-z_j|^q1_{(u_k\le \int_0^t\nu^k_sds,k=1,2,3)}I_t^{-3}q_\delta(y_{m+1},X_t^{(m+1)})\cr
&\phantom{\le c \int\int|y_j-z_j|^q1(}\times\prod_{i=1}^m[p_{2\delta+4\gamma_i^0(I_t+\sum_{k=1}^3u_k)}(z_i-y_i+b_i^0t)dz_i]\prod_{k=1}^3du_k\cr
&\le cI_t^{-3}q_\delta(y_{m+1},X_t^{(m+1)})\Bigl(t^q+\Bigl(\delta+I_t+\sum_{k=1}^3\int_0^t\nu_s^kds\Bigr)^{q/2}\Bigr)\prod_{k=1}^3\int_0^t\nu^k_sds.&(7.8)}$$
(3.14) and (3.16) imply 
$$\int \nu_s^p\,\N_0(d\nu)=(\gamma s)^{p-1}\Gamma(p+1)\hbox{ for }p\ge 0,\eqno(7.9)$$
and so by Jensen
$$\eqalignno{\int\Bigl(\int_0^t\nu_sds\Bigr)^{q/2+1}\N_0(d\nu)&\le t^{q/2+1}\int\int\nu_s^{q/2+1}\N_0(d\nu){ds\over t}\cr
&=\Gamma(2+q/2)t^{q/2}\int_0^t(\gamma s)^{q/2}ds\le ct^{q+1}.&(7.10)}$$
This bound, (7.9) with $p=1$, and (7.8), together with the expression for $E^\delta_1$, shows that 
$$\eqalignno{E_1^\delta&\le c\E_{z_{m+1}}(I_t^{-3}q_\delta(y_{m+1},X_t^{(m+1)})(t^q+\delta^{q/2}+I_t^{q/2})t^3)\cr
&\le c \E_{z_{m+1}}\Bigl(q_\delta(y_{m+1},X_t^{(m+1)})[(t^{q+3}+\delta^{q/2}t^3)\E_{z_{m+1}}(I_t^{-3}|X_t^{(m+1)})\cr
&\phantom{\le c \E_{z_{m+1}}\Bigl(q_\delta(y_{m+1},X_t^{(m+1)})[}+t^3\E_{z_{m+1}}(I_t^{q/2-3}|X_t^{(m+1)})]\Bigr)\cr
&\le c\E_{z_{m+1}}\Bigl(q_\delta(y_{m+1},X_t^{(m+1)}) \Bigr)[(t^q+\delta^{q/2})(z_{m+1}+t)^{-3}+t^{q/2}(t+z_{m+1})^{q/2-3}]\cr
&\le ct^{q-3}\E_{z_{m+1}}\Bigl(q_\delta(y_{m+1},X_t^{(m+1)})\Bigr)(1+(\delta^{1/2}/t)^{q})(1+z')^{q/2-3},&(7.11)}$$
where Lemma 3.2 is used in the next to last inequality.  

Let us jump ahead to $E_4^\delta$ which will be the dominant (and most interesting) term.  We use the decomposition and notation from Lemma~3.4 with $\rho=0$.  Let $S_n=\sum_{i=1}^ne_i(t)$, $R_n=\sum_{i=1}^nr_i(t)$, $p_k(z')=e^{-z'}{(z')^k\over k!}=P(N_0(t)=k)$, and $N^{(k)}=N(N-1)\dots(N-k+1)$.  What follows is an integration by parts formula on function space. 
Recalling that $\N_0(\nu_t>0)=(\gamma t)^{-1}$ from (3.14), we have from Lemma~3.4 and the exponential law of $\mu_t$ under $P^*_t$ (recall (3.16)) that
$$\eqalignno{\Bigl|&\E_{z_{m+1}}\Bigl(\int \Delta^3_t G^\delta_{t,z^{m)},y}(X,\nu^1,\nu^2,\nu^3)\prod_{i=1}^3(1_{(\nu_t^i>0)}\N_0(d\nu^i))\Bigr)\Bigr|\cr
&=(\gamma t)^{-3}\Bigl|\E_{z_{m+1}}\Bigl(G^\delta_{t,z^{(m)},y}(R_{N_0+3}+I_2(t),S_{N_0+3}+X_0'(t))\cr
&\phantom{=(\gamma t)^{-3}\Bigl|\E_{z_{m+1}}\Bigl(}-3G^\delta_{t,z^{(m)},y}(R_{N_0+2}+I_2(t),S_{N_0+2}+X_0'(t))\cr
&\phantom{=(\gamma t)^{-3}\Bigl|\E_{z_{m+1}}\Bigl(}+3G^\delta_{t,z^{(m)},y}(R_{N_0+1}+I_2(t),S_{N_0+1}+X_0'(t))\cr
&\phantom{=(\gamma t)^{-3}\Bigl|\E_{z_{m+1}}\Bigl(}-G^\delta_{t,z^{(m)},y}(R_{N_0}+I_2(t),S_{N_0}+X_0'(t))\Bigl)\Bigl|\cr
&=(\gamma t)^{-3}\Bigl|\sum_{n=0}^\infty(p_{n-3}(z')-3p_{n-2}(z')+3p_{n-1}(z')-p_n(z'))\cr
&\phantom{=(\gamma t)^{-3}\Bigl|\sum_{n=0}^\infty(}\times \E_{z_{m+1}}(G^\delta_{t,z^{(m)},y}(R_{n}+I_2(t),S_{n}+X_0'(t)))\Bigr|\cr
&=(\gamma t)^{-3}\Bigl|\sum_{n=0}^\infty p_n(z')(z')^{-3}[n^{(3)}-3n^{(2)}z'+3n(z')^2-(z')^3]\cr
&\phantom{=(\gamma t)^{-3}\Bigl|\sum_{n=0}^\infty(}\times \E_{z_{m+1}}(G^\delta_{t,z^{(m)},y}(R_{n}+I_2(t),S_{n}+X_0'(t)))\Bigr|\cr
&=z_{m+1}^{-3}\Bigl|\E_{z_{m+1}}\Bigl([N_0^{(3)}-3N_0^{(2)}z'+3N_0(z')^2-(z')^3]G^\delta_{t,z^{(m)},y}(I_t,X_t^{(m+1)})\Bigr)\Bigr|.&(7.12)}$$
In the last line we have again used Lemma~3.4 to reconstruct $(I_t,X^{(m+1)}_t)$.  

We also have
$$\eqalignno{\int&|y_j-z_j|^qG^\delta_{t,z^{(m)},y}(I_t,X_t^{(m+1)})dz^{(m)}\cr
&=q_\delta(y_{m+1},X_t^{(m+1)})\int |y_j-z_j|^qp_{\delta+2\gamma_j^0I_t}(z_j-y_j+b_j^0t)dz_j\cr
&\le c q_\delta(y_{m+1},X_t^{(m+1)})[t^q+\delta^{q/2}+I_t^{q/2}].&(7.13)}$$

Combine (7.12) and (7.13) to derive
$$\eqalignno{E_4^{\delta}&\le cz^{-3}_{m+1}\E_{z_{m+1}}\Bigl(|N_0^{(3)}-3N_0^{(2)}z'+3N_0(z')^2-(z')^3|\cr
&\phantom{|E_4^{\delta}|\le cz^{-3}_{m+1}\E_{z_{m+1}}\Bigl(}\times q_\delta(y_{m+1},X_t^{(m+1)})[t^{q}+\delta^{q/2}+I_t^{q/2}]\Bigr).&(7.14)}$$
Apply Jensen's inequality (as $q/2\le 1$) in Corollary~3.15 to see that
$$\E_{z_{m+1}}(I_t^{q/2}|N_0,X_t^{(m+1)})\le c[t^q+t^qN_0^{q/2}+t^{q/2}(X_t^{(m+1)})^{q/2}+t^{q/2}z_{m+1}^{q/2}],\eqno(7.15)$$
and so
$$\eqalignno{E_4^\delta&\le c t^{q-3}z'^{-3}\E_{z_{m+1}}\Bigl(|N_0^{(3)}-3N_0^{(2)}z'+3N_0(z')^2-(z')^3|q_\delta(y_{m+1}X_t^{(m+1)})\cr
&\phantom{\le c t^{q-3}z'^{-3}\E_{z_{m+1}}\Bigl(}\times[1+(\delta^{1/2}/t)^{q}+(X_t^{(m+1)}/t)^{q/2}+(z')^{q/2}+N_0^{q/2}]\Bigr)&(7.16)}$$

Next consider $E_2^\delta$.  Not surprisingly the argument is a combination of the ideas used to bound $E_1^\delta$ and $E^\delta_4$.  Define
$$\eqalign{H^\delta_{t,z^{(m)},y}(I,X,\nu^2,\nu^3)=&G^\delta_{t,z^{(m)},y}\Big(I+\int_0^t\nu^2_s+\nu^3_s\,ds,X\Big)-G^\delta_{t,z^{(m)},y}\Big(I+\int_0^t\nu^2_s\,ds,X\Big)\cr
&-G^\delta_{t,z^{(m)},y}\Big(I+\int_0^t\nu^3_s\,ds,X\Big)+G^\delta_{t,z^{(m)},y}(I,X).}$$

Now apply the decomposition in Lemma~3.4 with $\rho=1/2$ so that $N_{1/2}(t)$ is Poisson with mean $z'/2$.  Arguing as in the derivation of (7.12), but now with a simpler first order summation by parts (which we leave for the reader), we obtain
$$\eqalignno{\Bigl|&\E_{z_{m+1}}\Bigl(\int \Delta^3_t G^\delta_{t,z^{(m)},y}(X,\nu^1,\nu^2,\nu^3)1_{(\nu^1_t>0,\nu^2_t=\nu^3_t=0)}\prod_{k=1}^3\N_0(d\nu^k)\Bigr)\Bigr|\cr
&=(\gamma t)^{-1}\Bigl|\E_{z_{m+1}}\Bigl(\int\int H^\delta_{t,z^{(m)},y}\Big(\int_0^t\Xm_s+\nu^1_sds,\Xm_t+\nu^1_t,\nu^2,\nu^3\Big)\cr
&\phantom{=(\gamma t)^{-1}}-H^\delta_{t,z^{(m)},y}\Big(\int_0^t\Xm_sds,\Xm_t,\nu^2,\nu^3\Big)P^*_t(d\nu^1)1_{(\nu^2_t=\nu^3_t=0)}\prod_{k=2}^3\N_0(d\nu^k)\Bigr)\Bigr|\cr
&=(\gamma t)^{-1}\Bigl|\E_{z_{m+1}}\Bigl(\int\int H^\delta_{t,z^{(m)},y}(I_2(t)+R_{N_{1/2}+1},X_0'(t)+S_{N_{1/2}+1},\nu^2,\nu^3)\cr
&\phantom{=(\gamma t)^{-1}}-H^\delta_{t,z^{(m)},y}(I_2(t)+R_{N_{1/2}},X_0'(t)+S_{N_{1/2}},\nu^2,\nu^3)\prod_{k=2}^3\N_0(d\nu^k)\Bigr)\Bigr|\cr
&={2\over z_{m+1}}\Bigl|\int \E_{z_{m+1}}\Bigl((N_{1/2}-z'/2)H^\delta_{t,z^{(m)},y}(I(t),X^{(m+1)}_t,\nu^2,\nu^3)\Bigr)\prod_{k=2}^3\N_0(d\nu^k)\Bigr|.&(7.17)}$$

Now use (7.7) (with $k=2$) and argue as in (7.8) to see that
$$\eqalignno{\int&\int |y_j-z_j|^q|H^\delta_{t,z^{(m)},y}(I,X,\nu^2,\nu^3)|dz^{(m)}\prod_{k=2}^3\N_0(d\nu^k)\cr
&\le cI^{-2}q_\delta(y_{m+1}X)\int(t^q+(\delta+I+\sum_2^3\int_0^t \nu_s^kds)^{q/2})\prod_{k=2}^3\Bigl[\int_0^t\nu_s^kds\N_0(d\nu^k)\Bigr]\cr
&\le c I^{-2}q_\delta(y_{m+1},X)[t^q+\delta^{q/2}+I^{q/2}]t^2,&(7.18)}$$
where the last line uses (7.10).  

Take the absolute values inside the inside the integral in (7.17), multiply by $|y_j-z_j|^q$, integrate with respect to $z^{(m)}$, and use the above bound to conclude that
$$\eqalignno{E_2^\delta&\le cz_{m+1}^{-1}t^2\E_{z_{m+1}}\Bigl(|N_{1/2}-z'/2|q_\delta(y_{m+1},\Xm_t)\cr&\phantom{le cz_{m+1}^{-1}t^2\E_{z_{m+1}}\Bigl(|}\times\Bigl[(t^q+\delta^{q/2})I(t)^{-2}+I(t)^{q/2-2}\Bigr]\Bigr).&(7.19)}$$
If $r\ge 0$, 
the independence of $X'_0$ from $(N_{1/2},\{e_j\})$ and Lemma~3.2, applied to $X'_0$, imply that
$$\eqalignno{\E_{z_{m+1}}&(I(t)^{-r}|\Xm_t,N_{1/2})\cr
&\le \E\Bigl(\E\Bigl(\Bigl(\int_0^t X'_0(s)ds\Bigr)^{-r}|X'_0(t),N_{1/2},\{e_j\}\Bigr)\Bigr|\Xm_t,N_{1/2})\cr
&=\E\Bigl(\E\Bigl(\Bigl(\int_0^t X_0'(s)ds\Bigr)^{-r}\Bigl|X'_0(t)\Bigr)\Bigl|\Xm_t,N_{1/2}\Bigr)\cr
&\le c(t+z_{m+1}/2)^{-r}t^{-r}.&(7.20)}$$
The last line is where it is convenient that $\rho=1/2>0$.

Use (7.20) in (7.19) with $r=2$ and $2-q/2$.  After a bit of algebra this leads to
$$\eqalignno{E_2^\delta&\le c t^{q-3}(z')^{-1}\E_{z_{m+1}}(|N_{1/2}-z'/2|q_\delta(y_{m+1},\Xm_t))\cr
&\phantom{le c t^{q-3}(z')^{-1}\E_{z_{m+1}}(}\times\Bigr[\Bigl({\sqrt\delta\over t}\Bigr)^q(1+z')^{-2}+(1+z')^{q/2-2}\Bigr].&(7.21)}$$

The argument for $E_3^\delta$ is similar to the above. One works with
$$\tilde H_{t,z^{(m}),y}^\delta(I,X,\nu^3)=G^\delta_{t,z^{(m)},y}\Big(I+\int_0^t \nu^3_sds,X\Big)-G^\delta_{t,z^{(m)},y}(I,X).$$
The required third order difference of $G^\delta_{t,z^{(m)},y}$ on $\{\nu^1_t>0,\nu^2_t>0,\nu^3_t=0\}$ is now a second order difference of $\tilde H^\delta_{t,z^{(m)},y}$.  Minor modifications of the derivation of (7.21) lead to
$$\eqalignno{E_3^\delta&\le ct^{q-3}(z')^{-2}\E_{z_{m+1}}(|N_{1/2}^{(2)}-N_{1/2}z'+(z'/2)^2|q_\delta(y_{m+1},\Xm_t))\cr
&\phantom{\le ct^{q-3}(z')^{-2}\E_{z_{m+1}}(|}\times [(\sqrt\delta/t)^q(1+z')^{-1}+(1+z')^{q/2-1}].&(7.22)}$$

The above bounds in $E_i^\delta$ $i=1,\dots 4$ may be used in (7.5) and after the terms involving $\sqrt \delta/t$ are neglected (for $q=0$ these terms are bounded by their neighbours, and for $q>0$, if they do not approach $0$, the right side below must be infinite) we find
$$\eqalignno{\int& |y_j-z_j|^q|D^2_{z_{m+1}}p_t(z,y)|dz^{(m)}\cr
&\le ct^{q-3}\liminf_{\delta\to0}\E_{z_{m+1}}\Bigl(q_\delta(y_{m+1},\Xm_t)\Bigl[(1+z')^{q/2-1}&(7.23)\cr
&\phantom{\le ct^{q-3}}\times\Bigl((1+z')^{-2}+|N_{1/2}-z'/2|(1+z')^{-1}(z')^{-1}+|N_{1/2}^{(2)}-N_{1/2}z'+(z'/2)^2|(z')^{-2}\Bigr)\cr
&\phantom{\le ct^{q-3}}+|N_0^{(3)}-3N_0^{(2)}z'+3N_0(z')^{2}-(z')^3|(z')^{-3}[1+(\Xm_t/t)^{q/2}+(z')^{q/2}+N_0^{q/2}]\Bigr]\Bigr).}$$

The required bound follows from the above by a bit of algebra but as the reader may be fatigued at this point we point out the way. Trivial considerations show it suffices to show the following inequalities for $n_0,n_{1/2}\in\Z_+$ and $z'\ge0$:
$$\eqalignno{|&n_0^{(3)}-3n^{(2)}_0z'-3n_0(z')^2-(z')^3|(z')^{-3}\cr
&\le c[1_{(n_0\le 2)}(1+(z')^{-n_0})+1_{(n_0\ge 3)}(z')^{-3}(|n_0-z'|^3+3n_0|n_0-z'|+n_0),&(7.24)}$$
and
$$\eqalignno{[&(1+z')^{-2}+|n_{1/2}-z'/2|(1+z')^{-1}(z')^{-1}+|n_{1/2}^{(2)}-n_{1/2}z'+(z'/2)^2|(z')^{-2}\cr
&\le c[1_{(n_{1/2}\le 1)}+1_{(n_{1/2}=1)}{z'}^{-1}+1_{(n_{1/2}\ge 2)}{z'}^{-2}(n_{1/2}+(n_{1/2}-{z'\over 2})^2)].&(7.25)}$$
(7.24) is easy.  (7.25) reduces fairly directly to showing that for $n_{1/2}\ge 2$,
$$(1+z')^{-1}\le c(n_{1/2}+(n_{1/2}-z'/2)^2)(z')^{-2}.$$
If $z'\le 1$ this is trivial and for $z'>1$ consider $n_{1/2}\le z'/4$ and $n_{1/2}> z'/4$ separately.  This completes the proof of (a).

(b) The proof of this second order version of (a) is very similar to, but simpler than that of (a).  One now only has a second order difference and three $E_i^\delta$ terms to consider.  In fact we will not actually need $q>0$ in (a) but included it so that the reader will not complain about the missing details in the proof of (b) (where $q>0$ has been used in Proposition 4.12).  We do comment on the lack of $N_{1/2}$ in this bound.  

An argument similar to that leading to (7.23) shows that \hfil\break
$\int|z_j-y_j|^q|D^2_{z_{m+1}}p_t(z,y)|dz^{(m)}$ is bounded by
$$\eqalignno{c&t^{q-2}\liminf_{\delta\to0}\E_{z_{m+1}}\Bigl(q_\delta(y_{m+1},\Xm_t)\Bigl[(1+z')^{-1+q/2}\Bigl((1+z')^{-1}+|N_{1/2}-z'/2|(z')^{-1}\Bigr)\cr
&\phantom{}+|N_0^{(2)}-2N_0z'+(z')^2|(z')^{-2}(1+(\Xm_t/t)^{q/2}+(z')^{q/2}+N_0^{q/2})\Bigr]\Bigr)\cr
&\equiv ct^{q-2}\liminf_{\delta\to 0}\E_{z_{m+1}}\Bigl(q_\delta(y_{m+1},\Xm_t)[T_{1/2}+T_0]\Bigr).}$$
It is easy to check that 
$$T_0\le c\phi(z',\Xm_t/t,N_0),$$
and, using $N_{1/2}\le N_0$ from (3.23), that
$$T_{1/2}\le 2(1+z')^{-1+q/2}(1+N_0(z')^{-1}+1)\equiv\bar T_{1/2}.$$
Hence to prove (b), it remains to verify
$$\bar T_{1/2}\le c\phi(z',\Xm_t/t,N_0).$$
Trivial considerations reduce this to showing that $(1+z')^{-1+q/2}\le c\phi(z',\Xm_t/t,N_0)$ for $N_0\ge 2$.  This is easily verified by considering $N_0<z'/2$ and $N_0\ge z'/2$ separately. 

(c) Note that (7.3) is the first order version of (a) and (b) with $q=0$. The proof is substantially simpler, but, as it plays the pivotal role in the proof for the important $2$-dimensional case, we give the proof in Section~8.  (7.4) then follows immediately since the spatial homogeneity in the first $m$ variables, (3.7), implies 
$$p_t(z,y)=p_t(-y^{(m)},z_{m+1},-z^{(m)},y_{m+1}).\eqno(7.26)$$\qed

\ms

\ni{\bf Proof of Lemma 4.5(b).} Let $J$ denote the integral to be bounded in the statement of (b), and $p_n(w)=e^{-w}w^n/n!$ be the Poisson probabilities.  Let $\Gamma_n$ be a Gamma random variable with density 
$$g_n(x)=x^{n+b/\gamma-1}e^{-x}\Gamma(n+b/\gamma)^{-1},\eqno(7.27)$$ 
and recall $z'=z_{m+1}/\gamma t$. By integrating the bound from Lemma~7.1(b) in $z_{m+1}$ (using Fatou's Lemma) we see that
$$\eqalignno{J&\le c_{7.1}t^{q-2}\liminf_{\delta\to 0}\int z_{m+1}^p\E_{z_{m+1}}\Bigl(q_\delta(y_{m+1},\Xm_t)\phi(z',\Xm_t/t,N_0)\Bigr)\,z_{m+1}^{b/\gamma-1}dz_{m+1}.}$$
Our formula for the joint distribution of $(\Xm_t,N_0)$ (Lemma 3.6(a)) allows us to evaluate the above and after changing variables and the order of integration we see that if $y'=y/\gamma t$, then
$$\eqalignno{J&\le c_{7.1}t^{q-2+p}\liminf_{\delta\to 0}\int(q_\delta(y_{m+1},y)\sum_{n=0}^\infty p_n(y')\cr
&\phantom{\le c_{7.1}}\times\Bigl[\int_0^\infty g_n(z')[(1_{(n\le 1)}+1_{(n=1)}(z')^{-1})(1+(y')^{q/2}+(z')^{q/2})\cr
&\phantom{\phantom{\le c_{7.1}}\times\Bigl[}+1_{(n\ge 2)}(z')^{-2}[(n-z')^2+n][(y')^{q/2}+(z')^{q/2}+n^{q/2}](z')^pdz'\Bigr]y^{b/\gamma-1}dy\cr
&=c_{7.1}t^{q-2+p}\liminf_{\delta\to 0}\int q_\delta(y_{m+1},y)\sum_{n=0}^\infty p_n(y')\cr
&\phantom{=c_{7.1}}\times\Bigl[1_{(n\le 1)}\E((1+(y')^{q/2}+\Gamma_n^{q/2})\Gamma_n^p)+1_{(n=1)}\E((1+(y')^{q/2}+\Gamma_n^{q/2})\Gamma_n^{-1+p})\cr
&\phantom{\phantom{=c_{7.1}}\times\Bigl[}+1_{(n\ge2)}\E(((n-\Gamma_n)^2+n)\Gamma_n^{-2+p}((y')^{q/2}+\Gamma_n^{q/2}+n^{q/2}))\Bigr]y^{b/\gamma-1}dy.&(7.28)}$$
There is a constant $c_0$ (as in Convention 3.1) so that
$$\E(\Gamma_n^r)\le c_0(n\vee1)^r\hbox{ for all }|r|\le 4\hbox{ and $n\in\Z_+$ satisfying }r+n\ge {-3\over 4 M_0^2}.\eqno(7.29)$$
Indeed the above expectation is $\Gamma(n+b/\gamma+r)/\Gamma(n+b/\gamma)$, where 
$$r+n+b/\gamma\ge  {-3\over 4 M_0^2}+M_0^{-2}= (2M_0)^{-2}.$$
The result now follows by an elementary, and easily proved, property of the Gamma function.

Assume now the slightly stronger condition
$$|r|\le 3, n\in\Z_+\hbox{ and } r+n\ge {-1\over 2M_0^2}.\eqno(7.30)$$
Then $\Gamma_n=\Gamma_0+S_n$, where $S_n$ is a sum of $n$ i.i.d.\  mean one exponential random variables.  If $s$ and $s'$ are H\"older dual exponents, where $s$ is taken close enough to $1$ so that  the conditions of (7.29) remain valid with $rs$ in place of $r$, then
$$\eqalignno{\E((&\Gamma_n-n)^2\Gamma_n^r)\le \E((\Gamma_n-n)^{2s'})^{1/s'}\E(\Gamma_n^{sr})^{1/s}\cr
&\le cn(n\vee 1)^r,&(7.31)}$$
where we have used an elementary martingale estimate for $|S_n-n|$ and (7.29).  Here $c$ again is as in Convention 3.1.  

We now use (7.31) and (7.29) to bound the Gamma expectations in (7.28).  It is easy to check that our bounds on $p$ and $q$ imply the powers we will be bounding satisfy (7.30).  This leads to
$$\eqalignno{J&\le ct^{q-2+p}\liminf_{\delta\to0}\int q_\delta(y_{m+1},y)\sum_{n+0}^\infty p_n(y')\cr
&\phantom{\le ct^{q-2+p}}\times\Bigl[1_{(n\le 1)}(1+(y')^{q/2})+1_{(n\ge 2)}n^{-1+p}((y')^{q/2}+n^{q/2})\Bigr]y^{b/\gamma-1}\,dy\cr
&\le ct^{q-2+p}\liminf_{\delta\to0}\int q_\delta(y_{m+1},y)\Bigl[e^{-y'}(1+y')(1+(y')^{q/2})\cr
&\phantom{\le ct^{q-2+p}}+\E(1_{(N(y')\ge 2)}((y')^{q/2}N(y')^{p-1}+N(y')^{q/2-1+p}))\Bigr]y^{b/\gamma-1}dy.&(7.32)}$$
In the last line $N(y')$ is a Poisson random variable with mean $y'$. Well-known properties of the Poisson distribution show
that for a universal constant $c_2$
$$\E(N(y')^r1_{(N(y')\ge 2)})\le h_r(y')\equiv c_2(1+y')^r\hbox{ for all } y'\ge 0, |r|\le 2.\eqno(7.33)$$
For negative values of $r$ see Lemma~4.3(a) of [BP] where the constant depends on $r$ but the argument there easily shows for $r$ bounded one gets a uniform constant.  If 
$$h(y')=e^{-y'}(1+y')(1+(y')^{q/2})+h_{q/2-1+p}(y')+y'^{q/2}h_{p-1}(y'),$$
then clearly
$$h(y')\le c_3(1+y')^{q/2-1+p}.$$
As all of the powers appearing in (7.32) satisfy the bounds in (7.33), we may use (7.33) to bound
the left-hand side of (7.32) and arrive at
$$\eqalign{J&\le ct^{q-2+p}\liminf_{\delta\to0}\E_{y_{m+1}}(h(\Xm_\delta/\gamma t))\cr
&\le ct^{q-2+p}\liminf_{\delta\to0}\E_{y_{m+1}}((1+(\Xm_\delta/\gamma t)^{q/2-1+p})\cr
&=ct^{q-2+p}(1+y')^{q/2-1+p}\qq \hbox{ (Dominated Convergence)}\cr
&=ct^{q/2-1}(t+y_{m+1})^{q/2-1+p}.}$$\qed

\ms

\ni{\bf Proof of Lemma 4.5(c).} The spatial homogeneity (7.26) shows the integral being bounded equals
$$\int|y_j-(-z_j)|^q z_{m+1}|D^2_{z_{m+1}}p_t(y^{(m)},z_{m+1},-z^{(m)},y_{m+1})|dy^{(m)}\mu_{m+1}(dy_{m+1}).$$
This shows we can again use the upper bound in Lemma~7.1(b) to bound the integral over $y^{(m)}$ in the above. One then must integrate the resulting bound in $y_{m+1}$ instead of $z_{m+1}$.  This actually greatly simplifies the calculation just given as one can integrate $y_{m+1}$ at the beginning  and hence the $q_\delta$ term conveniently disappears (see the proof of (4.14) below). For example, if we neglect the insignificant $n\le 1$ contribution to $\phi$ in Lemma~7.1, the resulting integral is bounded by
$$ct^{q-1}(z')^{-1}\E(1_{(N_0\ge 2)}(N_0+(N_0-z')^2)((z')^{q/2}+N_0^{q/2}+(\Xm_t/t)^{q/2})).$$
This can be bounded using elementary estimates of the Poisson and H\"older's inequality, the latter being
much simpler than invoking Lemma~3.6.  We omit the details.  \qed

\ms 

\ni{\bf Proof of Lemma 4.5(a).}  (4.13) is the first order version of Lemma 4.5 (b) and we omit the proof which is much simpler. (4.14) is a bit different from (c).  Integrate (7.4) over $y_{m+1}$ to see that
$$\eqalignno{\int &y^p_{m+1}|D_{z_{m+1}}p_t(z,y)|\mu(dy)\cr
&\le c_{7.1}t^{-1}\liminf_{\delta\to0}\int y_{m+1}^p \E_{z_{m+1}}\Bigl(q_\delta(y_{m+1},\Xm_t)[(1+z')^{-1}+|N_0-z'|/z']\Bigr)\mu_{m+1}(dy_{m+1})\cr
&=c_{7.1}t^{-1} \liminf_{\delta\to0}\E_{z_{m+1}}\Bigl([(1+z')^{-1}+|N_0-z'|/z']\E_{\Xm_t}((\Xm_\delta)^p)\Bigr).}$$
Now use the moment bounds in Lemma~3.3(d,e) to bound the above by
$$ct^{-1}\E_{z_{m+1}}\Bigl([(1+z')^{-1}+|N_0-z'|/z'](\Xm_t)^p\Bigr).\eqno(7.34)$$
The first term is trivially bounded by the required expression using Lemma~3.3 again. Using the 
joint density formula (Lemma~3.6), the Gamma power bounds (7.29), and arguing as in the proof of (b) above, the term in (7.34) involving $N_0$ is at most
$$ct^{p-1}(z')^{-1}\E(|N-z'|(N\vee 1)^p),\eqno(7.35)$$
where $N=N(z')$ is a Poisson random variable with mean $z'$.  We have 
$$\E(|N-z'|N^p1_{(N>0)})\le c_0(z'\wedge(z')^{1/2+p})\hbox{ for all }z'>0 \hbox{ and }-1\le p\le 1/2.\eqno(7.36)$$
For $p\le0$ Lemma~3.3 of [BP] shows this (the uniformity for bounded $p$ is again clear).  For $1/2\ge p>0$ use Cauchy-Schwarz to prove (7.36).  Separating out the contribution from $N=0$, we see from (7.36) that (7.35) is at most
$$ct^{p-1}(z')^{-1}[e^{-z'}z'+(z'\wedge (z')^{1/2+p})]\le c t^{p-1}(e^{-z'}+1\wedge(z')^{p-1/2})\le ct^{p-1}(z'+1)^{p-1/2}.$$
The result follows.  
 \qed

\ms

\ni{\bf Proof of Lemma 4.5(f).}  By the spatial homogeneity in the first $m$ variables (7.26) we may use Lemma~4.4 to conclude
$$|D_{z_j}D^2_{z_{m+1}}p_t(z,y)|=\Bigl|\int D_{z_j}p_{t/2}(x,-z^{(m)},y_{m+1})D^2_{z_{m+1}}p_{t/2}(-y^{(m)},z_{m+1},x)\mu(dx)\Bigr|.$$
Therefore 
$$\eqalign{\int &z_{m+1}^p|D_{z_j}D^2_{z_{m+1}}p_t(z,y)|\mu(dz)\cr
&\le \int\int\Bigl[\int|D_{z_j}p_{t/2}(x,-z^{(m)},y_{m+1})|dz^{(m)}\Bigr]\cr
&\phantom{\le \int\int\Bigl[\int|D_{z_j}}\times z_{m+1}^p|D^2_{z_{m+1}}p_{t/2}(-y^{(m)},z_{m+1},x)|\mu(dx)\mu_{m+1}(dz_{m+1}).}$$
Use Lemma~4.2(a) to bound the first integral in square brackets and so bound the above by
$$\eqalign{c&t^{-1/2}\int z_{m+1}^p|D^2_{z_{m+1}}p_{t/2}(-y^{(m)},z_{m+1},x)|dx^{(m)}\,\mu_{m+1}(dz_{m+1})\cr
&\phantom{t^{-1/2}\int |D^2_{z_{m+1}}}\times(t+x_{m+1}+y_{m+1})^{-1/2}q_{t/2}(x_{m+1},y_{m+1})\mu_{m+1}(dx_{m+1}).}$$
The spatial homogeneity (7.26) implies
$$p_{t/2}(-y^{(m)},z_{m+1},x)=p_{t/2}(-x^{(m)},z_{m+1},y^{(m)},x_{m+1}),$$
and so we conclude from the above that
$$\eqalign{\int & z_{m+1}^p|D_{z_j}D^2_{z_{m+1}}p_t(z,y)|\mu(dz)\cr
&\le ct^{-1/2}\int\Bigl[\int z^p_{m+1} |D^2_{z_{m+1}}p_{t/2}(-x^{(m)},z_{m+1},y^{(m)},x_{m+1})|dx^{(m)}\mu_{m+1}(dz_{m+1})\Bigr]\cr
&\phantom{\le ct^{-1/2}\int\Bigl[\int z^p_{m+1} |D^2_{z_{m+1}}}\times(t+x_{m+1}+y_{m+1})^{-1/2}q_{t/2}(x_{m+1},y_{m+1})\mu_{m+1}(dx_{m+1})\cr
&\le ct^{-3/2}\int (t+x_{m+1})^{p-1}(t+x_{m+1}+y_{m+1})^{-1/2}q_{t/2}(x_{m+1},y_{m+1})\mu_{m+1}(dx_{m+1})\cr
&\le c t^{-3/2}\E_{y_{m+1}}\Bigl((t+X^{(m+1)}_t)^{p-3/2}\Bigr)\le ct^{p-3}.}$$
We have used Lemma 4.5(b) with $q=0$ in the next to last inequality and $p\le 3/2$ in the last line.  \qed

\ms

\ni{\bf Proof of Lemma 4.5(e).}   For (4.18), use Lemma~4.4 and the spatial homogeneity (7.26) to bound the left-hand side of (4.18) by
$$\eqalign{\int\int\Bigl[\int y_{m+1}^p&|D_{y_{m+1}}\hat p_{t/2}(-x^{(m)},y_{m+1},-y^{(m)},x_{m+1})|\cr
&\times z_{m+1}|D_{z_{m+1}}^2 p_{t/2}(0,z_{m+1},x^{(m)}-z^{(m)},x_{m+1})|\mu(dx)\Bigl]dz^{(m)}\mu_{m+1}(dy_{m+1}).}$$
Use the substitution (for $z^{(m)}$) $w=x^{(m)}-z^{(m)}$ and do the $dx^{(m)}\mu_{m+1}(dy_{m+1})$ integral first, using (4.13) to bound this integral by $c_{4.5}t^{p-1}$ (as $p\le 1/2$).  Now use (4.5) to bound the remaining $dw \mu_{m+1}(dx_{m+1})$ integral by $c_{4.5}t^{-1}$.

The derivation of (4.19) is almost the same as above.  One uses Lemma~4.2(b) now to bound the first integral.  

\ms

\ni{\bf Proof of Lemma 4.5(d).} The approach is similar to that in (b) as we integrate the bound in Lemma~7.1(a).  There is some simplification now even with the higher derivative as $q=0$.  We use the notation from that proof, so that $g_n$ is the Gamma density in (7.27), $p_n(w)$ are the Poisson probabilities with mean $w$, and $\Gamma_n$ is a random variable with density $g_n$.  Also let $B_n$ be a Binomial $(n,1/2)$ random variable independent of $\Gamma_n$.  To ease the transition to Lemma~4.6 we replace $t$ with $s$ in this calculation.  We also keep the notation $z'=z_{m+1}/\gamma s$, $y'=y/\gamma s$.  If
$$\phi_1(z',k)=(z'+1)^{-1}[1_{(k\le 1)}(1+(z')^{-k})+1_{(k\ge 2)}(z')^{-2}[k+(k-z'/2)^2]],$$
and
$$\phi_2(z',k)=1_{(k\le2)}(1+(z')^{-k})+1_{(k\ge 3)}(z')^{-3}[|k-z'|^3+3k|k-z'|+k],$$
then by Lemma~7.1(a) and Fatou's lemma, the integral we need to bound is at most
$$\eqalignno{&\liminf_{\delta\to0}cs^{-3}\Bigl(\int \E_{z_{m+1}}\Bigl(q_\delta(y_{m+1},\Xm_s)\phi_1(z',N_{1/2})\Bigr)z_{m+1}^{3/2}\mu_{m+1}(dz_{m+1})\cr
&\phantom{\liminf_{\delta\to0}cs^{-3}\Bigl(}+\int \E_{z_{m+1}}\Bigl(q_\delta(y_{m+1},\Xm_s)\phi_2(z',N_{0})\Bigr)z_{m+1}^{3/2}\mu_{m+1}(dz_{m+1})\Bigr)\cr
&:=\liminf_{\delta\to0} cs^{-3}[J_1(\delta)+J_2(\delta)].&(7.37)}$$

By Lemma 3.6(b),
$$\eqalignno{J_1(\delta)&\le cs^{3/2}\int_0^\infty q_\delta(y_{m+1},y)\sum_{n=0}^\infty p_n(y')\int_0^\infty\sum_{k=0}^n\choose n k 2^{-n}g_n(z')\phi_1(z',k)(z')^{3/2}dz' y^{b/\gamma-1}dy\cr
&=cs^{3/2}\int_0^\infty q_\delta(y_{m+1},y)\sum_{n=0}^\infty p_n(y')\E(\phi_1(\Gamma_n,B_n)\Gamma_n^{3/2})y^{b/\gamma-1}dy.&(7.38)}$$

The moment bounds (7.29) (the conditions there will be trivially satisfied now) give us
$$\eqalignno{\E(\phi_1(\Gamma_n,B_n)\Gamma_n^{3/2})&\le c\Bigl(P(B_n=0)\E(\Gamma_n^{3/2})+P(B_n=1)\E(\Gamma_n^{1/2}+\Gamma_n^{3/2})\cr
&\quad +\E(1_{(B_n\ge 2)}B_n)\E(\Gamma_n^{-3/2})+\E(1_{(B_n\ge 2)}(B_n-\Gamma_n/2)^2\Gamma_n^{-3/2})\Bigr)\cr
&\le c\Bigl(2^{-n}(1+n)(n^{1/2}+n^{3/2})+1_{(n\ge 2)}n^{-1/2}\cr
&\phantom{\le c\Bigl(2^{-n}(1+n)}+1_{(n\ge 2)}\E(\E((B_n-\Gamma_n/ 2)^2|\Gamma_n)\Gamma_n^{-3/2})\Bigr).&(7.39)}$$
The conditional expectation in the last term is $[(\Gamma_n-n)^2+\Gamma_n]/4$.  Therefore we may now use (7.31) and also (7.29) (as $n\ge 2$ and $r=-1/2$ or $-3/2$ the conditions there are satisfied) to 
see that for $n\ge2$
$$\eqalignno{\E&(\E((B_n-\Gamma_n/ 2)^2|\Gamma_n)\Gamma_n^{-3/2})
\le( \E((\Gamma_n-n)^2\Gamma_n^{-3/2})+\E(\Gamma_n^{-1/2}))/4\le cn^{-1/2}.}$$
Insert this bound into (7.39) and conclude that
$$\eqalign{\E(\phi_1(\Gamma_n,B_n)\Gamma_n^{3/2})&\le c(2^{-n}(1+n^{5/2})+1_{(n\ge 2)}n^{-1/2})\le c.}$$
Therefore we can sum over $n$ and integrate over $y$ in (7.38) and obtain
$$J_1(\delta)\le c_1s^{3/2},\eqno(7.40)$$
where as always $c_1$ satisfies Convention~3.1.  

Lemma 3.6(a) and the argument leading to (7.38) shows that 
$$J_2(\delta)\le 
cs^{3/2}\int_0^\infty q_\delta(y_{m+1},y)\sum_{n=0}^\infty p_n(y') \E(\phi_2(\Gamma_n,n)\Gamma_n^{3/2})y^{b/\gamma-1}dy.\eqno(7.41)$$
We have
$$\eqalign{\E&(\phi_2(\Gamma_n,n)\Gamma_n^{3/2})\cr
&=1_{(n\le 2)}\E((1+\Gamma_n^{-n})\Gamma_n^{3/2})+1_{(n\ge 3)}\E(\Gamma_n^{-3/2}(|n-\Gamma_n|^3+3n|n-\Gamma_n|+n)).}$$
Some simple Gamma distribution calculations like those in the proof of (b), and which the reader can easily provide (recall Convention~3.1), show
that the above is bounded by a constant depending only on $M_0$.   As before by using this bound in (7.41) and integrating out $n$ and $y$ we arrive at
$$J_2(\delta)\le c_2s^{3/2}.\eqno(7.42)$$
Insert the above bounds on $J_i(\delta)$ into (7.37) to complete the proof.  \qed

\ms
\ni{\bf Proof of Lemma 4.6.} Consider (4.22).  The functions $\phi_1$ and $\phi_2$ are as in the previous argument.  Argue just as in the derivation of (7.37) to bound the left-hand side of (4.22) by
$$ \eqalignno{&\liminf_{\delta\to0}cs^{-3}\Bigl[t^{b/\gamma}\int_0^{\gamma t}\E_{z_{m+1}}\Bigl(q_\delta(y_{m+1},\Xm_s)\phi_1(z',N_{1/2})\Bigr)dz_{m+1}\cr
&\phantom{\liminf_{\delta\to0}cs^{-3}\Bigl[}+t^{b/\gamma}\int_0^{\gamma t}\E_{z_{m+1}}\Bigl(q_\delta(y_{m+1},\Xm_s)\phi_2(z',N_{0})\Bigr)dz_{m+1}\Bigr]\cr
&:= \liminf_{\delta\to0}cs^{-3}[K_1(\delta)+K_2(\delta)].&(7.43)}$$
Note we are integrating with respect to $z_{m+1}$ and not $\mu_{m+1}(dz_{m+1})$ as in the previous calculation. Lemma~3.6(b) implies that  
$$\eqalignno{K_1(\delta)&\le c s\int_0^\infty q_\delta(y_{m+1},y)\sum_{n=0}^\infty p_n(y')\cr
&\phantom{\le c s\int_0^\infty q_\delta}\times(t/s)^{b/\gamma}\int_0^{t/s}\E(\phi_1(z',B_n))e^{-z'}{(z')^n\over \Gamma(n+b/\gamma)}dz'y^{b/\gamma-1}dy\cr
&\le c s\int_0^\infty q_\delta(y_{m+1},y)\sum_{n=0}^\infty p_n(y')\int_0^{t/s}\E(\phi_1(z',B_n)){(z')^n\over (n+1)^2}dz'y^{b/\gamma-1}dy.&(7.44)}$$
We have bounded $\Gamma(n+b/\gamma)^{-1}$ in a rather crude manner in the last line.

For $0<z'\le t/s\le 1$ we have
$$\eqalignno{\E&(\phi_1(z',B_n))(z')^n(n+1)^{-2}\cr
&\le c\Bigl[\P(B_n=0)+\P(B_n=1)(1+(z')^{-1})\cr
&\phantom{\le c\Bigl[}+1_{(n\ge 2)}(z'+1)^{-1}(z')^{-2}\E(B_n+(B_n-z'/2)^2)\Bigr](z')^n(n+1)^{-2}\cr
&\le c\Bigl[2^{-n}((z')^n+n(1+(z')^{n-1})+1_{(n\ge 2)}(z')^{n-2}(n+(n-z')^2)(n+1)^{-2}\Bigr]\le c.}$$
This, together with (7.44), shows that
$$K_1(\delta)\le ct.$$

Next use Lemma~3.6(a) to see that 
$$\eqalignno{K_2(\delta)&\le cs \int_0^\infty q_\delta(y_{m+1},y)\sum_{n=0}^\infty p_n(y')\cr
&\phantom{\le c s \int_0^\infty q_\delta}\times(t/s)^{b/\gamma}\int_0^{t/s}\phi_2(z',n)e^{-z'}{(z')^n\over \Gamma(n+b/\gamma)}dz'y^{b/\gamma-1}dy\cr
&\le c s\int_0^\infty q_\delta(y_{m+1},y)\sum_{n=0}^\infty p_n(y')\int_0^{t/s}\phi_2(z',n){(z')^n\over (n+1)^3}dz'y^{b/\gamma-1}dy.}$$
As above, an elementary argument shows that for $0<z'\le 1$, $\phi_2(z',n)(z')^n(n+1)^{-3}$ is uniformly bounded in $n$, $z'$ and also $(b,\gamma)$ as in Convention~3.1.  Hence, we may infer
$$K_2(\delta)\le ct.$$
Put the bounds on $K_i(\delta)$ into (7.43) to complete the proof of (4.22).

We omit the proof of (4.21) which is the first order analogue of (4.22) and is considerably easier.
\qed

\ms

We need a probability estimate for Lemma~4.7. As usual $X^{(m+1)}$ is the Feller branching diffusion with generator (3.1).

\proclaim Lemma 7.2. (a) $\P_z(X^{(m+1)}_t\ge w)\le (w/z)^{b/2\gamma}\exp\Bigl\{{-(\sqrt z-\sqrt w)^2\over \gamma t}\Bigr\}\hbox{ for all }w> z\ge 0$.

\ni(b) $\P_z(X^{(m+1)}_t\le w)\le (w/z)^{b/2\gamma}\exp\Bigl\{{-(\sqrt z-\sqrt w)^2\over \gamma t}\Bigr\}\hbox{ for all }0\le w\le z$.

\proof This is a simple estimate using the Laplace transform in Lemma~3.3(c). Write $X_t$ for $X^{(m+1)}_t$.   If $-(\gamma t)^{-1}<\lam\le 0$, then
$$\P_z(X_t\ge w)\le e^{\lam w}\E_z(e^{-\lam X_t})=e^{\lam w}(1+\lam t\gamma)^{-b/\gamma}\exp\Bigl\{{-z\lam\over 1+\lam \gamma t}\Bigr\}.$$
If $\lam\ge 0$, then 
$$\P_z(X_t\le w)\le e^{\lam w}\E_z(e^{-\lam X_t})=e^{\lam w}(1+\lam t\gamma)^{-b/\gamma}\exp\Bigl\{{-z\lam\over 1+\lam \gamma t}\Bigr\}.$$
Now set $\lam={\sqrt{z/w}-1\over \gamma t}$ in both cases. This is in $(-(\gamma t)^{-1}, 0)$ if $0\le z<w$ and in $[0,\infty)$ if $0\le z\ge w$.  A bit of algebra gives the bounds.\qed

\ni{\bf Proof of Lemma 4.7.} We will again integrate the bound in Lemma~7.2(b) over $z_{m+1}$ but as $q=0$ we will use the function
$$\psi(z',n)=1_{(n\le 1)}+1_{(n=1)}(z')^{-1}+1_{(n\ge 2)}(z')^{-2}[(n-z')^2+n].$$
We then have from Lemma~7.2(b), that for each $z'>0$, 
$$\eqalignno{J&(z_{m+1})\cr
&\equiv \int\int(1_{(y_{m+1}\le w\le z_{m+1})}+1_{(z_{m+1}\le w\le y_{m+1})})z^p_{m+1}\cr
&\phantom{\equiv \int\int(1_{(y_{m+1}\le w\le z_{m+1})}}\times|D^2_{z_{m+1}}p_t(z,y)|dz^{(m)}\mu_{m+1}(dy_{m+1})\cr
&\le ct^{-2}z^p_{m+1}\liminf_{\delta\to0}\E_{z_{m+1}}\Bigl(\psi(z',N_0)\int q_{\delta}(y_{m+1},\Xm_t)(1_{(y_{m+1}\le w\le z_{m+1})}\cr
&\phantom{\le ct^{-2}z^p_{m+1}\liminf_{\delta\to0}\E_{z_{m+1}}\Bigl(\psi(z',N_0)\int }+1_{(z_{m+1}\le w\le y_{m+1})})\mu_{m+1}(dy_{m+1})\cr
&= ct^{-2}z^p_{m+1}\E_{z_{m+1}}\Bigl(\psi(z',N_0)[1_{(w\le z_{m+1})}1_{(\Xm_t\le w)}\cr
&\phantom{= ct^{-2}z^p_{m+1}\E_{z_{m+1}}\Bigl(\psi(z',N_0)[}+1_{(w\ge z_{m+1})}1_{(\Xm_t\ge w)}]\Bigr)\cr
&\le ct^{-2} z_{m+1}^p \E_{z_{m+1}}(\psi(z',N_0)^2)^{1/2}[1_{(w\le z_{m+1})}\P_{z_{m+1}}(\Xm_t\le w)^{1/2}\cr
&\phantom{\le ct^{-2}\E_{z_{m+1}}(\psi(z',N_0)^2)^{1/2}[}+1_{(w\ge z_{m+1})}\P_{z_{m+1}}(\Xm_t\ge w)^{1/2}].&(7.45)}$$
In the third line we have used the a.s. and, hence weak, convergence of $\Xm_\delta$ to $\Xm_0$ as $\delta\to0$ and the fact that $\Xm_t\neq w$ a.s.

We have
$$\eqalignno{\E(\psi(z',N_0)^2)&\le c\Bigl((1+z')e^{-z'}+(z')^{-1}e^{-z'}\cr
&\phantom{=(1+z')e^{-z'}}+(z')^{-4}[\E(1_{(N_0\ge 2)}(N_0-z')^4)+\E(1_{(N_0\ge 2)}N_0^2)]\Bigr).}$$
An elementary calculation (consider small $z'$ and large $z'$ separately) shows that the term in square brackets is at most $c(z')^2$. Therefore we deduce that
$$\E(\psi(z',N_0)^2)^{1/2}\le c(z')^{-1}.\eqno(7.46)$$
If we set $w'=w/\gamma t$, then this, together with (7.45) and Lemma~7.2 allows us to conclude that
$$\eqalignno{\int& J(z_{m+1})\mu_{m+1}(dz_{m+1})\cr
&\le ct^{-2}\int z_{m+1}^{p-1+b/\gamma}(z')^{-1}(w/z_{m+1})^{b/4\gamma}\exp\Bigl({-(\sqrt{ z_{m+1}}-\sqrt w)^2\over 2\gamma t}\Bigr) dz_{m+1}\cr
&=ct^{p-2+b/\gamma}(w')^{b/4\gamma}\int (z')^{p-2+3b/4\gamma}\exp\Bigl({-(\sqrt{ z'}-\sqrt w')^2\over 2}\Bigr) dz'\cr
&\equiv ct^{p-2+b/\gamma}(w')^{b/4\gamma}K_{p-2+3b/4\gamma}(w').&(7.47)}$$
A simple calculation using the obvious substitution $x=(\sqrt{ z'}-\sqrt w')^2$ shows that for any $\eps>0$ there is a $c_0(\eps)$ such that 
$$K_r(w')\le c_0[1_{(w'\le 1)}+(w')^{r+1/2}1_{(w'>1)}] \hbox{ for all $w'\ge 0$ and }-1+\eps\le r\le \eps^{-1}.$$
Our bounds on $p$ and Convention 3.1 imply that $r=p-2+3b/4\gamma\in[-1+3(4M_0^2)^{-1},M_0^2]$.  Therefore the left-hand side of (7.47) is at most 
$$\eqalignno{c&t^{p-2+b/\gamma}(w')^{b/4\gamma}(1_{(w'\le 1)}+1_{(w'>1)}(w')^{p-3/2+3b/4\gamma})\cr
&\le c(1_{(w\le \gamma t)}t^{p-2+b/\gamma}+1_{(w>\gamma t)}t^{-1/2}w^{p-3/2+b/\gamma}).}$$\qed

\subsec{8. A Remark on the Two-dimensional Case}  

As has already been noted, the proof of Proposition~2.2 (by far the most challenging step) simplifies substantially if $d=2$.  As this is the case required in [DGHSS], we now describe this simplification in a bit more detail. 

Recall the three cases (i)--(iii) for $d=2$ listed following Theorem~1.4.  As noted there, the case $\sE=\emptyset$ is covered by 
Theorem~A of [BP] (with $d=2$) without removing $(0,0)$ from the state space, so we will focus mainly on the other two cases here (but see the last paragraph below). 
In these cases the localization in Theorem~2.1 reduces the problem to the study of the martingale problem for a perturbation of the constant coefficient operator
$$\sA^0=\sum_{i=1}^2 b_i^0D_{x_i}+\gamma_i^0x_2D_{x_i}^2,\eqno(8.1)$$
with resolvent $R_\lam$ and semigroup $P_t$.  Our job is to establish Proposition~2.2 for this resolvent.
 
 For $f\in \sD_0$, we have 
$$\Vert\sA^0R_\lam f\Vert_2=\Vert\lam R_\lam f-f\Vert_2\le 2\Vert f\Vert_2,\eqno(8.2)$$
the latter by Corollary~2.12.  We may therefore remove the term $x_2(R_\lam f)_{22}$ from the summation in (2.6) because $L^2$-boundedness of this term will follow from the other three and (8.2).
Recall the required boundedness of any of the terms was reduced to (2.28) by Cotlar's Lemma.   (Proposition 2.10 was only used to ensure the operator $T_t$ was bounded on $L^2$ so that Cotlar's Lemma may be employed in the proof of Proposition~2.2.  This boundedness, as well as the bound of $ct^{-1}$, is also implied by (2.28) with $s=t$ and the elementary Lemma~3.16. So we only need consider (2.28).)
For the two derivatives involving $x_1$, (2.28) was fairly easily checked in Lemma~4.3 thanks to the ``explicit'' formulae (4.7) and (4.8), and the bound in Lemma~3.3(f).  

It remains only to check (2.28) for $D_{x_2}$.  This was done in Proposition~4.9, using only Lemma 4.6(a) and in fact only used (4.14) for $p=0$ and (4.13) for $p\le 0$.  These proofs in turn were fairly simple consequences of part (c) of the key Lemma~7.1. (Admittedly the proof of (4.13) was omitted, being much simpler than that of (4.15).)  As (7.4) was a trivial consequence of (7.3) (recall (7.26)), we have essentially reduced the
two-dimensional case to the proof of (7.3).  To justify our earlier statements, that this really is much simpler than that of (7.1), we give the proof.  
 At the risk of slightly lengthening the argument we will take this opportunity to explicitly write an integration by parts formula which was implicit (and hidden) in the 
more complicated setting of Lemma~7.1.  Recall that $m=1$ (the proofs below are the same for general $m$), $I_t=\int_0^tX_s^{(2)}\,ds$, and (see (3.26))
$$G_{t,z_1}f(I,X)=\int f(x_1,X)p_{2\gamma^0_1I}(x_1-z_1-b_1^0t)\,dx_1.$$  $N_0=N_0(t) $ is the Poisson variable in Lemma~3.4.  

\proclaim Proposition 8.1 (Integration by Parts Formula).  If $f:\R_+\times\R\to\R$ is bounded and Borel, then 
$$D_{z_2}P_tf(z)=(\gamma t)^{-1}\E_{z_2}\Bigl({(N_0-(z_2/\gamma t))\over z_2/\gamma t}G_{t,z_1}f(I_t, X_t^{(2)})\Bigr)+E_1(t,z,f),\eqno(8.3)$$
where $E_1$ is given by (8.6) below and satisfies
$$\eqalignno{|E_1(t,z,f)|\le \E_{z_2}\Bigl(\int\int\int &|f(x_1,X_t^{(2)})|4I_t^{-1}p_{4\gamma^0_1(u+I_t)}(x_1-z_1-b^0_1t)\cr
&\times1_{(u\le\int_0^t\nu_sds)}dx_1\,du\,d\N_0(\nu)\Bigr).&(8.4)}$$

\proof By (3.25) 
$$\eqalignno{D_{z_2}P_tf(z)&=\E_{z_2}\Bigl(\int [G_{t,z_1}f\Bigl(\int X_s^{(2)}+\nu_sds,X_t^{(2)}\Bigr)- G_{t,z_1}f\Bigl(\int X_s^{(2)},X_t^{(2)}\Bigr)]1_{(\nu_t=0)}\N_0(d\nu)\Bigr)\cr
&\quad+\E_{z_2}\Bigl(\int [G_{t,z_1}f\Bigl(\int X_s^{(2)}+\nu_sds,X_t^{(2)}+\nu_t\Bigr)\cr
&\phantom{quad+\E_{z_2}\Bigl(\int [}- G_{t,z_1}f\Bigl(\int X_s^{(2)}ds,X_t^{(2)}\Bigr)]1_{(\nu_t>0)}\N_0(d\nu)\Bigr)\cr
&\equiv E_1(t,z,f)+E_2(t,z,f).&(8.5)}$$
Now use 
$${\partial p_t\over \partial t}(z)=(z^2t^{-1}-1)(2t)^{-1}p_t(z),$$
and (by some calculus)
$$\Bigl|{\partial p_t\over \partial t}(z)\Bigr|\le {4\over t}p_{2t}(z),$$
together with the Fundamental Theorem of Calculus, to obtain
$$\eqalignno{E_1(t,z,f)=&\E_{z_2}\Bigl(\int\int\int f(x_1,X_t^{(2)})\Bigl[{(x_1-z_1-b^0_1t)^2\over 2\gamma_1^0(u+I_t)}-1\Bigr](2(u+I_t))^{-1}\cr
&\ \times p_{2\gamma_1^0(u+I_t)}(x_1-z_1-b_1^0t)1_{(u\le\int_0^t\nu_sds)}dx_1\,du1_{(\nu_t=0)}d\N_0(\nu)\Bigr),&(8.6)}$$
and
$$\eqalign{|E_1(t,z,f)|\le \E_{z_2}\Bigl(\int\int\int &|f(x_1,X_t^{(2)})|4(u+I_t)^{-1}p_{4\gamma_1^0(u+I_t)}(x_1-z_1-b_1^0t)\cr
&\times1_{(u\le \int_0^t\nu_sds)}dx_1\,du\,d\N_0(\nu).}$$
The latter inequality gives (8.4).

For $E_2$ we use the decomposition in Lemma~3.4 with $\rho=0$.  $S_n$ and $R_n$ are the sum of the first $n$ 
of the $e_i$ and $r_i$, respectively, and we continue to write $p_n(w)=e^{-w}w^n/n!$ and $z'=z_2/\gamma t$.  Then  Lemma~3.4 allows us to write
$$\eqalignno{E_2(t,z,f)&=\N_0(\nu_t>0)\E_{z_{2}}\Bigl(G_{t,z_2}f(I_2(t)+R_{N_0+1},X_0'(t)+S_{N_0+1})\cr
&\phantom{\N_0(\nu_t>0)\E_{z_{2}}\Bigl(G_{t,z_2}f}-G_{t,z_2}f(I_2(t)+R_{N_0},X_0'(t)+S_{N_0})\Bigr)\cr
&=(\gamma t)^{-1}\sum_{n=0}^\infty(p_{n-1}(z')-p_n(z'))\E_{z_2}(G_{t,z_2}f(I_2(t)+R_n,X'_0(t)+S_n))\cr
&= (\gamma t)^{-1}\sum_{n=0}^\infty p_{n}(z')(n-z')(z')^{-1}\E_{z_2}(G_{t,z_2}f(I_2(t)+R_n,X'_0(t)+S_n))\cr
&=(\gamma t)^{-1}\E_{z_2}((N_0-z')(z')^{-1}G_{t,z_2}f(I_t,X^{(2)}_t)).&(8.7)}$$
In the last line we have again used Lemma 3.4 to reconstruct $X^{(2)}$.   (8.7) and (8.5) complete the proof.  \qed

\proclaim Remark 8.2. {\rm Since $|G_{t,z_1}f|\le \Vert f\Vert_\infty$, the above implies the sup
norm bound
$$\eqalign{|&D_{z_2}P_tf(z)|\cr
&\le (\gamma t)^{-1}\E_{z_2}\Bigl({|N_0-z'|\over z'}\Bigr)\Vert f\Vert_\infty +4\Vert f\Vert_\infty \E_{z_2}(I_t^{-1})\int\int_0^t\nu_s\,ds\,d\N_0(\nu)\Bigr)\cr
&\le (\gamma t)^{-1}2\Vert f\Vert_\infty+4t^{-1}\Vert f\Vert_\infty.}$$
We have used Lemma 3.3(f) and (3.17) in the last.  This gives a derivation of (4.4).  More importantly we can use the above to derive an $L^1$ bound which will allow us to take $f=\delta_y$.  Recall that $q_t(x,y)$ is the transition density of $X^{(2)}$ with respect to $y^{b/\gamma-1}dy$. 
}

\proclaim Corollary 8.3.  If $f:\R_+\times \R\to\R$ is bounded and Borel, then 
$$\eqalign {\int &|D_{z_2}P_tf(z)|dz_1\cr
&\le c_{8.3}t^{-1}\E_{z_2}\Bigl(\int|f(z_1,X_t^{(2)})|dz_1\Bigl[{|N_0-z_2/\gamma t|\over z_2/\gamma t}+\Bigl({z_2\over \gamma t}+1\Bigr)^{-1}\Bigr]\Bigr).}$$

\proof Note first that $\int |G_{t,z_1}f(I,X)|dz_1\le\int |f(x_1,X)|dx_1$, and then integrate over $z_1$ in Proposition~7.1 to see that the above integral is at most ($z'=z_2/\gamma t$ as usual)
$$\eqalignno{(\gamma t)^{-1}&\E_{z_2}\Bigl({|N_0-z'|\over z'}\int |f(x_1,X_t^{(2)})|dx_1\Bigr)\cr
&+\E_{z_2}\Bigl(\int\int|f(x_1,X_t^{(2)})|dx_14I_t^{-1}\Bigr)\Bigl(\int\int_0^t\nu_s\,ds\,d\N_0(\nu)\Bigr).&(8.8)}$$
Use Lemma 3.3(f) and (3.17) again to bound the last term by
$$4c_{3.2}(t+z_2)^{-1}\E_{z_2}\Bigl(\int |f(x_1,X_t^{(2)})|dx_1\Bigr).$$
Use this in (8.8) to derive the required bound.\qed

\ni {\bf Proof of (7.3).} Let $f^{y,\delta}(z_1,z_2)=p_\delta(z_1-y_1)q_\delta(y_2,z_2)$ (bounded in $z$ by Lemma~3.3(a)).  Then (3.30) shows that $\lim _{\delta\to 0}D_{z_2}P_t f^{y,\delta}(z)=D_{z_2}p_t(z,y)$.  Apply Fatou's Lemma and Corollary~8.3 to conclude
$$\eqalign{\int|&D_{z_2}p_t(z,y)|dz_1\cr
&\le \liminf_{\delta\to0}c_{7.3}t^{-1}\E_{z_2}\Bigl(\int |f^{y,\delta}(z_1,X_t^{(2)})|dz_1\Bigl[{|N_0-z_2/\gamma t|\over z_2/\gamma t}+\Bigl({z_2\over \gamma t}+1\Bigr)^{-1}\Bigr]\Bigr)\cr
&\le \liminf_{\delta\to 0}c_{7.3}t^{-1}\E_{z_2}\Bigl(q_\delta (y_2,X_t^{(2)})\Bigl[{|N_0-z_2/\gamma t|\over z_2/\gamma t}+\Bigl({z_2\over \gamma t}+1\Bigr)^{-1}\Bigr]\Bigr).}$$
The required result follows. \qed
If we wanted to include the case $\sE=\emptyset$ to make the above ``short proof" self-contained, then we need to consider Proposition~2.2 and hence (4.1) and (4.2) for the case 
$${\sA'}^0=\sum_{i=1}^2b_i^0D_{x_i}+\gamma_i^0x_iD_{x_i}^2.$$ 
The associated semigroup $P_t=\prod_{j=1}^2Q^j_t$ is a product of one-dimensional Feller branching (with immigration) semigroups with transition densities given by (3.3).  As in the the last part of the proof of Proposition~2.14 at the end of Section~4, (4.1) and (4.2) reduce easily to checking (4.1) and (4.2) for each 
one dimensional $Q_t^i$.  In the first part of the proof of Proposition 2.14 (in Section~4) we saw that these easily followed for each differential operator by projecting down the corresponding result for $\sA^0$ (as in (8.1)) to the second coordinate.  This was checked in the ``short" proof above for the the first order operators.
It therefore only remains to check (4.1) and (4.2) for $\tilde D_x=xD^2_{x}$ and $q_t$ in place of $p_t$.  As in the proof of Proposition~4.12, we must verify (4.33), (4.34), (4.24), and (4.25) for this operator and one-dimensional density.  These, however, can be done by direct calculation using the series expansion (3.3)--the arguments are much simpler and involve direct summation by parts with Poisson probabilities and elementary Poisson bounds.

\subsec{References}

\refer{ABBP} S. R. Athreya, M.T. Barlow, R.F. Bass, and E.A. Perkins, Degenerate stochastic differential equations and super-Markov chains.  {\it Prob. Th. Rel. Fields \bf 123} (2002), 484--520.

\refer{Ba} R.F. Bass, {\it Probabilistic Techniques in Analysis}, 
Springer, Berlin 1995.

\refer{BP} R.F. Bass and E.A. Perkins, Degenerate stochastic differential equations with H\"older continuous coefficients and super-Markov chains. 
{\it Trans. Amer. Math. Soc. \bf 355} (2003) 373--405.

\refer{DF} D.A. Dawson and K. Fleischmann, Catalytic and mutually catalytic branching.  In {\it Infinite Dimensional Stochastic analysis}, Ned. Acak. Wet., Vol. 52, R. Neth. Acad. Arts Sci., Amsterdam, 2000, pp. 145--170.  

\refer{DFX} D. A. Dawson, K. Fleischmann, and J. Xiong, Strong uniqueness for cyclically catalytic symbiotic branching diffusions. {\it Statist. Probab. Lett. \bf 73} (2005) 251--257.

\refer{DGHSS} \ \ D.A. Dawson, A. Greven, F. den Hollander, Rongfeng Sun, and J.M. Swart, The renormalization transformation for two-type branching models, to appear {\it Ann. de l'Inst. H. Poincar\'e, Prob. et Stat.}

\refer{DP1} D.A. Dawson and E.A. Perkins, On the uniqueness problem for catalytic branching networks and other singular diffusions. 
{\it Illinois J. Math. \bf 50} (2006) 323--383.

\refer{DP2} D.A. Dawson and E. A. Perkins, Long-time behaviour and coexistence in a mutually catalytic branching model. {\it Ann. Probab. \bf 26} (1998) 1088--1138.

\refer {ES} M. Eigen and P. Schuster, {\it The Hypercycle: a Principle of Natural Self-organization}, Springer, Berlin, 1979.

\refer{Fe} C. Fefferman, Recent progress in classical Fourier analysis. {\it 
Proceedings of the International Con\-gress of Mathematicians},  Vol. 1, pp. 95--118. Montr\'eal, Canadian  Math. Congress, 1975.

\refer{FX} K. Fleischmann and J. Xiong, A cyclically catalytic super-Brownian motion.  {\it Ann. Probab. \bf 29} (2001) 820--861.

\refer{K} S. Kliem, Degenerate stochastic differential equations 
for catalytic branching networks, preprint.
 
\refer{M} L. Mytnik, Uniqueness for a mutually catalytic branching model. {\it Prob. Th. Rel. Fields \bf 112} (1998) 245-253.  

\refer{P} E.A. Perkins, Dawson-Watanabe Superprocesses and Measure-Valued Diffusions, in {\it Lectures on Probability and Statistics, Ecole d'\'Et\'e de Probabilit\'es de Saint-Flour XXIX (1999)}, LNM vol. 1781, Springer-Verlag, Berlin, 2002, pp. 125--324.

\refer{RY} D. Revuz and M. Yor, {\it Continuous Martingales and Brownian Motion}, Berlin, Springer-Verlag, 1991. 

\refer{SV} D. W. Stroock and S. R. S. Varadhan, {\it Multidimensional Diffusion Processes}, Springer-Verlag, Berlin 1979.

\refer{T} A. Torchinsky, {\it Real-variable methods in harmonic analysis.} 
Academic Press,  Orlando, FL, 1986. 

\bye